\documentclass[11pt,twoside]{article}
\usepackage{amsthm}
\usepackage{amsmath}
\usepackage{amssymb}
\usepackage{array}
\usepackage{graphicx}
\usepackage{graphics}
\usepackage{epsfig}
\usepackage{amsfonts,amscd}
\usepackage{exscale}
\usepackage{eucal}
\usepackage{latexsym,amsmath}
\usepackage{indentfirst}
\numberwithin{equation}{section}
\newtheorem{Prop}{\bf Proposition}[section]

\newtheorem{Rem}{\bf Remark}[section]

\theoremstyle{definition} \theoremstyle{plain}

\begin{document}
\def \b{\Box}

\begin{center}
{\Large \bf Programs in $C^{++}$ for matrix computations in min plus algebra}\\[0.1cm]
\end{center}
\begin{center}
{\Large \bf Mihai Ivan and Gheorghe Ivan} \\[0.1cm]
\end{center}

\setcounter{page}{1}

 \pagestyle{myheadings}

{\small {\bf Abstract}. The main purpose of this paper is to
propose six programs in $C^{++}$ for matrix computations and
solving recurrent equations systems with entries in min plus
algebra.}
{}

\section{Introduction}
\smallskip
\indent\indent Idempotent mathematics is based on replacing the
usual arithmetic operations with a new set of basic operations,
that is on replacing numerical fields by idempotent semirings.
Exotic semirings such as the max plus algebra ${\bf R}_{max}$ and
min plus algebra ${\bf R}_{min}$ have been introduced in
connection with various fields: graph theory, Markov decision
processes, discrete event systems theory, see \cite{bcoq, kusa,
mohr, litv, ivan}.

In min plus algebra (resp., max plus algebra), the arithmetic
addition of the conventional algebra is replaced by the point-wise
minimization (resp., maximization), denoted here by the symbol
$\oplus$ and the arithmetic multiplication is replaced by the
point-wise addition, represented by the symbol $\otimes$.

Several results in min plus algebra can be successfully applied to
a number of problems in networking, see \cite{both}.

 The paper is organized as follows. The semiring of
matrices with entries in min plus algebra is presented in Section
2. In Section 3 we give six programs in language $~C^{++}~$  for
matrix computations in
 min plus algebra.

\section{Semirings. Matrices over min plus algebra}
\indent\indent We start this section by recalling of some
necessary backgrounds on semirings for our purposes (see
\cite{bcoq, mohr, olro} and references therein for more details).

Let $S$ be a nonempty set endowed with two binary operations, {\it
addition} (denoted with $\oplus$) and {\it multiplication} (denoted
with $\otimes$). The algebraic structure $(S,\oplus, \otimes,
\varepsilon, e )$ is a \textit{semiring}, if it fulfills the
following conditions:

$(1)~~(S,\oplus, \varepsilon)~$ is a commutative monoid with
$\varepsilon$ as the neutral element for $\oplus;$

$(2)~(S, \otimes,  e )~$ is a  monoid with $\varepsilon$ as the
identity element for $\otimes;$

$(3)~~\otimes $ distributes over $\oplus;$

$(4)~~\varepsilon~$ is an absorbing element for $\otimes$, that is
$~a\otimes \varepsilon=\varepsilon \otimes a= \varepsilon,~\forall
a\in S.$

A semiring where addition is idempotent (that is, $~a\oplus
a=a,~\forall a \in S$) is called an {\it idempotent semiring}. If
$\otimes$ is commutative, we say that $S$ is a {\it commutative
semiring}.

Let $(S, \oplus, \otimes, \varepsilon, e )$ be an (idempotent)
semiring.  For each pair of positive integer $(m,n)$, let $
M_{m\times n}(S)$ be denote the set of $m\times n$ matrices with
entries in $S$. The operations $\oplus$ and $\otimes$ on $S$ induce
corresponding operations on $ M_{m\times n}(S)$ in the obvious way.
Indeed, if
$A=(A_{ij}), B=(B_{ij})\in M_{m\times n}(S)$ then we have:\\[-0.2cm]
\begin{equation}
A\oplus B= ((A\oplus B)_{ij})~~~\hbox{where}~~~ (A\oplus B)_{ij}:=
A_{ij} \oplus B_{ij}.\label{(2.1)}
\end{equation}
If $A=(A_{ij})\in M_{m\times n}(S)$ and $B=(B_{jk})\in M_{n\times p}(S)$ then we have:\\[-0.2cm]
\begin{equation}
A\otimes B= ((A\otimes B)_{ik}),~i =\overline{1,m},
~k=\overline{1,p}~~\hbox{where}~ (A\otimes B)_{ik}:=
\bigoplus\limits_{j=1}^{n} A_{ij}\otimes B_{jk}.\label{(2.2)}
\end{equation}

The product of a matrix $A=(A_{ij}) \in M_{m\times n}(S)$ with a
scalar
$\alpha\in S$ is given by:\\[-0.2cm]
\begin{equation}
\alpha \otimes A= ((\alpha\otimes A)_{ij})
~~~\hbox{where}~~~(\alpha\otimes A)_{ij}:= \alpha  \otimes
A_{ij}.\label{(2.3)}
\end{equation}

 The set $M_{n\times n}(S)$ contains two special matrices with entries in
$S$, namely the zero matrix $O_{\oplus n}$, which  has all its
entries equal to $\varepsilon$, and the identity matrix $I_{\otimes
n}$, which  has the diagonal entries equal to $e$ and the other
entries equal to $\varepsilon$.

It is easy to check that the following proposition holds.
\begin{Prop}
$( M_{n\times n}(S), \oplus, \otimes, O_{\oplus n}, I_{\otimes n}) $
is an idempotent semiring, where the operations $\oplus$ and
$\otimes$ are given in $(2.1)$ and $(2.2)$.\hfill$\Box$
\end{Prop}

\markboth{M. Ivan, Gh. Ivan}{Programs in $C^{++}$ for matrix
computations in min plus algebra}

We call $( M_{n\times n}(S), \oplus, \otimes, O_{\oplus n},
I_{\otimes n}) $ the {\it semiring of $n\times n$ matrices with
entries in $S$}. In particular, if $S:={\bf R}_{min}= ({\bf R}\cup
\{+\infty\},\oplus:=min, \otimes:=+, \varepsilon:=+\infty, e:= 0$ is
called the {\it semiring of $n\times n$ matrices over $ {\bf
R}_{min}$}.

When $S={\bf R}_{min}$, the operations $\oplus$ and
$\otimes$ given in $(2.1)$ and $(2.2)$, becomes:\\[-0.2cm]
\begin{equation}
(A\oplus B)_{ij}:= min\{A_{ij}, B_{ij}\}~~~\hbox{and}~~~ (A\otimes
B)_{ik}:= min_{1\leq k\leq n} \{A_{ij}+ B_{jk}\}.\label{(2.4)}
\end{equation}

The operation $\otimes$ on $ M_{m\times n}({\bf R}_{min})$
given in $(2.3)$ becomes:\\[-0.2cm]
\begin{equation}
\alpha \otimes A= ((\alpha\otimes A)_{ij})
~~~\hbox{where}~~~(\alpha\otimes A)_{ij}:= \alpha +
A_{ij}.\label{(2.5)}
\end{equation}

In conventional algebra we know that for a matrix $A=(A_{ij})\in
M_{n\times n}({\bf R})$, the determinant $\det(A)$ is given
by\\[-0.2cm]
\[
\det(A)=\sum_{\sigma\in P_{n}} sgn(\sigma)
\prod\limits_{i=1}^{n}A_{i\sigma(i)},
\]
where $P_{n}$ is the set of all permutations of the set $\{1, 2,
..., n\}$ and $ sgn(\sigma)$ is the signature of the permutation
$\sigma$.

In min plus algebra (resp., max plus algebra) the determinant has
no direct analogue because of the absence of addition inverses.
The concept of permanent of a matrix partially play the role of
determinant. It is defined similarly of the determinant but with
the $ sgn(\sigma)$ simply omitted, see \cite{olro}.

For $ A=(A_{ij})\in M_{n\times n}(S)$, the {\it permanent} of $A$,
denoted by $perm(A)$ is defined by\\[-0.2cm]
\[
perm(A)= \bigoplus_{\sigma\in P_{n}}
\bigotimes\limits_{i=1}^{n}A_{i\sigma(i)}.\\[-0.2cm]
\]

Let $P_{n}^{e}$ (resp.,$P_{n}^{o}$) be the set of even (resp.,
odd) permutations of the set $\{1, 2, ..., n\}$. The {\it
bideterminant} of $ A=(A_{ij})\in M_{n\times n}(S)$ is the pair $(
\Delta_{1}(A), \Delta_{2}(A))$,  where\\[-0.2cm]
\begin{equation}
\Delta _{1}(A)= \bigoplus_{\sigma\in P_{n}^{e}}
\bigotimes\limits_{i=1}^{n}A_{i\sigma(i)}~~~\hbox{and}~~~\Delta
_{2}(A)= \bigoplus_{\sigma\in P_{n}^{o}}
\bigotimes\limits_{i=1}^{n}A_{i\sigma(i)}.\label{(2.6)}
\end{equation}

The element $ \Delta_{1}(A)$ (resp., $ \Delta_{2}(A)$) is called
{\it plus-determinant} (resp., {\it minus-determinant}) of $A$ and
is denoted with  $ \Delta_{1}(A):=\det_{\oplus}^{+}(A) $ (resp.,
$\Delta_{2}(A):=\det_{\oplus}^{-}(A)$).

Note that $perm(A)=\Delta_{1}(A)\oplus \Delta_{2}(A).$

When $S={\bf R}_{min}$, the bideterminant of $ A $  given in $(2.6)$, becomes:\\[-0.2cm]
\begin{equation}
\Delta _{1}(A)= \min_{\sigma\in P_{n}^{e}}
\sum\limits_{i=1}^{n}A_{i\sigma(i)}~~~\hbox{and}~~~\Delta _{2}(A)=
\min_{\sigma\in P_{n}^{o}}
\sum\limits_{i=1}^{n}A_{i\sigma(i)}.\label{(2.7)}
\end{equation}

\section{Programs in $C^{++}$}

In this section we give programs written in the language $C^{++}$
for the basic operations with matrices over ${\bf
R}_{min}$ and for solving a recurrent linear system, and namely:\\[0.2cm]
${\bf 1.}~~~$ the sum of two matrices $A, B\in M_{m\times n}({\bf
R}_{min})$;\\[0.1cm]
${\bf 2.}~~~$ the product of two matrices $A\in M_{m\times n}({\bf
R}_{min})$ and $B\in M_{n\times p}({\bf R}_{min})$;\\[0.1cm]
${\bf 3.}~~~$ the product of a a scalar $\alpha\in {\bf R}_{min}$
with a matrix $A\in M_{m\times n}({\bf R}_{min})$;\\[0.1cm]
${\bf 4.}~~~$ the power  of a matrix $A\in M_{n\times n}({\bf R}_{min})$;\\[0.1cm]
${\bf 5.}~~~~ \Delta_{1}(A)$ and $\Delta_{1}(A)$ for $A\in M_{n\times n}({\bf R}_{min})$ when $n=2$ and $n=3$;\\[0.1cm]
${\bf 6.}~~~$ the solving a linear system of the form:\\[-0.2cm]
\[
X(k+1)=A\otimes X(k),~~~k\geq 0,
\]
where $~A\in M_{n\times n}({\bf R}_{min})$ and $X(k)\in M(n,1;{\bf
R}_{min}).$
\begin{Rem}
{\rm In \cite{ivmg} are given five programs in $C^{++}$ for matrix
computations in ${\bf R}_{max}$}.\hfill$\Box$
\end{Rem}

$\bullet~~~~~$ {\it We first give the Form1 Designer of program}.\\[0.2cm]
name space Operations$_{-}$with$_{-}$matrices$_{-}$in$_{-}$min$_{-}$plus$_{-}$algebra\\[0.1cm]
$\{$\\[0.1cm]
    \hspace*{0.5cm}partial class Form1\\
    \hspace*{0.5cm}$\{$\\
        \hspace*{1cm}/// $<$summary$>$\\
        \hspace*{1cm}/// Required designer variable.\\
        \hspace*{1cm}/// $<$/summary$>$\\
        \hspace*{1cm}private System.ComponentModel.IContainer components = null;\\[0.2cm]
        \hspace*{1cm}/// $<$summary$>$\\
        \hspace*{1cm}/// Clean up any resources being used.\\
        \hspace*{1cm}/// $<$/summary$>$\\
        \hspace*{1cm}/// $<$param name="disposing"$>$true if managed resources should be disposed; otherwise, false.$<$/param$>$\\
        \hspace*{1cm}protected override void Dispose(bool disposing)\\
        \hspace*{1cm}$\{$\\
          \hspace*{1.5cm}  if (disposing $\&\&$ (components ! = null))\\
            \hspace*{1.5cm}$\{$\\
              \hspace*{1.7cm}  components.Dispose();\\
            \hspace*{1.5cm}$\}$\\
            \hspace*{1.5cm}base.Dispose(disposing);\\
        \hspace*{1cm}$\}$\\[0.2cm]
        \hspace*{1cm}$\#$region Windows Form Designer generated code\\[0.2cm]
        \hspace*{1cm}/// $<$summary$>$\\
        \hspace*{1cm}/// Required method for Designer support-do not modify\\
        \hspace*{1cm}/// the contents of this method with the code editor.\\
        \hspace*{1cm}/// $<$/summary$>$\\
        \hspace*{1cm}private void InitializeComponent()\\
        \hspace*{1cm}$\{$\\
        \hspace*{1.5cm}System.ComponentModel.ComponentResourceManager resources = new System.\\
 ComponentModel.ComponentResourceManager(typeof(Form1));\\
            \hspace*{1.5cm}this.tabControl1 = new System.Windows.Forms.TabControl();\\
            \hspace*{1.5cm}this.tabAddition = new System.Windows.Forms.TabPage();\\
            \hspace*{1.5cm}this.btReset = new System.Windows.Forms.Button();\\
            \hspace*{1.5cm}this.btComputationSum = new System.Windows.Forms.Button();\\
            \hspace*{1.5cm}this.btGenerating = new
            System.Windows.Forms.Button();\\
            \hspace*{1.5cm}this.label6 = new System.Windows.Forms.Label();\\
            \hspace*{1.5cm}this.pictureBox1 = new System.Windows.Forms.PictureBox();\\
            \hspace*{1.5cm}this.panel1 = new System.Windows.Forms.Panel();\\
            \hspace*{1.5cm}this.pictureBox2 = new System.Windows.Forms.PictureBox();\\
            \hspace*{1.5cm}this.dataGridResAddition = new System.Windows.Forms.DataGridView();\\
            \hspace*{1.5cm}this.panel2 = new System.Windows.Forms.Panel();\\
            \hspace*{1.5cm}this.label4 = new System.Windows.Forms.Label();\\
            \hspace*{1.5cm}this.label1 = new System.Windows.Forms.Label();\\
            \hspace*{1.5cm}this.dataGridB = new System.Windows.Forms.DataGridView();\\
            \hspace*{1.5cm}this.dataGridA = new System.Windows.Forms.DataGridView();\\
            \hspace*{1.5cm}this.textColumn = new System.Windows.Forms.TextBox();\\
            \hspace*{1.5cm}this.textLine = new System.Windows.Forms.TextBox();\\
            \hspace*{1.5cm}this.label3 = new System.Windows.Forms.Label();\\
            \hspace*{1.5cm}this.label2 = new System.Windows.Forms.Label();\\
            \hspace*{1.5cm}this.tabMultiplication = new System.Windows.Forms.TabPage();\\
            \hspace*{1.5cm}this.pictureBox3 = new System.Windows.Forms.PictureBox();\\
            \hspace*{1.5cm}this.btResetMultiplication = new System.Windows.Forms.Button();\\
            \hspace*{1.5cm}this.btComputingProduct = new System.Windows.Forms.Button();\\
            \hspace*{1.5cm}this.btGenerating2 = new System.Windows.Forms.Button();\\
            \hspace*{1.5cm}this.label5 = new System.Windows.Forms.Label();\\
            \hspace*{1.5cm}this.label13 = new System.Windows.Forms.Label();\\
            \hspace*{1.5cm}this.textcolumnB2 = new System.Windows.Forms.TextBox();\\
            \hspace*{1.5cm}this.textLineB2 = new System.Windows.Forms.TextBox();\\
            \hspace*{1.5cm}this.label11 = new System.Windows.Forms.Label();\\
            \hspace*{1.5cm}this.label12 = new System.Windows.Forms.Label();\\
            \hspace*{1.5cm}this.panel3 = new System.Windows.Forms.Panel();\\
            \hspace*{1.5cm}this.pictureBox4 = new System.Windows.Forms.PictureBox();\\
            \hspace*{1.5cm}this.dataGridProduct = new System.Windows.Forms.DataGridView();\\
            \hspace*{1.5cm}this.panel4 = new System.Windows.Forms.Panel();\\
            \hspace*{1.5cm}this.label7 = new System.Windows.Forms.Label();\\
            \hspace*{1.5cm}this.label8 = new System.Windows.Forms.Label();\\
            \hspace*{1.5cm}this.dataGridB2 = new System.Windows.Forms.DataGridView();\\
            \hspace*{1.5cm}this.dataGridA2 = new System.Windows.Forms.DataGridView();\\
            \hspace*{1.5cm}this.textColumnA2 = new System.Windows.Forms.TextBox();\\
            \hspace*{1.5cm}this.textLineA2 = new System.Windows.Forms.TextBox();\\
            \hspace*{1.5cm}this.label9 = new System.Windows.Forms.Label();\\
            \hspace*{1.5cm}this.label10 = new System.Windows.Forms.Label();\\
            \hspace*{1.5cm}this.tabMultiplication$_{-}$scalar = new System.Windows.Forms.TabPage();\\
            \hspace*{1.5cm}this.pictureBox6 = new System.Windows.Forms.PictureBox();\\
            \hspace*{1.5cm}this.btReset$_{-}$Product$_{-}$Scalar = new System.Windows.Forms.Button();\\
            \hspace*{1.5cm}this.btProductScalar = new System.Windows.Forms.Button();\\
            \hspace*{1.5cm}this.butonsGenerating3 = new System.Windows.Forms.Button();\\
            \hspace*{1.5cm}this.label28 = new System.Windows.Forms.Label();\\
            \hspace*{1.5cm}this.textScalar = new System.Windows.Forms.TextBox();\\
            \hspace*{1.5cm}this.label16 = new System.Windows.Forms.Label();\\
            \hspace*{1.5cm}this.panel6 = new System.Windows.Forms.Panel();\\
            \hspace*{1.5cm}this.pictureBox5 = new System.Windows.Forms.PictureBox();\\
            \hspace*{1.5cm}this.dataGridProductScalar = new System.Windows.Forms.DataGridView();\\
            \hspace*{1.5cm}this.panel5 = new System.Windows.Forms.Panel();\\
            \hspace*{1.5cm}this.label17 = new System.Windows.Forms.Label();\\
            \hspace*{1.5cm}this.dataGridA3 = new System.Windows.Forms.DataGridView();\\
            \hspace*{1.5cm}this.textcolumnA3 = new System.Windows.Forms.TextBox();\\
            \hspace*{1.5cm}this.textLineA3 = new System.Windows.Forms.TextBox();\\
            \hspace*{1.5cm}this.label14 = new System.Windows.Forms.Label();\\
            \hspace*{1.5cm}this.label15 = new System.Windows.Forms.Label();\\
            \hspace*{1.5cm}this.tabLifting$_{-}$at$_{-}$power = new System.Windows.Forms.TabPage();\\
            \hspace*{1.5cm}this.pictureBox9 = new System.Windows.Forms.PictureBox();\\
            \hspace*{1.5cm}this.btReset$_{-}$lifting$_{-}$at$_{-}$power = new System.Windows.Forms.Button();\\
            \hspace*{1.5cm}this.btComputPower = new System.Windows.Forms.Button();\\
            \hspace*{1.5cm}this.btGeneratingA4 = new System.Windows.Forms.Button();\\
            \hspace*{1.5cm}this.label29 = new System.Windows.Forms.Label();\\
            \hspace*{1.5cm}this.label21 = new System.Windows.Forms.Label();\\
            \hspace*{1.5cm}this.panel7 = new System.Windows.Forms.Panel();\\
            \hspace*{1.5cm}this.label27 = new System.Windows.Forms.Label();\\
            \hspace*{1.5cm}this.label25 = new System.Windows.Forms.Label();\\
            \hspace*{1.5cm}this.label26 = new System.Windows.Forms.Label();\\
            \hspace*{1.5cm}this.labelPower = new System.Windows.Forms.Label();\\
            \hspace*{1.5cm}this.label24 = new System.Windows.Forms.Label();\\
            \hspace*{1.5cm}this.dataGridMatrix$_{-}$at$_{-}$power$_{-}$n = new System.Windows.Forms.DataGridView();\\
            \hspace*{1.5cm}this.panel8 = new System.Windows.Forms.Panel();\\
            \hspace*{1.5cm}this.label20 = new System.Windows.Forms.Label();\\
            \hspace*{1.5cm}this.dataGridA4 = new System.Windows.Forms.DataGridView();\\
            \hspace*{1.5cm}this.textPower = new System.Windows.Forms.TextBox();\\
            \hspace*{1.5cm}this.textlineA4 = new System.Windows.Forms.TextBox();\\
            \hspace*{1.5cm}this.label18 = new System.Windows.Forms.Label();\\
            \hspace*{1.5cm}this.label19 = new System.Windows.Forms.Label();\\
            \hspace*{1.5cm}this.tabComputDet = new System.Windows.Forms.TabPage();\\
            \hspace*{1.5cm}this.pictureBox10 = new System.Windows.Forms.PictureBox();\\
            \hspace*{1.5cm}this.bt$_{-}$Reset$_{-}$values$_{-}$det = new System.Windows.Forms.Button();\\
            \hspace*{1.5cm}this.btComputingDeterminant = new System.Windows.Forms.Button();\\
            \hspace*{1.5cm}this.btGeneratingA5 = new System.Windows.Forms.Button();\\
            \hspace*{1.5cm}this.label30 = new System.Windows.Forms.Label();\\
            \hspace*{1.5cm}this.pictureBox8 = new System.Windows.Forms.PictureBox();\\
            \hspace*{1.5cm}this.pictureBox7 = new System.Windows.Forms.PictureBox();\\
            \hspace*{1.5cm}this.textDet$_{-}$minus = new System.Windows.Forms.TextBox();\\
            \hspace*{1.5cm}this.textDet$_{-}$plus = new System.Windows.Forms.TextBox();\\
            \hspace*{1.5cm}this.panel10 = new System.Windows.Forms.Panel();\\
            \hspace*{1.5cm}this.label23 = new System.Windows.Forms.Label();\\
            \hspace*{1.5cm}this.dataGridA5 = new System.Windows.Forms.DataGridView();\\
            \hspace*{1.5cm}this.textlineA5 = new System.Windows.Forms.TextBox();\\
            \hspace*{1.5cm}this.label22 = new System.Windows.Forms.Label();\\
            \hspace*{1.5cm}this.tabPage1 = new System.Windows.Forms.TabPage();\\
            \hspace*{1.5cm}this.panel12 = new System.Windows.Forms.Panel();\\
            \hspace*{1.5cm}this.label41 = new System.Windows.Forms.Label();\\
            \hspace*{1.5cm}this.dataGridX0 = new System.Windows.Forms.DataGridView();\\
            \hspace*{1.5cm}this.panel9 = new System.Windows.Forms.Panel();\\
            \hspace*{1.5cm}this.label$_{-}$k = new System.Windows.Forms.Label();\\
            \hspace*{1.5cm}this.label32 = new System.Windows.Forms.Label();\\
            \hspace*{1.5cm}this.label33 = new System.Windows.Forms.Label();\\
            \hspace*{1.5cm}this.label37 = new System.Windows.Forms.Label();\\
            \hspace*{1.5cm}this.dataGridSolutionXk = new System.Windows.Forms.DataGridView();\\
            \hspace*{1.5cm}this.panel11 = new System.Windows.Forms.Panel();\\
            \hspace*{1.5cm}this.label38 = new System.Windows.Forms.Label();\\
            \hspace*{1.5cm}this.dataGridA6 = new System.Windows.Forms.DataGridView();\\
            \hspace*{1.5cm}this.textk = new System.Windows.Forms.TextBox();\\
            \hspace*{1.5cm}this.textlineA6 = new System.Windows.Forms.TextBox();\\
            \hspace*{1.5cm}this.label39 = new System.Windows.Forms.Label();\\
            \hspace*{1.5cm}this.label40 = new System.Windows.Forms.Label();\\
            \hspace*{1.5cm}this.btResetSystem = new System.Windows.Forms.Button();\\
            \hspace*{1.5cm}this.btComputSystem = new System.Windows.Forms.Button();\\
            \hspace*{1.5cm}this.btGeneratingMatrices = new System.Windows.Forms.Button();\\
            \hspace*{1.5cm}this.tabControl1.SuspendLayout();\\
            \hspace*{1.5cm}this.tabAddition.SuspendLayout();\\[0.2cm]
 ((System.ComponentModel.ISupportInitialize)(this.pictureBox1)).BeginInit();\\
            \hspace*{1.5cm}this.panel1.SuspendLayout();\\[0.2cm]
 ((System.ComponentModel.ISupportInitialize)(this.pictureBox2)).BeginInit();\\[0.2cm]
 ((System.ComponentModel.ISupportInitialize)(this.dataGridRezAddition)).BeginInit();\\
            \hspace*{1.5cm}this.panel2.SuspendLayout();\\[0.2cm]
 ((System.ComponentModel.ISupportInitialize)(this.dataGridB)).BeginInit();\\[0.2cm]
 ((System.ComponentModel.ISupportInitialize)(this.dataGridA)).BeginInit();\\
            \hspace*{1.5cm}this.tabMultiplication.SuspendLayout();\\[0.2cm]
 ((System.ComponentModel.ISupportInitialize)(this.pictureBox3)).BeginInit();\\
            \hspace*{1.5cm}this.panel3.SuspendLayout();\\[0.2cm]
 ((System.ComponentModel.ISupportInitialize)(this.pictureBox4)).BeginInit();\\[0.2cm]
 ((System.ComponentModel.ISupportInitialize)(this.dataGridProduct)).BeginInit();\\
            \hspace*{1.5cm}this.panel4.SuspendLayout();\\[0.2cm]
 ((System.ComponentModel.ISupportInitialize)(this.dataGridB2)).BeginInit();\\[0.2cm]
 ((System.ComponentModel.ISupportInitialize)(this.dataGridA2)).BeginInit();\\
            \hspace*{1.5cm}this.tabMultiplication$_{-}$scalar.SuspendLayout();\\[0.2cm]
 ((System.ComponentModel.ISupportInitialize)(this.pictureBox6)).BeginInit();\\
            \hspace*{1.5cm}this.panel6.SuspendLayout();\\[0.2cm]
 ((System.ComponentModel.ISupportInitialize)(this.pictureBox5)).BeginInit();\\[0.2cm]
 ((System.ComponentModel.ISupportInitialize)(this.dataGridProductScalar)).BeginInit();\\
            \hspace*{1.5cm}this.panel5.SuspendLayout();\\[0.2cm]
 ((System.ComponentModel.ISupportInitialize)(this.dataGridA3)).BeginInit();\\
            \hspace*{1.5cm}this.tabLifting$_{-}$at$_{-}$power.SuspendLayout();\\[0.2cm]
 ((System.ComponentModel.ISupportInitialize)(this.pictureBox9)).BeginInit();\\
            \hspace*{1.5cm}this.panel7.SuspendLayout();\\[0.2cm]
 ((System.ComponentModel.ISupportInitialize)(this.dataGridMatrix$_{-}$at$_{-}$power$_{-}$n)).BeginInit();\\
            \hspace*{1.5cm}this.panel8.SuspendLayout();\\[0.2cm]
 ((System.ComponentModel.ISupportInitialize)(this.dataGridA4)).BeginInit();\\
            \hspace*{1.5cm}this.tabComputDet.SuspendLayout();\\[0.2cm]
 ((System.ComponentModel.ISupportInitialize)(this.pictureBox10)).BeginInit();\\[0.2cm]
 ((System.ComponentModel.ISupportInitialize)(this.pictureBox8)).BeginInit();\\[0.2cm]
 ((System.ComponentModel.ISupportInitialize)(this.pictureBox7)).BeginInit();\\
            \hspace*{1.5cm}this.panel10.SuspendLayout();\\[0.2cm]
 ((System.ComponentModel.ISupportInitialize)(this.dataGridA5)).BeginInit();\\
            \hspace*{1.5cm}this.tabPage1.SuspendLayout();\\
            \hspace*{1.5cm}this.panel12.SuspendLayout();\\[0.2cm]
 ((System.ComponentModel.ISupportInitialize)(this.dataGridX0)).BeginInit();\\
            \hspace*{1.5cm}this.panel9.SuspendLayout();\\[0.2cm]
 ((System.ComponentModel.ISupportInitialize)(this.dataGridSolutionXk)).BeginInit();\\
            \hspace*{1.5cm}this.panel11.SuspendLayout();\\[0.2cm]
 ((System.ComponentModel.ISupportInitialize)(this.dataGridA6)).BeginInit();\\
            \hspace*{1.5cm}this.SuspendLayout();\\
            \hspace*{1.5cm}//\\
            \hspace*{1.5cm}// tabControl1\\
            \hspace*{1.5cm}//\\
            \hspace*{1.5cm}this.tabControl1.Controls.Add(this.tabAddition);\\
            \hspace*{1.5cm}this.tabControl1.Controls.Add(this.tabIMultiplication);\\
            \hspace*{1.5cm}this.tabControl1.Controls.Add(this.tabMultiplication$_{-}$scalar);\\
            \hspace*{1.5cm}this.tabControl1.Controls.Add(this.tabLifting$_{-}$at$_{-}$power);\\
            \hspace*{1.5cm}this.tabControl1.Controls.Add(this.tabComputDet);\\
            \hspace*{1.5cm}this.tabControl1.Controls.Add(this.tabPage1);\\
            \hspace*{1.5cm}this.tabControl1.Location = new System.Drawing.Point(12, 5);\\
            \hspace*{1.5cm}this.tabControl1.Name = "tabControl1";\\
            \hspace*{1.5cm}this.tabControl1.SelectedIndex = 0;\\
            \hspace*{1.5cm}this.tabControl1.Size = new System.Drawing.Size(1000, 696);\\
            \hspace*{1.5cm}this.tabControl1.TabIndex = 0;\\
            \hspace*{1.5cm}//\\
            \hspace*{1.5cm}// tabAddition\\
            \hspace*{1.5cm}//\\
            \hspace*{1.5cm}this.tabAddition.BackColor = System.Drawing.Color.FromArgb(((int)(((byte)(192)))),\\
  ((int)(((byte)(255)))), ((int)(((byte)(192)))));\\
            \hspace*{1.5cm}this.tabAddition.Controls.Add(this.btReset);\\
            \hspace*{1.5cm}this.tabAddition.Controls.Add(this.btComputingSum);\\
            \hspace*{1.5cm}this.tabAddition.Controls.Add(this.btGenerating);\\
            \hspace*{1.5cm}this.tabAddition.Controls.Add(this.label6);\\
            \hspace*{1.5cm}this.tabAddition.Controls.Add(this.pictureBox1);\\
            \hspace*{1.5cm}this.tabAddition.Controls.Add(this.panel1);\\
            \hspace*{1.5cm}this.tabAddition.Controls.Add(this.panel2);\\
            \hspace*{1.5cm}this.tabAddition.Controls.Add(this.textColumn);\\
            \hspace*{1.5cm}this.tabAddition.Controls.Add(this.textLine);\\
            \hspace*{1.5cm}this.tabAddition.Controls.Add(this.label3);\\
            \hspace*{1.5cm}this.tabAddition.Controls.Add(this.label2);\\
            \hspace*{1.5cm}this.tabAddition.Location = new System.Drawing.Point(4, 25);\\
            \hspace*{1.5cm}this.tabAddition.Name = "tabAddition";\\
            \hspace*{1.5cm}this.tabAddition.Padding = new System.Windows.Forms.Padding(3);\\
            \hspace*{1.5cm}this.tabAddition.Size = new System.Drawing.Size(992, 667);\\
            \hspace*{1.5cm}this.tabAddition.TabIndex = 0;\\
            \hspace*{1.5cm}this.tabAddition.Text = "Addition";\\
            \hspace*{1.5cm}//\\
            \hspace*{1.5cm}// btReset\\
            \hspace*{1.5cm}//\\
            \hspace*{1.5cm}this.btReset.BackColor = System.Drawing.Color.Pink;\\
            \hspace*{1.5cm}this.btReset.ForeColor = System.Drawing.Color.Black;\\
            \hspace*{1.5cm}this.btReset.Location = new System.Drawing.Point(546, 296);\\
            \hspace*{1.5cm}this.btReset.Name = "btReset";\\
            \hspace*{1.5cm}this.btReset.Size = new System.Drawing.Size(125, 65);\\
            \hspace*{1.5cm}this.btReset.TabIndex = 38;\\
            \hspace*{1.5cm}this.btReset.Text = "Reset values";\\
            \hspace*{1.5cm}this.btReset.UseVisualStyleBackColor = false;\\
            \hspace*{1.5cm}this.btReset.Click + = new System.EventHandler(this.btReset$_{-}$Click);\\
            \hspace*{1.5cm}//\\
            \hspace*{1.5cm}// btComputingSum\\
            \hspace*{1.5cm}//\\
            \hspace*{1.5cm}this.btComputingSum.BackColor = System.Drawing.Color.PowderBlue;\\
            \hspace*{1.5cm}this.btComputingSum.ForeColor = System.Drawing.Color.Black;\\
            \hspace*{1.5cm}this.btComputingSum.Location = new System.Drawing.Point(546, 178);\\
            \hspace*{1.5cm}this.btComputingSum.Name = "btComputingSum";\\
            \hspace*{1.5cm}this.btComputingSum.Size = new System.Drawing.Size(125, 65);\\
            \hspace*{1.5cm}this.btComputingSum.TabIndex = 37;\\
            \hspace*{1.5cm}this.btComputingSum.Text = "Computing Matrix";\\
            \hspace*{1.5cm}this.btComputingSum.UseVisualStyleBackColor = false;\\
            \hspace*{1.5cm}this.btComputingSum.Click + = new System.EventHandler(this.btComputingSum$_{-}$Click);
            \hspace*{1.5cm}//\\
            \hspace*{1.5cm}// btGenerating\\
            \hspace*{1.5cm}//\\
            \hspace*{1.5cm}this.btGenerating.BackColor = System.Drawing.Color.PowderBlue;\\
            \hspace*{1.5cm}this.btGenerating.ForeColor = System.Drawing.Color.Black;\\
            \hspace*{1.5cm}this.btGenerating.Location = new System.Drawing.Point(546, 24);\\
            \hspace*{1.5cm}this.btGenerating.Name = "btGenerating";\\
            \hspace*{1.5cm}this.btGenerating.Size = new System.Drawing.Size(125, 65);\\
            \hspace*{1.5cm}this.btGenerating.TabIndex = 36;\\
            \hspace*{1.5cm}this.btGenerating.Text = "Generating matrix";\\
            \hspace*{1.5cm}this.btGenerating.UseVisualStyleBackColor = false;\\
            \hspace*{1.5cm}this.btGenerating.Click + = new System.EventHandler(this.btGenerating$_{-}$Click);\\
            \hspace*{1.5cm}//\\
            \hspace*{1.5cm}// label6\\
            \hspace*{1.5cm}//\\
            \hspace*{1.5cm}this.label6.AutoSize = true;\\
            \hspace*{1.5cm}this.label6.Font = new System.Drawing.Font("Microsoft Sans Serif", 9.75F,\\
   System.Drawing.FontStyle.Regular, System.Drawing.GraphicsUnit.Point,((byte)(0)));\\
            \hspace*{1.5cm}this.label6.ForeColor = System.Drawing.Color.Red;\\
            \hspace*{1.5cm}this.label6.Location = new System.Drawing.Point(25, 132);\\
            \hspace*{1.5cm}this.label6.Name = "label6";\\
            \hspace*{1.5cm}this.label6.Size = new System.Drawing.Size(39, 16);\\
            \hspace*{1.5cm}this.label6.TabIndex = 35;\\
            \hspace*{1.5cm}this.label6.Text = "Rem:";\\
            \hspace*{1.5cm}//\\
            \hspace*{1.5cm}// pictureBox1\\
            \hspace*{1.5cm}//\\
            \hspace*{1.5cm}this.pictureBox1.Image =\\
            ((System.Drawing.Image)(resources.GetObject("pictureBox1.Image")));\\
            \hspace*{1.5cm}this.pictureBox1.Location = new System.Drawing.Point(74, 128);\\
            \hspace*{1.5cm}this.pictureBox1.Name = "pictureBox1";\\
            \hspace*{1.5cm}this.pictureBox1.Size = new System.Drawing.Size(56, 27);\\
            \hspace*{1.5cm}this.pictureBox1.TabIndex = 34;\\
            \hspace*{1.5cm}this.pictureBox1.TabStop = false;\\
            \hspace*{1.5cm}//\\
            \hspace*{1.5cm}// panel1\\
            \hspace*{1.5cm}//\\
            \hspace*{1.5cm}this.panel1.Controls.Add(this.pictureBox2);\\
            \hspace*{1.5cm}this.panel1.Controls.Add(this.dataGridResAddition);\\
            \hspace*{1.5cm}this.panel1.Location = new System.Drawing.Point(698, 178);\\
            \hspace*{1.5cm}this.panel1.Name = "panel1";\\
            \hspace*{1.5cm}this.panel1.Size = new System.Drawing.Size(257, 297);\\
            \hspace*{1.5cm}this.panel1.TabIndex = 11;\\
            \hspace*{1.5cm}//\\
            \hspace*{1.5cm}// pictureBox2\\
            \hspace*{1.5cm}//\\
            \hspace*{1.5cm}this.pictureBox2.BackColor = System.Drawing.Color.Red;\\
            \hspace*{1.5cm}this.pictureBox2.Image = global::Operations$_{-}$with$_{-}$matrices.Properties.Resources.a$_{-}$b;\\
            \hspace*{1.5cm}this.pictureBox2.Location = new System.Drawing.Point(17, 10);\\
            \hspace*{1.5cm}this.pictureBox2.Name = "pictureBox2";\\
            \hspace*{1.5cm}this.pictureBox2.Size = new System.Drawing.Size(115, 33);\\
            \hspace*{1.5cm}this.pictureBox2.TabIndex = 13;\\
            \hspace*{1.5cm}this.pictureBox2.TabStop = false;\\
            \hspace*{1.5cm}//\\
            \hspace*{1.5cm}// dataGridResAddition\\
            \hspace*{1.5cm}//\\
            \hspace*{1.5cm}this.dataGridResAddition.BackgroundColor = System.Drawing.Color.WhiteSmoke;\\
            \hspace*{1.5cm}this.dataGridResAddition.ColumnHeadersHeightSizeMode = System.Windows.Forms.\\
  DataGridViewColumnHeadersHeightSizeMode.AutoSize;\\
            \hspace*{1.5cm}this.dataGridResAddition.Location = new System.Drawing.Point(17, 49);\\
            \hspace*{1.5cm}this.dataGridResAddition.Name = "dataGridResAddition";\\
            \hspace*{1.5cm}this.dataGridResAddition.Size = new System.Drawing.Size(223, 227);\\
            \hspace*{1.5cm}this.dataGridResAddition.TabIndex = 8\\
            \hspace*{1.5cm}//\\
            \hspace*{1.5cm}// panel2\\
            \hspace*{1.5cm}//\\
            \hspace*{1.5cm}this.panel2.Controls.Add(this.label4);\\
            \hspace*{1.5cm}this.panel2.Controls.Add(this.label1);\\
            \hspace*{1.5cm}this.panel2.Controls.Add(this.dataGridB);\\
            \hspace*{1.5cm}this.panel2.Controls.Add(this.dataGridA);\\
            \hspace*{1.5cm}this.panel2.Location = new System.Drawing.Point(16, 178);\\
            \hspace*{1.5cm}this.panel2.Name = "panel2";\\
            \hspace*{1.5cm}this.panel2.Size = new System.Drawing.Size(497, 297);\\
            \hspace*{1.5cm}this.panel2.TabIndex = 8;\\
            \hspace*{1.5cm}//\\
            \hspace*{1.5cm}// label4\\
            \hspace*{1.5cm}//\\
            \hspace*{1.5cm}this.label4.AutoSize = true;\\
            \hspace*{1.5cm}this.label4.ForeColor = System.Drawing.Color.MidnightBlue;\\
            \hspace*{1.5cm}this.label4.Location = new System.Drawing.Point(266, 16);\\
            \hspace*{1.5cm}this.label4.Name = "label4";\\
            \hspace*{1.5cm}this.label4.Size = new System.Drawing.Size(72, 16);\\
            \hspace*{1.5cm}this.label4.TabIndex = 11;\\
            \hspace*{1.5cm}this.label4.Text = "Matrix B";\\
            \hspace*{1.5cm}//\\
            \hspace*{1.5cm}// label1\\
            \hspace*{1.5cm}//\\
            \hspace*{1.5cm}this.label1.AutoSize = true;\\
            \hspace*{1.5cm}this.label1.ForeColor = System.Drawing.Color.MidnightBlue;\\
            \hspace*{1.5cm}this.label1.Location = new System.Drawing.Point(27, 16);\\
            \hspace*{1.5cm}this.label1.Name = "label1";\\
            \hspace*{1.5cm}this.label1.Size = new System.Drawing.Size(72, 16);\\
            \hspace*{1.5cm}this.label1.TabIndex = 10;\\
            \hspace*{1.5cm}this.label1.Text = "Matrix A";\\
            \hspace*{1.5cm}//\\
            \hspace*{1.5cm}// dataGridB\\
            \hspace*{1.5cm}//\\
            \hspace*{1.5cm}this.dataGridB.BackgroundColor = System.Drawing.Color.WhiteSmoke;\\
            \hspace*{1.5cm}this.dataGridB.ColumnHeadersHeightSizeMode = System.Windows.Forms.\\
   DataGridViewColumnHeadersHeightSizeMode.AutoSize;\\
            \hspace*{1.5cm}this.dataGridB.Location = new System.Drawing.Point(257, 49);\\
            \hspace*{1.5cm}this.dataGridB.Name = "dataGridB";\\
            \hspace*{1.5cm}this.dataGridB.Size = new System.Drawing.Size(223, 227);\\
            \hspace*{1.5cm}this.dataGridB.TabIndex = 6;\\
            \hspace*{1.5cm}//\\
            \hspace*{1.5cm}// dataGridA\\
            \hspace*{1.5cm}//\\
            \hspace*{1.5cm}this.dataGridA.BackgroundColor = System.Drawing.Color.WhiteSmoke;\\
            \hspace*{1.5cm}this.dataGridA.ColumnHeadersHeightSizeMode = System.Windows.Forms.\\
   DataGridViewColumnHeadersHeightSizeMode.AutoSize;\\
            \hspace*{1.5cm}this.dataGridA.Location = new System.Drawing.Point(12, 49);\\
            \hspace*{1.5cm}this.dataGridA.Name = "dataGridA";\\
            \hspace*{1.5cm}this.dataGridA.Size = new System.Drawing.Size(223, 227);\\
            \hspace*{1.5cm}this.dataGridA.TabIndex = 0;\\
            \hspace*{1.5cm}//\\
            \hspace*{1.5cm}// textColumn\\
            \hspace*{1.5cm}//\\
            \hspace*{1.5cm}this.textColumn.Location = new System.Drawing.Point(253, 70);\\
            \hspace*{1.5cm}this.textColoumn.Name = "textColumn";\\
            \hspace*{1.5cm}this.textColumn.Size = new System.Drawing.Size(41, 22);\\
            \hspace*{1.5cm}this.textColumn.TabIndex = 7;\\
            \hspace*{1.5cm}//\\
            \hspace*{1.5cm}// textLine\\
            \hspace*{1.5cm}//\\
            \hspace*{1.5cm}this.textLine.Location = new System.Drawing.Point(253, 21);\\
            \hspace*{1.5cm}this.textLine.Name = "textLine";\\
            \hspace*{1.5cm}this.textLine.Size = new System.Drawing.Size(41, 22);\\
            \hspace*{1.5cm}this.textLine.TabIndex = 6;\\
            \hspace*{1.5cm}//\\
            \hspace*{1.5cm}// label3\\
            \hspace*{1.5cm}//\\
            \hspace*{1.5cm}this.label3.Anchor = ((System.Windows.Forms.AnchorStyles)(((System.\\
   Windows. Forms.AnchorStyles.Top $|$ System.Windows.Forms.AnchorStyles. Bottom)$|$\\
   System. Windows. Forms.AnchorStyles.Left)));\\
            \hspace*{1.5cm}this.label3.AutoSize = true;\\
            \hspace*{1.5cm}this.label3.Location = new System.Drawing.Point(25, 73);\\
            \hspace*{1.5cm}this.label3.Name = "label3";\\
            \hspace*{1.5cm}this.label3.Size = new System.Drawing.Size(196, 16);\\
            \hspace*{1.5cm}this.label3.TabIndex = 5;\\
            \hspace*{1.5cm}this.label3.Text = "Introduce number of columns !";\\
            \hspace*{1.5cm}//\\
            \hspace*{1.5cm}// label2\\
            \hspace*{1.5cm}//\\
            \hspace*{1.5cm}this.label2.Anchor = ((System.Windows.Forms.AnchorStyles)(((System.\\
   Windows.Forms.AnchorStyles. Top $|$ System.Windows.Forms.AnchorStyles.Bottom)$|$\\
   System.Windows.Forms.AnchorStyles.Left)));\\
            \hspace*{1.5cm}this.label2.AutoSize = true;\\
            \hspace*{1.5cm}this.label2.Location = new System.Drawing.Point(22, 24);\\
            \hspace*{1.5cm}this.label2.Name = "label2";\\
            \hspace*{1.5cm}this.label2.Size = new System.Drawing.Size(166, 16);\\
            \hspace*{1.5cm}this.label2.TabIndex = 4;\\
            \hspace*{1.5cm}this.label2.Text = "Introduce number of lines !";\\
            \hspace*{1.5cm}//\\
            \hspace*{1.5cm}// tabMultiplication\\
            \hspace*{1.5cm}//\\
            \hspace*{1.5cm}this.tabMultiplication.BackColor = System.Drawing.Color.FromArgb(((int)(((byte)(192)))), ((int)(((byte)(255)))),
  ((int)(((byte)(192)))));\\
            \hspace*{1.5cm}this.tabMultiplication.Controls.Add(this.pictureBox3);\\
            \hspace*{1.5cm}this.tabMultiplication.Controls.Add(this.btResetMultiplication);\\
            \hspace*{1.5cm}this.tabMultiplication.Controls.Add(this.btComputingProduct);\\
            \hspace*{1.5cm}this.tabMultiplication.Controls.Add(this.btGenerating2);\\
            \hspace*{1.5cm}this.tabMultiplication.Controls.Add(this.label5);\\
            \hspace*{1.5cm}this.tabMultiplication.Controls.Add(this.label13);\\
            \hspace*{1.5cm}this.tabMultiplication.Controls.Add(this.textcolumnB2);\\
            \hspace*{1.5cm}this.tabMultiplication.Controls.Add(this.textLineB2);\\
            \hspace*{1.5cm}this.tabMultiplication.Controls.Add(this.label11);\\
            \hspace*{1.5cm}this.tabMultiplication.Controls.Add(this.label12);\\
            \hspace*{1.5cm}this.tabMultiplication.Controls.Add(this.panel3);\\
            \hspace*{1.5cm}this.tabMultiplication.Controls.Add(this.panel4);\\
            \hspace*{1.5cm}this.tabMultiplication.Controls.Add(this.textColumnA2);\\
            \hspace*{1.5cm}this.tabMultiplication.Controls.Add(this.textLineA2);\\
            \hspace*{1.5cm}this.tabMultiplication.Controls.Add(this.label9);\\
            \hspace*{1.5cm}this.tabMultiplication.Controls.Add(this.label10);\\
            \hspace*{1.5cm}this.tabMultiplication.Location = new System.Drawing.Point(4, 25);\\
            \hspace*{1.5cm}this.tabMultiplication.Name = "tabMultiplication";\\
            \hspace*{1.5cm}this.tabMultiplication.Padding = new System.Windows.Forms.Padding(3);\\
            \hspace*{1.5cm}this.tabMultiplication.Size = new System.Drawing.Size(992, 667);\\
            \hspace*{1.5cm}this.tabMultiplication.TabIndex = 1;\\
            \hspace*{1.5cm}this.tabMultiplication.Text = "Multiplication";\\
            \hspace*{1.5cm}//\\
            \hspace*{1.5cm}// pictureBox3\\
            \hspace*{1.5cm}//\\
            \hspace*{1.5cm}this.pictureBox3.Image =((System.Drawing.Image)(resources.GetObject\\
  ("pictureBox3.Image")));\\
            \hspace*{1.5cm}this.pictureBox3.Location = new System.Drawing.Point(70, 220);\\
            \hspace*{1.5cm}this.pictureBox3.Name = "pictureBox3";\\
            \hspace*{1.5cm}this.pictureBox3.Size = new System.Drawing.Size(56, 27);\\
            \hspace*{1.5cm}this.pictureBox3.TabIndex = 40;\\
            \hspace*{1.5cm}this.pictureBox3.TabStop = false;\\
            \hspace*{1.5cm}//\\
            \hspace*{1.5cm}// btResetMultiplication\\
            \hspace*{1.5cm}//\\
            \hspace*{1.5cm}this.btResetMultiplication.BackColor = System.Drawing.Color.Pink;\\
            \hspace*{1.5cm}this.btResetMultiplication.ForeColor = System.Drawing.Color.Black;\\
            \hspace*{1.5cm}this.btResetMultiplication.Location = new System.Drawing.Point(548, 356);\\
            \hspace*{1.5cm}this.btResetMultiplication.Name = "btResetMultiplication";\\
            \hspace*{1.5cm}this.btResetMultiplication.Size = new System.Drawing.Size(125, 65);\\
            \hspace*{1.5cm}this.btResetMultiplication.TabIndex = 39;\\
            \hspace*{1.5cm}this.btResetMultiplication.Text = "Reset values";\\
            \hspace*{1.5cm}this.btResetMultiplication.UseVisualStyleBackColor = false;\\
            \hspace*{1.5cm}this.btResetMultiplication.Click + = new System.EventHandler(this.\\
  btResetMultiplication$_{-}$Click);\\
            \hspace*{1.5cm}//\\
            \hspace*{1.5cm}// btComputingProduct\\
            \hspace*{1.5cm}//\\
            \hspace*{1.5cm}this.btComputingProduct.BackColor = System.Drawing.Color.PowderBlue;\\
            \hspace*{1.5cm}this.btCalculareProdus.ForeColor = System.Drawing.SystemColors.ControlText;\\
            \hspace*{1.5cm}this.btComputingProduct.Location = new System.Drawing.Point(548, 253);\\
            \hspace*{1.5cm}this.btComputingProduct.Name = "btComputingProduct";\\
            \hspace*{1.5cm}this.btComputingProduct.Size = new System.Drawing.Size(125, 65);\\
            \hspace*{1.5cm}this.btComputingProduct.TabIndex = 38;\\
            \hspace*{1.5cm}this.btComputingProduct.Text = "Computing Matrix";\\
            \hspace*{1.5cm}this.btComputingProduct.UseVisualStyleBackColor = false;\\
            \hspace*{1.5cm}this.btComputingProduct.Click + = new System.EventHandler(this.\\
  btComputingProduct$_{-}$Click);\\
            \hspace*{1.5cm}//\\
            \hspace*{1.5cm}// btGenerating2\\
            \hspace*{1.5cm}//\\
            \hspace*{1.5cm}this.btGenerating2.BackColor = System.Drawing.Color.PowderBlue;\\
            \hspace*{1.5cm}this.btGenerating2.ForeColor = System.Drawing.Color.Black;\\
            \hspace*{1.5cm}this.btGenerating2.Location = new System.Drawing.Point(548, 145);\\
            \hspace*{1.5cm}this.btGenerating2.Name = "btGenerating2";\\
            \hspace*{1.5cm}this.btGenerating2.Size = new System.Drawing.Size(125, 65);\\
            \hspace*{1.5cm}this.btGenerating2.TabIndex = 37;\\
            \hspace*{1.5cm}this.btGenerating2.Text = "Generating matrix";\\
            \hspace*{1.5cm}this.btGenerating2.UseVisualStyleBackColor = false;\\
            \hspace*{1.5cm}this.btGenerating2.Click + = new System.EventHandler(this.btGenerating2$_{-}$Click);\\
            \hspace*{1.5cm}//\\
            \hspace*{1.5cm}// label5\\
            \hspace*{1.5cm}//\\
            \hspace*{1.5cm}this.label5.AutoSize = true;\\
            \hspace*{1.5cm}this.label5.ForeColor = System.Drawing.Color.Red;\\
            \hspace*{1.5cm}this.label5.Location = new System.Drawing.Point(25, 224);\\
            \hspace*{1.5cm}this.label5.Name = "label5";\\
            \hspace*{1.5cm}this.label5.Size = new System.Drawing.Size(39, 16);\\
            \hspace*{1.5cm}this.label5.TabIndex = 35;\\
            \hspace*{1.5cm}this.label5.Text = "Rem:";\\
            \hspace*{1.5cm}//\\
            \hspace*{1.5cm}// label13\\
            \hspace*{1.5cm}//\\
            \hspace*{1.5cm}this.label13.AutoSize = true;\\
            \hspace*{1.5cm}this.label13.ForeColor = System.Drawing.Color.Red;\\
            \hspace*{1.5cm}this.label13.Location = new System.Drawing.Point(14, 103);\\
            \hspace*{1.5cm}this.label13.Name = "label13";\\
            \hspace*{1.5cm}this.label13.Size = new System.Drawing.Size(808, 16);\\
            \hspace*{1.5cm}this.label13.TabIndex = 28;\\
            \hspace*{1.5cm}this.label13.Text = "Rem: Matrix A can by multiply with  matrix B, only if number of columns of " +
 "matrix A is equal with number of lines of matrix B!";\\
            \hspace*{1.5cm}//\\
            \hspace*{1.5cm}// textcolumnB2\\
            \hspace*{1.5cm}//\\
            \hspace*{1.5cm}this.textcolumnB2.Location = new System.Drawing.Point(344, 191);\\
            \hspace*{1.5cm}this.textcolumnB2.Name = "textcolumnB2";\\
            \hspace*{1.5cm}this.textcolumnB2.Size = new System.Drawing.Size(41, 22);\\
            \hspace*{1.5cm}this.textcolumnB2.TabIndex = 27;\\
            \hspace*{1.5cm}//\\
            \hspace*{1.5cm}// textLineB2\\
            \hspace*{1.5cm}//\\
            \hspace*{1.5cm}this.textLineB2.Enabled = false;\\
            \hspace*{1.5cm}this.textLineB2.Location = new System.Drawing.Point(344, 145);\\
            \hspace*{1.5cm}this.textLineB2.Name = "textLineB2";\\
            \hspace*{1.5cm}this.textLineB2.Size = new System.Drawing.Size(41, 22);\\
            \hspace*{1.5cm}this.textLineB2.TabIndex = 26;\\
            \hspace*{1.5cm}//\\
            \hspace*{1.5cm}// label11\\
            \hspace*{1.5cm}//\\
            \hspace*{1.5cm}this.label11.Anchor = ((System.Windows.Forms.AnchorStyles)(((System.Windows.\\
 Forms. AnchorStyles.Top $|$ System.Windows.Forms.AnchorStyles.Bottom)$|$ System.Windows.\\
 Forms.AnchorStyles.Left)));\\
            \hspace*{1.5cm}this.label11.AutoSize = true;\\
            \hspace*{1.5cm}this.label11.Location = new System.Drawing.Point(22, 194);\\
            \hspace*{1.5cm}this.label11.Name = "label11";\\
            \hspace*{1.5cm}this.label11.Size = new System.Drawing.Size(303, 16);\\
            \hspace*{1.5cm}this.label11.TabIndex = 25;\\
            \hspace*{1.5cm}this.label11.Text = "Introduce number of columns of matrix B !";\\
            \hspace*{1.5cm}//\\
            \hspace*{1.5cm}// label12\\
            \hspace*{1.5cm}//\\
            \hspace*{1.5cm}his.label12.Anchor = ((System.Windows.Forms.AnchorStyles)(((System.Windows.\\
 Forms.AnchorStyles.Top $|$ System.Windows.Forms. AnchorStyles.Bottom) $|$ System.Windows.\\
 Forms.AnchorStyles.Left)));\\
            \hspace*{1.5cm}this.label12.AutoSize = true;\\
            \hspace*{1.5cm}this.label12.Location = new System.Drawing.Point(22, 148);\\
            \hspace*{1.5cm}this.label12.Name = "label12";\\
            \hspace*{1.5cm}this.label12.Size = new System.Drawing.Size(273, 16);\\
            \hspace*{1.5cm}this.label12.TabIndex = 24;\\
            \hspace*{1.5cm}this.label12.Text = "Introduce number of lines of matrix B !";\\
            \hspace*{1.5cm}//\\
            \hspace*{1.5cm}// panel3\\
            \hspace*{1.5cm}//\\
            \hspace*{1.5cm}this.panel3.Controls.Add(this.pictureBox4);\\
            \hspace*{1.5cm}this.panel3.Controls.Add(this.dataGridProduct);\\
            \hspace*{1.5cm}this.panel3.Location = new System.Drawing.Point(698, 253);\\
            \hspace*{1.5cm}this.panel3.Name = "panel3";\\
            \hspace*{1.5cm}this.panel3.Size = new System.Drawing.Size(257, 297);\\
            \hspace*{1.5cm}this.panel3.TabIndex = 21;\\
            \hspace*{1.5cm}//\\
            \hspace*{1.5cm}// pictureBox4\\
            \hspace*{1.5cm}//\\
            \hspace*{1.5cm}this.pictureBox4.Image = global::Operations$_{-}$with$_{-}$matrices.Properties.Resources.AB;\\
            \hspace*{1.5cm}this.pictureBox4.Location = new System.Drawing.Point(17, 14);\\
            \hspace*{1.5cm}this.pictureBox4.Name = "pictureBox4";\\
            \hspace*{1.5cm}this.pictureBox4.Size = new System.Drawing.Size(126, 29);\\
            \hspace*{1.5cm}this.pictureBox4.TabIndex = 9;\\
            \hspace*{1.5cm}this.pictureBox4.TabStop = false;\\
            \hspace*{1.5cm}//\\
            \hspace*{1.5cm}// dataGridProduct\\
            \hspace*{1.5cm}//\\
            \hspace*{1.5cm}this.dataGridProduct.BackgroundColor = System.Drawing.Color.WhiteSmoke;\\
            \hspace*{1.5cm}this.dataGridProduct.ColumnHeadersHeightSizeMode = System.Windows.Forms.\\
 DataGridViewColumnHeadersHeightSizeMode.AutoSize;\\
            \hspace*{1.5cm}this.dataGridProduct.Location = new System.Drawing.Point(17, 49);\\
            \hspace*{1.5cm}this.dataGridProduct.Name = "dataGridProduct";\\
            \hspace*{1.5cm}this.dataGridProduct.Size = new System.Drawing.Size(223, 227);\\
            \hspace*{1.5cm}this.dataGridProduct.TabIndex = 8;\\
            \hspace*{1.5cm}//\\
            \hspace*{1.5cm}// panel4\\
            \hspace*{1.5cm}//\\
            \hspace*{1.5cm}this.panel4.Controls.Add(this.label7);\\
            \hspace*{1.5cm}this.panel4.Controls.Add(this.label8);\\
            \hspace*{1.5cm}this.panel4.Controls.Add(this.dataGridB2);\\
            \hspace*{1.5cm}this.panel4.Controls.Add(this.dataGridA2);\\
            \hspace*{1.5cm}this.panel4.Location = new System.Drawing.Point(16, 253);\\
            \hspace*{1.5cm}this.panel4.Name = "panel4";\\
            \hspace*{1.5cm}this.panel4.Size = new System.Drawing.Size(497, 297);\\
            \hspace*{1.5cm}this.panel4.TabIndex = 18;\\
            \hspace*{1.5cm}//\\
            \hspace*{1.5cm}// label7\\
            \hspace*{1.5cm}//\\
            \hspace*{1.5cm}this.label7.AutoSize = true;\\
            \hspace*{1.5cm}this.label7.ForeColor = System.Drawing.Color.MidnightBlue;\\
            \hspace*{1.5cm}this.label7.Location = new System.Drawing.Point(266, 16);\\
            \hspace*{1.5cm}this.label7.Name = "label7";\\
            \hspace*{1.5cm}this.label7.Size = new System.Drawing.Size(72, 16);\\
            \hspace*{1.5cm}this.label7.TabIndex = 11;\\
            \hspace*{1.5cm}this.label7.Text = "Matrix B";\\
            \hspace*{1.5cm}//\\
            \hspace*{1.5cm}// label8\\
            \hspace*{1.5cm}//\\
            \hspace*{1.5cm}this.label8.AutoSize = true;\\
            \hspace*{1.5cm}this.label8.ForeColor = System.Drawing.Color.MidnightBlue;\\
            \hspace*{1.5cm}this.label8.Location = new System.Drawing.Point(27, 16);\\
            \hspace*{1.5cm}this.label8.Name = "label8";\\
            \hspace*{1.5cm}this.label8.Size = new System.Drawing.Size(72, 16);\\
            \hspace*{1.5cm}this.label8.TabIndex = 10;\\
            \hspace*{1.5cm}this.label8.Text = "Matrix A";\\
            \hspace*{1.5cm}//\\
            \hspace*{1.5cm}// dataGridB2\\
            \hspace*{1.5cm}//\\
            \hspace*{1.5cm}this.dataGridB2.BackgroundColor = System.Drawing.Color.WhiteSmoke;\\
            \hspace*{1.5cm}this.dataGridB2.ColumnHeadersHeightSizeMode = System.Windows.Forms.\\
 DataGridViewColumnHeadersHeightSizeMode.AutoSize;\\
            \hspace*{1.5cm}this.dataGridB2.Location = new System.Drawing.Point(257, 49);\\
            \hspace*{1.5cm}this.dataGridB2.Name = "dataGridB2";\\
            \hspace*{1.5cm}this.dataGridB2.Size = new System.Drawing.Size(223, 227);\\
            \hspace*{1.5cm}this.dataGridB2.TabIndex = 6;\\
            \hspace*{1.5cm}//\\
            \hspace*{1.5cm}// dataGridA2\\
            \hspace*{1.5cm}//\\
            \hspace*{1.5cm}this.dataGridA2.BackgroundColor = System.Drawing.Color.WhiteSmoke;\\
            \hspace*{1.5cm}this.dataGridA2.ColumnHeadersHeightSizeMode = System.Windows.Forms.\\
 DataGridViewColumnHeadersHeightSizeMode.AutoSize;\\
            \hspace*{1.5cm}this.dataGridA2.Location = new System.Drawing.Point(12, 49);\\
            \hspace*{1.5cm}this.dataGridA2.Name = "dataGridA2";\\
            \hspace*{1.5cm}this.dataGridA2.Size = new System.Drawing.Size(223, 227);\\
            \hspace*{1.5cm}this.dataGridA2.TabIndex = 0;\\
            \hspace*{1.5cm}//\\
            \hspace*{1.5cm}// textColumnA2\\
            \hspace*{1.5cm}//\\
            \hspace*{1.5cm}this.textColumnA2.Location = new System.Drawing.Point(344, 63);\\
            \hspace*{1.5cm}this.textColumnA2.Name = "textColumnA2";\\
            \hspace*{1.5cm}this.textColumnA2.Size = new System.Drawing.Size(41, 22);\\
            \hspace*{1.5cm}this.textColumnA2.TabIndex = 17;\\
            \hspace*{1.5cm}this.textColumnA2.Leave + = new System.EventHandler(this.textColumnA2$_{-}$Leave);\\
            \hspace*{1.5cm}//\\
            \hspace*{1.5cm}// textLineA2\\
            \hspace*{1.5cm}//\\
            \hspace*{1.5cm}this.textLineA2.Location = new System.Drawing.Point(344, 17);\\
            \hspace*{1.5cm}this.textLineA2.Name = "textLineA2";\\
            \hspace*{1.5cm}this.textLineA2.Size = new System.Drawing.Size(41, 22);\\
            \hspace*{1.5cm}this.textLineA2.TabIndex = 16;\\
            \hspace*{1.5cm}//\\
            \hspace*{1.5cm}// label9\\
            \hspace*{1.5cm}//\\
            \hspace*{1.5cm}this.label9.Anchor = ((System.Windows.Forms.AnchorStyles)(((System.Windows.\\
Forms.AnchorStyles.Top $|$ System.Windows.Forms.AnchorStyles.Bottom) $|$ System.Windows.\\
Forms.AnchorStyles.Left)));\\
            \hspace*{1.5cm}this.label9.AutoSize = true;\\
            \hspace*{1.5cm}this.label9.ForeColor = System.Drawing.Color.Black;\\
            \hspace*{1.5cm}this.label9.Location = new System.Drawing.Point(22, 66);\\
            \hspace*{1.5cm}this.label9.Name = "label9";\\
            \hspace*{1.5cm}this.label9.Size = new System.Drawing.Size(303,16);\\
            \hspace*{1.5cm}this.label9.TabIndex = 15;\\
            \hspace*{1.5cm}this.label9.Text = "Introduce number of columns of matrix A !";\\
            \hspace*{1.5cm}//\\
            \hspace*{1.5cm}// label10\\
            \hspace*{1.5cm}//\\
            \hspace*{1.5cm}this.label10.Anchor = ((System.Windows.Forms.AnchorStyles)(((System.Windows.\\
Forms.AnchorStyles.Top $|$ System.Windows.Forms.AnchorStyles.Bottom)$|$ System.Windows.\\
Forms.AnchorStyles.Left)));\\
            \hspace*{1.5cm}this.label10.AutoSize = true;\\
            \hspace*{1.5cm}this.label10.Location = new System.Drawing.Point(22, 17);\\
            \hspace*{1.5cm}this.label10.Name = "label10";\\
            \hspace*{1.5cm}this.label10.Size = new System.Drawing.Size(273, 16);\\
            \hspace*{1.5cm}this.label10.TabIndex = 14;\\
            \hspace*{1.5cm}this.label10.Text = "Introduce  number of lines of matrix A !";\\
            \hspace*{1.5cm}//\\
            \hspace*{1.5cm}// tabMultiplication$_{-}$scalar\\
            \hspace*{1.5cm}//\\
            \hspace*{1.5cm}this.tabMultiplication$_{-}$scalar.BackColor = System.Drawing.Color.FromArgb\\
(((int)(((byte)(192)))),((int)(((byte)(255)))), ((int)(((byte)(192)))));\\
            \hspace*{1.5cm}this.tabMultiplication$_{-}$scalar.Controls.Add(this.pictureBox6);\\
            \hspace*{1.5cm}this.tabMultiplication$_{-}$scalar.Controls.Add(this.btReset$_{-}$Product$_{-}$Scalar);\\
            \hspace*{1.5cm}this.tabMultiplication$_{-}$scalar.Controls.Add(this.btProductScalar);\\
            \hspace*{1.5cm}this.tabMultiplication$_{-}$scalar.Controls.Add(this.butonsGenerating3);\\
            \hspace*{1.5cm}this.tabMultiplication$_{-}$scalar.Controls.Add(this.label28);\\
            \hspace*{1.5cm}this.tabMultiplication$_{-}$scalar.Controls.Add(this.textScalar);\\
            \hspace*{1.5cm}this.tabMultiplication$_{-}$scalar.Controls.Add(this.label16);\\
            \hspace*{1.5cm}this.tabMultiplication$_{-}$scalar.Controls.Add(this.panel6);\\
            \hspace*{1.5cm}this.tabMultiplication$_{-}$scalar.Controls.Add(this.panel5);\\
            \hspace*{1.5cm}this.tabMultiplication$_{-}$scalar.Controls.Add(this.textcolumnA3);\\
            \hspace*{1.5cm}this.tabMultiplication$_{-}$scalar.Controls.Add(this.textLineA3);\\
            \hspace*{1.5cm}this.tabMultiplication$_{-}$scalar.Controls.Add(this.label14);\\
            \hspace*{1.5cm}this.tabMultiplication$_{-}$scalar.Controls.Add(this.label15);\\
            \hspace*{1.5cm}this.tabMultiplication$_{-}$scalar.Location = new System.Drawing.Point(4, 25);\\
            \hspace*{1.5cm}this.tabMultiplication$_{-}$scalar.Name = "tabMultiplication$_{-}$scalar";\\
            \hspace*{1.5cm}this.tabMultiplication$_{-}$scalar.Padding = new System.Windows.Forms.Padding(3);\\
            \hspace*{1.5cm}this.tabMultiplication$_{-}$scalar.Size = new System.Drawing.Size(992, 667);\\
            \hspace*{1.5cm}\hspace*{1.5cm}this.tabMultiplication$_{-}$scalar.TabIndex = 2;\\
            \hspace*{1.5cm}this.tabMultiplication$_{-}$scalar.Text = "Multiplication with a scalar";\\
            \hspace*{1.5cm}//\\
            \hspace*{1.5cm}// pictureBox6\\
            \hspace*{1.5cm}//\\
            \hspace*{1.5cm}this.pictureBox6.Image = ((System.Drawing.Image)(resources.GetObject\\
("pictureBox6.Image")));\\
            \hspace*{1.5cm}this.pictureBox6.Location = new System.Drawing.Point(56, 162);\\
            \hspace*{1.5cm}this.pictureBox6.Name = "pictureBox6";\\
            \hspace*{1.5cm}this.pictureBox6.Size = new System.Drawing.Size(56, 27);\\
            \hspace*{1.5cm}this.pictureBox6.TabIndex = 43;\\
            \hspace*{1.5cm}this.pictureBox6.TabStop = false;\\
            \hspace*{1.5cm}//\\
            \hspace*{1.5cm}// btReset$_{-}$Product$_{-}$Scalar\\
            \hspace*{1.5cm}//\\
            \hspace*{1.5cm}this.btReset$_{-}$Product$_{-}$Scalar.BackColor = System.Drawing.Color.Pink;\\
            \hspace*{1.5cm}this.btReset$_{-}$Product$_{-}$Scalar.ForeColor = System.Drawing.Color.Black;\\
            \hspace*{1.5cm}this.btReset$_{-}$Product$_{-}$Scalar.Location = new System.Drawing.Point(316, 421);\\
            \hspace*{1.5cm}this.btReset$_{-}$Product$_{-}$Scalar.Name = "btReset$_{-}$Product$_{-}$Scalar";\\
            \hspace*{1.5cm}this.btReset$_{-}$Product$_{-}$Scalar.Size = new System.Drawing.Size(125, 65);\\
            \hspace*{1.5cm}this.btReset$_{-}$Product$_{-}$Scalar.TabIndex = 42;\\
            \hspace*{1.5cm}this.btReset$_{-}$Product$_{-}$Scalar.Text = "Reset values";\\
            \hspace*{1.5cm}this.btReset$_{-}$Product$_{-}$Scalar.UseVisualStyleBackColor = false;\\
            \hspace*{1.5cm}this.btReset$_{-}$Product$_{-}$Scalar.Click + = new System.EventHandler(this.\\
btReset$_{-}$Product$_{-}$Scalar$_{-}$Click);\\
            \hspace*{1.5cm}//\\
            \hspace*{1.5cm}// btProductScalar\\
            \hspace*{1.5cm}//\\
            \hspace*{1.5cm}this.btProductScalar.BackColor = System.Drawing.Color.PowderBlue;\\
            \hspace*{1.5cm}this.btProdusctScalar.ForeColor = System.Drawing.SystemColors.ControlText;\\
            \hspace*{1.5cm}this.btProductScalar.Location = new System.Drawing.Point(316, 318);\\
            \hspace*{1.5cm}this.btProductScalar.Name = "btProductScalar";\\
            \hspace*{1.5cm}this.btProductScalar.Size = new System.Drawing.Size(125, 65);\\
            \hspace*{1.5cm}this.btProductScalar.TabIndex = 41;\\
            \hspace*{1.5cm}this.btProductScalar.Text = "Computing Matrix";\\
            \hspace*{1.5cm}this.btProductScalar.UseVisualStyleBackColor = false;\\
            \hspace*{1.5cm}this.btProductScalar.Click + = new System.EventHandler(this.btProductScalar$_{-}$Click);\\
            \hspace*{1.5cm}//\\
            \hspace*{1.5cm}// butonsGenerating3\\
            \hspace*{1.5cm}//\\
            \hspace*{1.5cm}this.butonsGenerating3.BackColor = System.Drawing.Color.PowderBlue;\\
            \hspace*{1.5cm}this.butonsGenerating3.ForeColor = System.Drawing.SystemColors.ControlText;\\
            \hspace*{1.5cm}this.butonsGenerating3.Location = new System.Drawing.Point(316, 210);\\
            \hspace*{1.5cm}this.butonsGenerating3.Name = "butonsGenerating3";\\
            \hspace*{1.5cm}this.butonsGenerating3.Size = new System.Drawing.Size(125, 65);\\
            \hspace*{1.5cm}this.butonsGenerating3.TabIndex = 40;\\
            \hspace*{1.5cm}this.butonsGenerating3.Text = "Generating of matrix";\\
            \hspace*{1.5cm}this.butonsGenerating3.UseVisualStyleBackColor = false;\\
            \hspace*{1.5cm}this.butonsGenerating3.Click + = new System.EventHandler(this.butonsGenerating3$_{-}$Click);\\
            \hspace*{1.5cm}//\\
            \hspace*{1.5cm}// label28\\
            \hspace*{1.5cm}//\\
            \hspace*{1.5cm}this.label28.AutoSize = true;\\
            \hspace*{1.5cm}this.label28.ForeColor = System.Drawing.Color.Red;\\
            \hspace*{1.5cm}this.label28.Location = new System.Drawing.Point(12, 169);\\
            \hspace*{1.5cm}this.label28.Name = "label28";\\
            \hspace*{1.5cm}this.label28.Size = new System.Drawing.Size(39, 16);\\
            \hspace*{1.5cm}this.label28.TabIndex = 35;\\
            \hspace*{1.5cm}this.label28.Text = "Rem:";\\
            \hspace*{1.5cm}//\\
            \hspace*{1.5cm}// textScalar\\
            \hspace*{1.5cm}//\\
            \hspace*{1.5cm}this.textScalar.Location = new System.Drawing.Point(143, 121);\\
            \hspace*{1.5cm}this.textScalar.Name = "textScalar";\\
            \hspace*{1.5cm}this.textScalar.Size = new System.Drawing.Size(41, 22);\\
            \hspace*{1.5cm}this.textScalar.TabIndex = 19;\\
            \hspace*{1.5cm}//\\
            \hspace*{1.5cm}// label16\\
            \hspace*{1.5cm}//\\
            \hspace*{1.5cm}this.label16.Anchor = ((System.Windows.Forms.AnchorStyles)(((System.Windows.\\
Forms.AnchorStyles.Top $|$ System.Windows.Forms.AnchorStyles.Bottom) $|$ System.Windows.\\
Forms.AnchorStyles.Left)));\\
            \hspace*{1.5cm}this.label16.AutoSize = true;\\
            \hspace*{1.5cm}this.label16.Location = new System.Drawing.Point(12, 127);\\
            \hspace*{1.5cm}this.label16.Name = "label16";\\
            \hspace*{1.5cm}this.label16.Size = new System.Drawing.Size(125, 16);\\
            \hspace*{1.5cm}this.label16.TabIndex = 18;\\
            \hspace*{1.5cm}this.label16.Text = "Introduce the scalar !";\\
            \hspace*{1.5cm}//\\
            \hspace*{1.5cm}// panel6\\
            \hspace*{1.5cm}//\\
            \hspace*{1.5cm}this.panel6.Controls.Add(this.pictureBox5);\\
            \hspace*{1.5cm}this.panel6.Controls.Add(this.dataGridProductScalar);\\
            \hspace*{1.5cm}this.panel6.Location = new System.Drawing.Point(498, 201);\\
            \hspace*{1.5cm}this.panel6.Name = "panel6";\\
            \hspace*{1.5cm}this.panel6.Size = new System.Drawing.Size(257, 297);\\
            \hspace*{1.5cm}this.panel6.TabIndex = 17;\\
            \hspace*{1.5cm}//\\
            \hspace*{1.5cm}// pictureBox5\\
            \hspace*{1.5cm}//\\
            \hspace*{1.5cm}this.pictureBox5.BackColor = System.Drawing.Color.Red;\\
            \hspace*{1.5cm}this.pictureBox5.Image = global::Operations$_{-}$with$_{-}$matrices.Properties.Resources.scalar;\\
            \hspace*{1.5cm}this.pictureBox5.Location = new System.Drawing.Point(17, 10);\\
            \hspace*{1.5cm}this.pictureBox5.Name = "pictureBox5";\\
            \hspace*{1.5cm}this.pictureBox5.Size = new System.Drawing.Size(147, 31);\\
            \hspace*{1.5cm}this.pictureBox5.TabIndex = 13;\\
            \hspace*{1.5cm}this.pictureBox5.TabStop = false;\\
            \hspace*{1.5cm}//\\
            \hspace*{1.5cm}// dataGridProductScalar\\
            \hspace*{1.5cm}//\\
            \hspace*{1.5cm}this.dataGridProductScalar.BackgroundColor = System.Drawing.Color.WhiteSmoke;\\
            \hspace*{1.5cm}this.dataGridProductScalar.ColumnHeadersHeightSizeMode = System.Windows.Forms.\\
DataGridViewColumnHeadersHeightSizeMode.AutoSize;\\
            \hspace*{1.5cm}this.dataGridProductScalar.Location = new System.Drawing.Point(17, 49);\\
            \hspace*{1.5cm}this.dataGridProductScalar.Name = "dataGridProductScalar";\\
            \hspace*{1.5cm}this.dataGridProductScalar.Size = new System.Drawing.Size(223, 227);\\
            \hspace*{1.5cm}this.dataGridProductScalar.TabIndex = 8;\\
            \hspace*{1.5cm}//\\
            \hspace*{1.5cm}// panel5\\
            \hspace*{1.5cm}//\\
            \hspace*{1.5cm}this.panel5.Controls.Add(this.label17);\\
            \hspace*{1.5cm}this.panel5.Controls.Add(this.dataGridA3);\\
            \hspace*{1.5cm}this.panel5.Location = new System.Drawing.Point(13, 210);\\
            \hspace*{1.5cm}this.panel5.Name = "panel5";\\
            \hspace*{1.5cm}this.panel5.Size = new System.Drawing.Size(269, 297);\\
            \hspace*{1.5cm}this.panel5.TabIndex = 15;\\
            \hspace*{1.5cm}//\\
            \hspace*{1.5cm}// label17\\
            \hspace*{1.5cm}//\\
            \hspace*{1.5cm}this.label17.AutoSize = true;\\
            \hspace*{1.5cm}this.label17.ForeColor = System.Drawing.Color.MidnightBlue;\\
            \hspace*{1.5cm}this.label17.Location = new System.Drawing.Point(27, 16);\\
            \hspace*{1.5cm}this.label17.Name = "label17";\\
            \hspace*{1.5cm}this.label17.Size = new System.Drawing.Size(72, 16);\\
            \hspace*{1.5cm}this.label17.TabIndex = 10;\\
            \hspace*{1.5cm}this.label17.Text = "Matrix A";\\
            \hspace*{1.5cm}//\\
            \hspace*{1.5cm}// dataGridA3\\
            \hspace*{1.5cm}//\\
            \hspace*{1.5cm}this.dataGridA3.BackgroundColor = System.Drawing.Color.WhiteSmoke;\\
            \hspace*{1.5cm}this.dataGridA3.ColumnHeadersHeightSizeMode = System.Windows.Forms.\\
DataGridViewColumnHeadersHeightSizeMode.AutoSize;\\
            \hspace*{1.5cm}this.dataGridA3.Location = new System.Drawing.Point(12, 49);\\
            \hspace*{1.5cm}this.dataGridA3.Name = "dataGridA3";\\
            \hspace*{1.5cm}this.dataGridA3.Size = new System.Drawing.Size(223, 227);\\
            \hspace*{1.5cm}this.dataGridA3.TabIndex = 0;\\
            \hspace*{1.5cm}//\\
            \hspace*{1.5cm}// textcolumnA3\\
            \hspace*{1.5cm}//\\
            \hspace*{1.5cm}this.textcolumnA3.Location = new System.Drawing.Point(318, 68);\\
            \hspace*{1.5cm}this.textcolumnA3.Name = "textcolumnA3";\\
            \hspace*{1.5cm}this.textcolumnA3.Size = new System.Drawing.Size(41, 22);\\
            \hspace*{1.5cm}this.textcolumnA3.TabIndex = 13;\\
            \hspace*{1.5cm}//\\
            \hspace*{1.5cm}// textLineA3\\
            \hspace*{1.5cm}//\\
            \hspace*{1.5cm}this.textLineA3.Location = new System.Drawing.Point(318, 19);\\
            \hspace*{1.5cm}this.textLineA3.Name = "textLineA3";\\
            \hspace*{1.5cm}this.textLineA3.Size = new System.Drawing.Size(41, 22);\\
            \hspace*{1.5cm}this.textLineA3.TabIndex = 12;\\
            \hspace*{1.5cm}//\\
            \hspace*{1.5cm}// label14\\
            \hspace*{1.5cm}//\\
            \hspace*{1.5cm}this.label14.Anchor = ((System.Windows.Forms.AnchorStyles)(((System.Windows.\\
 Forms.AnchorStyles.Top $|$ System.Windows.Forms.AnchorStyles.Bottom) $|$ System.Windows.\\
 Forms.AnchorStyles.Left)));\\
            \hspace*{1.5cm}this.label14.AutoSize = true;\\
            \hspace*{1.5cm}this.label14.Location = new System.Drawing.Point(12, 71);\\
            \hspace*{1.5cm}this.label14.Name = "label14";\\
            \hspace*{1.5cm}this.label14.Size = new System.Drawing.Size(300, 16);\\
            \hspace*{1.5cm}this.label14.TabIndex = 11;\\
            \hspace*{1.5cm}this.label14.Text = "Introduce number of columns of matrix A!";\\
            \hspace*{1.5cm}//\\
            \hspace*{1.5cm}// label15\\
            \hspace*{1.5cm}//\\
            \hspace*{1.5cm}this.label15.Anchor = ((System.Windows.Forms.AnchorStyles)(((System.Windows.\\
 Forms.AnchorStyles.Top $|$ System.Windows.Forms.AnchorStyles.Bottom)$|$ System.Windows.\\
 Forms.AnchorStyles.Left)));\\
            \hspace*{1.5cm}this.label15.AutoSize = true;\\
            \hspace*{1.5cm}this.label15.Location = new System.Drawing.Point(12, 22);\\
            \hspace*{1.5cm}this.label15.Name = "label15";\\
            \hspace*{1.5cm}this.label15.Size = new System.Drawing.Size(273, 16);\\
            \hspace*{1.5cm}this.label15.TabIndex = 10;\\
            \hspace*{1.5cm}this.label15.Text = "Introduce number of lines of matrix A !";\\
            \hspace*{1.5cm}//\\
            \hspace*{1.5cm}// tabLifting$_{-}$at$_{-}$power\\
            \hspace*{1.5cm}//\\
            \hspace*{1.5cm}this.tabLifting$_{-}$at$_{-}$power.BackColor = System.Drawing.Color.FromArgb(((int)\\
 (((byte)(192)))),((int)(((byte)(255)))),((int)(((byte)(192)))));\\
            \hspace*{1.5cm}this.tabLifting$_{-}$at$_{-}$power.Controls.Add(this.pictureBox9);\\
            \hspace*{1.5cm}this.tabLifting$_{-}$at$_{-}$power.Controls.Add(this.btReset$_{-}$lifting$_{-}$at$_{-}$power);\\
            \hspace*{1.5cm}this.tabLifting$_{-}$at$_{-}$power.Controls.Add(this.btComputPower);\\
            \hspace*{1.5cm}this.tabLifting$_{-}$at$_{-}$power.Controls.Add(this.btGeneratingA4);\\
            \hspace*{1.5cm}this.tabLifting$_{-}$at$_{-}$power.Controls.Add(this.label29);\\
            \hspace*{1.5cm}this.tabLifting$_{-}$at$_{-}$power.Controls.Add(this.label21);\\
            \hspace*{1.5cm}this.tabLifting$_{-}$at$_{-}$power.Controls.Add(this.panel7);\\
            \hspace*{1.5cm}this.tabLifting$_{-}$at$_{-}$power.Controls.Add(this.panel8);\\
            \hspace*{1.5cm}this.tabLifting$_{-}$at$_{-}$power.Controls.Add(this.textPower);\\
            \hspace*{1.5cm}this.tabLifting$_{-}$at$_{-}$power.Controls.Add(this.textlineA4);\\
            \hspace*{1.5cm}this.tabLifting$_{-}$at$_{-}$power.Controls.Add(this.label18);\\
            \hspace*{1.5cm}this.tabLifting$_{-}$at$_{-}$power.Controls.Add(this.label19);\\
            \hspace*{1.5cm}this.tabLifting$_{-}$at$_{-}$power.Location = new System.Drawing.Point(4, 25);\\
            \hspace*{1.5cm}this.tabLifting$_{-}$at$_{-}$power.Name = "tabLifting$_{-}$at$_{-}$power";\\
            \hspace*{1.5cm}this.tabLifting$_{-}$at$_{-}$power.Padding = new System.Windows.Forms.Padding(3);\\
            \hspace*{1.5cm}this.tabLifting$_{-}$at$_{-}$power.Size = new System.Drawing.Size(992, 667);\\
            \hspace*{1.5cm}this.tabLifting$_{-}$at$_{-}$power.TabIndex = 3;\\
            \hspace*{1.5cm}this.tabLifting$_{-}$at$_{-}$power.Text = "Lifting at power";\\
            \hspace*{1.5cm}//\\
            \hspace*{1.5cm}// pictureBox9\\
            \hspace*{1.5cm}//\\
            \hspace*{1.5cm}this.pictureBox9.Image = ((System.Drawing.Image)(resources.GetObject\\
 ("pictureBox9.Image")));\\
            \hspace*{1.5cm}this.pictureBox9.Location = new System.Drawing.Point(59, 155);\\
            \hspace*{1.5cm}this.pictureBox9.Name = "pictureBox9";\\
            \hspace*{1.5cm}this.pictureBox9.Size = new System.Drawing.Size(56, 27);\\
            \hspace*{1.5cm}this.pictureBox9.TabIndex = 46;\\
            \hspace*{1.5cm}this.pictureBox9.TabStop = false;\\
            \hspace*{1.5cm}//\\
            \hspace*{1.5cm}//btReset$_{-}$lifting$_{-}$at$_{-}$power\\
            \hspace*{1.5cm}//\\
            \hspace*{1.5cm}this.btReset$_{-}$lifting$_{-}$at$_{-}$power.BackColor = System.Drawing.Color.Pink;\\
            \hspace*{1.5cm}this.btReset$_{-}$lifting$_{-}$at$_{-}$power.ForeColor = System.Drawing.Color.Black;\\
            \hspace*{1.5cm}this.btReset$_{-}$lifting$_{-}$at$_{-}$power.Location = new System.Drawing.Point(319, 420);\\
            \hspace*{1.5cm}this.btReset$_{-}$lifting$_{-}$at$_{-}$power.Name = "btReset$_{-}$lifting$_{-}$at$_{-}$power";\\
            \hspace*{1.5cm}this.btReset$_{-}$lifting$_{-}$at$_{-}$power.Size = new System.Drawing.Size(125, 65);\\
            \hspace*{1.5cm}this.btReset$_{-}$lifting$_{-}$at$_{-}$power.TabIndex = 45;\\
            \hspace*{1.5cm}this.btReset$_{-}$lifting$_{-}$at$_{-}$power.Text = "Reset values";\\
            \hspace*{1.5cm}this.btReset$_{-}$lifting$_{-}$at$_{-}$power.UseVisualStyleBackColor = false;\\
            \hspace*{1.5cm}this.btReset$_{-}$lifting$_{-}$at$_{-}$power.Click + = new System.EventHandler(this.\\
 btReset$_{-}$lifting$_{-}$at$_{-}$power$_{-}$Click);\\
            \hspace*{1.5cm}//\\
            \hspace*{1.5cm}// btComputPower\\
            \hspace*{1.5cm}//\\
            \hspace*{1.5cm}this.btComputPower.BackColor = System.Drawing.Color.PowderBlue;\\
            \hspace*{1.5cm}this.btComputPower.ForeColor = System.Drawing.SystemColors.ControlText;\\
            \hspace*{1.5cm}this.btComputPower.Location = new System.Drawing.Point(319, 317);\\
            \hspace*{1.5cm}this.btComputPower.Name = "btComputPower";\\
            \hspace*{1.5cm}this.btComputPower.Size = new System.Drawing.Size(125, 65);\\
            \hspace*{1.5cm}this.btComputPower.TabIndex = 44;\\
            \hspace*{1.5cm}this.btComputPower.Text = "Comput Matrix";\\
            \hspace*{1.5cm}this.btComputPower.UseVisualStyleBackColor = false;\\
            \hspace*{1.5cm}this.btComputPower.Click + = new System.EventHandler(this.btComputPower$_{-}$Click);\\
            \hspace*{1.5cm}//\\
            \hspace*{1.5cm}// btGeneratingA4\\
            \hspace*{1.5cm}//\\
            \hspace*{1.5cm}this.btGeneratingA4.BackColor = System.Drawing.Color.PowderBlue;\\
            \hspace*{1.5cm}this.btGeneratingA4.ForeColor = System.Drawing.SystemColors.ControlText;\\
            \hspace*{1.5cm}this.btGeneratingA4.Location = new System.Drawing.Point(319, 209);\\
            \hspace*{1.5cm}this.btGeneratingA4.Name = "btGeneratingA4";\\
            \hspace*{1.5cm}this.btGeneratingA4.Size = new System.Drawing.Size(125, 65);\\
            \hspace*{1.5cm}this.btGeneratingA4.TabIndex = 43;\\
            \hspace*{1.5cm}this.btGeneratingA4.Text = "Generating matrix";\\
            \hspace*{1.5cm}this.btGeneratingA4.UseVisualStyleBackColor = false;\\
            \hspace*{1.5cm}this.btGeneratingA4.Click += new System.EventHandler(this.btGeneratingA4$_{-}$Click);\\
            \hspace*{1.5cm}//\\
            \hspace*{1.5cm}// label29\\
            \hspace*{1.5cm}//\\
            \hspace*{1.5cm}this.label29.AutoSize = true;\\
            \hspace*{1.5cm}this.label29.ForeColor = System.Drawing.Color.Red;\\
            \hspace*{1.5cm}this.label29.Location = new System.Drawing.Point(19, 160);\\
            \hspace*{1.5cm}this.label29.Name = "label29";\\
            \hspace*{1.5cm}this.label29.Size = new System.Drawing.Size(39, 16);\\
            \hspace*{1.5cm}this.label29.TabIndex = 33;\\
            \hspace*{1.5cm}this.label29.Text = "Rem:";\\
            \hspace*{1.5cm}//\\
            \hspace*{1.5cm}// label21\\
            \hspace*{1.5cm}//\\
            \hspace*{1.5cm}this.label21.AutoSize = true;\\
            \hspace*{1.5cm}this.label21.ForeColor = System.Drawing.Color.Red;\\
            \hspace*{1.5cm}this.label21.Location = new System.Drawing.Point(21, 127);\\
            \hspace*{1.5cm}this.label21.Name = "label21";\\
            \hspace*{1.5cm}this.label21.Size = new System.Drawing.Size(314, 16);\\
            \hspace*{1.5cm}this.label21.TabIndex = 23;\\
            \hspace*{1.5cm}this.label21.Text = "Rem: Can be lifting at power only quadratic matrices";\\
            \hspace*{1.5cm}//\\
            \hspace*{1.5cm}// panel7\\
            \hspace*{1.5cm}//\\
            \hspace*{1.5cm}this.panel7.Controls.Add(this.label27);\\
            \hspace*{1.5cm}this.panel7.Controls.Add(this.label25);\\
            \hspace*{1.5cm}this.panel7.Controls.Add(this.label26);\\
            \hspace*{1.5cm}this.panel7.Controls.Add(this.labelPower);\\
            \hspace*{1.5cm}this.panel7.Controls.Add(this.label24);\\
            \hspace*{1.5cm}this.panel7.Controls.Add(this.dataGridMatrix$_{-}$at$_{-}$power$_{-}$n);\\
            \hspace*{1.5cm}this.panel7.Location = new System.Drawing.Point(501, 200);\\
            \hspace*{1.5cm}this.panel7.Name = "panel7";\\
            \hspace*{1.5cm}this.panel7.Size = new System.Drawing.Size(257, 297);\\
            \hspace*{1.5cm}this.panel7.TabIndex = 22;\\
            \hspace*{1.5cm}//\\
            \hspace*{1.5cm}// label27\\
            \hspace*{1.5cm}//\\
            \hspace*{1.5cm}this.label27.AutoSize = true;\\
            \hspace*{1.5cm}this.label27.ForeColor = System.Drawing.Color.Navy;\\
            \hspace*{1.5cm}this.label27.Location = new System.Drawing.Point(34, 25);\\
            \hspace*{1.5cm}this.label27.Name = "label27";\\
            \hspace*{1.5cm}this.label27.Size = new System.Drawing.Size(60, 16);\\
            \hspace*{1.5cm}this.label27.TabIndex = 18;\\
            \hspace*{1.5cm}this.label27.Text = "Matrix";\\
            \hspace*{1.5cm}//\\
            \hspace*{1.5cm}// label25\\
            \hspace*{1.5cm}//\\
            \hspace*{1.5cm}this.label25.AutoSize = true;\\
            \hspace*{1.5cm}this.label25.Font = new System.Drawing.Font("Microsoft Sans Serif", 8.25F,\\
  System.Drawing.FontStyle.Regular, System.Drawing.GraphicsUnit.Point, ((byte)(0)));\\
            \hspace*{1.5cm}this.label25.ForeColor = System.Drawing.Color.Navy;\\
            \hspace*{1.5cm}this.label25.Location = new System.Drawing.Point(128, 13);\\
            \hspace*{1.5cm}this.label25.Name = "label25";\\
            \hspace*{1.5cm}this.label25.Size = new System.Drawing.Size(10, 13);\\
            \hspace*{1.5cm}this.label25.TabIndex = 17;\\
            \hspace*{1.5cm}this.label25.Text = ")";\\
            \hspace*{1.5cm}//\\
            \hspace*{1.5cm}// label26\\
            \hspace*{1.5cm}//\\
            \hspace*{1.5cm}this.label26.AutoSize = true;\\
            \hspace*{1.5cm}this.label26.Font = new System.Drawing.Font("Microsoft Sans Serif", 8.25F,\\
  System.Drawing.FontStyle.Regular, System.Drawing.GraphicsUnit.Point, ((byte)(0)));\\
            \hspace*{1.5cm}this.label26.ForeColor = System.Drawing.Color.Navy;\\
            \hspace*{1.5cm}this.label26.Location = new System.Drawing.Point(107, 13);\\
            \hspace*{1.5cm}this.label26.Name = "label26";\\
            \hspace*{1.5cm}this.label26.Size = new System.Drawing.Size(10, 13);\\
            \hspace*{1.5cm}this.label26.TabIndex = 16;\\
            \hspace*{1.5cm}this.label26.Text = "(";\\
            \hspace*{1.5cm}//\\
            \hspace*{1.5cm}// labelPower\\
            \hspace*{1.5cm}//\\
            \hspace*{1.5cm}this.labelPower.AutoSize = true;\\
            \hspace*{1.5cm}this.labelPower.Font = new System.Drawing.Font("Microsoft Sans Serif", 8.25F, System.Drawing.FontStyle.
  Regular, System.Drawing.GraphicsUnit.Point, ((byte)(0)));\\
            \hspace*{1.5cm}this.labelPower.ForeColor = System.Drawing.Color.Navy;\\
            \hspace*{1.5cm}this.labelPower.Location = new System.Drawing.Point(114, 13);\\
            \hspace*{1.5cm}this.labelPower.Name = "labelPower";\\
            \hspace*{1.5cm}this.labelPower.Size = new System.Drawing.Size(13, 13);\\
            \hspace*{1.5cm}this.labelPower.TabIndex = 15;\\
            \hspace*{1.5cm}this.labelPower.Text = "1";\\
            \hspace*{1.5cm}//\\
            \hspace*{1.5cm}// label24\\
            \hspace*{1.5cm}//\\
            \hspace*{1.5cm}this.label24.AutoSize = true;\\
            \hspace*{1.5cm}this.label24.ForeColor = System.Drawing.Color.Navy;\\
            \hspace*{1.5cm}this.label24.Location = new System.Drawing.Point(100, 25);\\
            \hspace*{1.5cm}this.label24.Name = "label24";\\
            \hspace*{1.5cm}this.label24.Size = new System.Drawing.Size(17, 16);\\
            \hspace*{1.5cm}this.label24.TabIndex = 14;\\
            \hspace*{1.5cm}this.label24.Text = "A";\\
            \hspace*{1.5cm}//\\
            \hspace*{1.5cm}// dataGridMatrix$_{-}$at$_{-}$power$_{-}$n\\
            \hspace*{1.5cm}//\\
            \hspace*{1.5cm}this.dataGridMatrix$_{-}$at$_{-}$power$_{-}$n.BackgroundColor = System.Drawing.Color.WhiteSmoke;\\
            \hspace*{1.5cm}this.dataGridMatrix$_{-}$at$_{-}$power$_{-}$n.ColumnHeadersHeightSizeMode =\\
  System.Windows.Forms.DataGridViewColumnHeadersHeightSizeMode.AutoSize;\\
            \hspace*{1.5cm}this.dataGridMatrix$_{-}$at$_{-}$power$_{-}$n.Location = new System.Drawing.Point(17, 49);\\
            \hspace*{1.5cm}this.dataGridMatrix$_{-}$at$_{-}$power$_{-}$n.Name = "dataGridMatrix$_{-}$at$_{-}$power$_{-}$n";\\
            \hspace*{1.5cm}this.dataGridMatrix$_{-}$at$_{-}$power$_{-}$n.Size = new System.Drawing.Size(223, 227);\\
            \hspace*{1.5cm}this.dataGridMatrix$_{-}$at$_{-}$power$_{-}$n.TabIndex = 8;\\
            \hspace*{1.5cm}//\\
            \hspace*{1.5cm}// panel8\\
            \hspace*{1.5cm}//\\
            \hspace*{1.5cm}this.panel8.Controls.Add(this.label20);\\
            \hspace*{1.5cm}this.panel8.Controls.Add(this.dataGridA4);\\
            \hspace*{1.5cm}this.panel8.Location = new System.Drawing.Point(16, 209);\\
            \hspace*{1.5cm}this.panel8.Name = "panel8";\\
            \hspace*{1.5cm}this.panel8.Size = new System.Drawing.Size(269, 297);\\
            \hspace*{1.5cm}this.panel8.TabIndex = 20;\\
            \hspace*{1.5cm}//\\
            \hspace*{1.5cm}// label20\\
            \hspace*{1.5cm}//\\
            \hspace*{1.5cm}this.label20.AutoSize = true;\\
            \hspace*{1.5cm}this.label20.ForeColor = System.Drawing.Color.MidnightBlue;\\
            \hspace*{1.5cm}this.label20.Location = new System.Drawing.Point(27, 16);\\
            \hspace*{1.5cm}this.label20.Name = "label20";\\
            \hspace*{1.5cm}this.label20.Size = new System.Drawing.Size(72, 16);\\
            \hspace*{1.5cm}this.label20.TabIndex = 10;\\
            \hspace*{1.5cm}this.label20.Text = "Matrix A";\\
            \hspace*{1.5cm}//\\
            \hspace*{1.5cm}// dataGridA4\\
            \hspace*{1.5cm}//\\
            \hspace*{1.5cm}this.dataGridA4.BackgroundColor = System.Drawing.Color.WhiteSmoke;\\
            \hspace*{1.5cm}this.dataGridA4.ColumnHeadersHeightSizeMode = System.Windows.Forms.\\
   DataGridViewColumnHeadersHeightSizeMode.AutoSize;\\
            \hspace*{1.5cm}this.dataGridA4.Location = new System.Drawing.Point(12, 49);\\
            \hspace*{1.5cm}this.dataGridA4.Name = "dataGridA4";\\
            \hspace*{1.5cm}this.dataGridA4.Size = new System.Drawing.Size(223, 227);\\
            \hspace*{1.5cm}this.dataGridA4.TabIndex = 0;\\
            \hspace*{1.5cm}//\\
            \hspace*{1.5cm}// textPower\\
            \hspace*{1.5cm}//\\
            \hspace*{1.5cm}this.textPower.Location = new System.Drawing.Point(311, 76);\\
            \hspace*{1.5cm}this.textPower.Name = "textPower";\\
            \hspace*{1.5cm}this.textPower.Size = new System.Drawing.Size(41, 22);\\
            \hspace*{1.5cm}this.textPower.TabIndex = 18;\\
            \hspace*{1.5cm}//\\
            \hspace*{1.5cm}// textlineA4\\
            \hspace*{1.5cm}//\\
            \hspace*{1.5cm}this.textlineA4.Location = new System.Drawing.Point(382, 24);\\
            \hspace*{1.5cm}this.textlineA4.Name = "textlineA4";\\
            \hspace*{1.5cm}this.textlineA4.Size = new System.Drawing.Size(41, 22);\\
            \hspace*{1.5cm}this.textlinieA4.TabIndex = 17;\\
            \hspace*{1.5cm}//\\
            \hspace*{1.5cm}// label18\\
            \hspace*{1.5cm}//\\
            \hspace*{1.5cm}this.label18.Anchor = ((System.Windows.Forms.AnchorStyles)(((System.Windows.\\
   Forms.AnchorStyles.Top $|$ System.Windows.Forms.AnchorStyles.Bottom) $|$ System.Windows.\\
   Forms.AnchorStyles.Left)));\\
            \hspace*{1.5cm}this.label18.AutoSize = true;\\
            \hspace*{1.5cm}this.label18.Location = new System.Drawing.Point(19, 76);\\
            \hspace*{1.5cm}this.label18.Name = "label18";\\
            \hspace*{1.5cm}this.label18.Size = new System.Drawing.Size(286, 16);\\
            \hspace*{1.5cm}this.label18.TabIndex = 16;\\
            \hspace*{1.5cm}this.label18.Text = "Introduce power of lifting for matrix A!";\\
            \hspace*{1.5cm}//\\
            \hspace*{1.5cm}// label19\\
            \hspace*{1.5cm}//\\
            \hspace*{1.5cm}this.label19.Anchor = ((System.Windows.Forms.AnchorStyles)(((System.Windows.\\
  Forms.AnchorStyles.Top $|$ System.Windows.Forms.AnchorStyles.Bottom) $|$ System.Windows.\\
  Forms.AnchorStyles.Left)));\\
            \hspace*{1.5cm}this.label19.AutoSize = true;\\
            \hspace*{1.5cm}this.label19.Location = new System.Drawing.Point(19, 27);\\
            \hspace*{1.5cm}this.label19.Name = "label19";\\
            \hspace*{1.5cm}this.label19.Size = new System.Drawing.Size(357, 16);\\
            \hspace*{1.5cm}this.label19.TabIndex = 15;\\
            \hspace*{1.5cm}this.label19.Text = "Introduce number of lines and columns for matrix A !";\\
            \hspace*{1.5cm}//\\
            \hspace*{1.5cm}// tabComputDet\\
            \hspace*{1.5cm}//\\
            \hspace*{1.5cm}this.tabComputDet.BackColor = System.Drawing.Color.FromArgb(((int)\\
 (((byte)(192)))), ((int)(((byte)(255)))),((int)(((byte)(192)))));\\
            \hspace*{1.5cm}this.tabComputDet.Controls.Add(this.pictureBox10);\\
            \hspace*{1.5cm}this.tabComputDet.Controls.Add(this.bt$_{-}$Reset$_{-}$values$_{-}$det);\\
            \hspace*{1.5cm}this.tabComputDet.Controls.Add(this.btComputingDeterminant);\\
            \hspace*{1.5cm}this.tabComputDet.Controls.Add(this.btGeneratingA5);\\
            \hspace*{1.5cm}this.tabComputDet.Controls.Add(this.label30);\\
            \hspace*{1.5cm}this.tabComputDet.Controls.Add(this.pictureBox8);\\
            \hspace*{1.5cm}this.tabComputDet.Controls.Add(this.pictureBox7);\\
            \hspace*{1.5cm}this.tabComputDet.Controls.Add(this.textDet$_{-}$minus);\\
            \hspace*{1.5cm}this.tabComputDet.Controls.Add(this.textDet$_{-}$plus);\\
            \hspace*{1.5cm}this.tabComputDet.Controls.Add(this.panel10);\\
            \hspace*{1.5cm}this.tabComputDet.Controls.Add(this.textlineA5);\\
            \hspace*{1.5cm}this.tabComputDet.Controls.Add(this.label22);\\
            \hspace*{1.5cm}this.tabComputDet.Location = new System.Drawing.Point(4, 25);\\
            \hspace*{1.5cm}this.tabComputDet.Name = "tabComputDet";\\
            \hspace*{1.5cm}this.tabComputDet.Padding = new System.Windows.Forms.Padding(3);\\
            \hspace*{1.5cm}this.tabComputDet.Size = new System.Drawing.Size(992, 667);\\
            \hspace*{1.5cm}this.tabComputDet.TabIndex = 4;\\
            \hspace*{1.5cm}this.tabComputDet.Text = "Computing of determinant";\\
            \hspace*{1.5cm}//\\
            \hspace*{1.5cm}// pictureBox10\\
            \hspace*{1.5cm}//\\
            \hspace*{1.5cm}this.pictureBox10.Image = ((System.Drawing.Image)(resources.GetObject\\
 ("pictureBox10.Image")));\\
            \hspace*{1.5cm}this.pictureBox10.Location = new System.Drawing.Point(53, 82);\\
            \hspace*{1.5cm}this.pictureBox10.Name = "pictureBox10";\\
            \hspace*{1.5cm}this.pictureBox10.Size = new System.Drawing.Size(56, 27);\\
            \hspace*{1.5cm}this.pictureBox10.TabIndex = 49;\\
            \hspace*{1.5cm}this.pictureBox10.TabStop = false;\\
            \hspace*{1.5cm}//\\
            \hspace*{1.5cm}// bt$_{-}$Reset$_{-}$values$_{-}$det\\
            \hspace*{1.5cm}//\\
            \hspace*{1.5cm}this.bt$_{-}$Reset$_{-}$values$_{-}$det.BackColor = System.Drawing.Color.Pink;\\
            \hspace*{1.5cm}this.bt$_{-}$Reset$_{-}$values$_{-}$det.ForeColor = System.Drawing.Color.Black;\\
            \hspace*{1.5cm}this.bt$_{-}$Reset$_{-}$values$_{-}$det.Location = new System.Drawing.Point(328, 369);\\
            \hspace*{1.5cm}this.bt$_{-}$Reset$_{-}$values$_{-}$det.Name = "bt$_{-}$Reset$_{-}$values$_{-}$det";\\
            \hspace*{1.5cm}this.bt$_{-}$Reset$_{-}$values$_{-}$det.Size = new System.Drawing.Size(125, 65);\\
            \hspace*{1.5cm}this.bt$_{-}$Reset$_{-}$values$_{-}$det.TabIndex = 48;\\
            \hspace*{1.5cm}this.bt$_{-}$Reset$_{-}$values$_{-}$det.Text = "Reset values";\\
            \hspace*{1.5cm}this.bt$_{-}$Reset$_{-}$values$_{-}$det.UseVisualStyleBackColor = false;\\
            \hspace*{1.5cm}this.bt$_{-}$Reset$_{-}$values$_{-}$det.Click += new System.EventHandler(this.\\
 bt$_{-}$Reset$_{-}$values$_{-}$det$_{-}$Click);\\
            \hspace*{1.5cm}//\\
            \hspace*{1.5cm}// btComputingDeterminant\\
            \hspace*{1.5cm}//\\
            \hspace*{1.5cm}this.btComputingDeterminant.BackColor = System.Drawing.Color.PowderBlue;\\
            \hspace*{1.5cm}this.btComputingDeterminant.ForeColor = System.Drawing.SystemColors.ControlText;\\
            \hspace*{1.5cm}this.btComputingDeterminant.Location = new System.Drawing.Point(328, 266);\\
            \hspace*{1.5cm}this.btComputingDeterminant.Name = "btComputingDeterminant";\\
            \hspace*{1.5cm}this.btComputingDeterminant.Size = new System.Drawing.Size(125, 65);\\
            \hspace*{1.5cm}this.btComputingDeterminant.TabIndex = 47;\\
            \hspace*{1.5cm}this.btComputingDeterminant.Text = "Computing determinant";\\
            \hspace*{1.5cm}this.btComputingDeterminant.UseVisualStyleBackColor = false;\\
            \hspace*{1.5cm}this.btComputingDeterminant.Click + = new System.EventHandler(this.\\
 btComputingDeterminant$_{-}$Click);\\
            \hspace*{1.5cm}//\\
            \hspace*{1.5cm}// btGeneratingA5\\
            \hspace*{1.5cm}//\\
            \hspace*{1.5cm}this.btGeneratingA5.BackColor = System.Drawing.Color.PowderBlue;\\
            \hspace*{1.5cm}this.btGeneratingA5.ForeColor = System.Drawing.SystemColors.ControlText;\\
            \hspace*{1.5cm}this.btGeneratingA5.Location = new System.Drawing.Point(328, 158);\\
            \hspace*{1.5cm}this.btGeneratingA5.Name = "btGeneratingA5";\\
            \hspace*{1.5cm}this.btGeneratingA5.Size = new System.Drawing.Size(125, 65);\\
            \hspace*{1.5cm}this.btGeneratingA5.TabIndex = 46;\\
            \hspace*{1.5cm}this.btGeneratingA5.Text = "Generating matrix";\\
            \hspace*{1.5cm}this.btGeneratingA5.UseVisualStyleBackColor = false;\\
            \hspace*{1.5cm}this.btGeneratingA5.Click + = new System.EventHandler(this.btGeneratingA5$_{-}$Click);\\
            \hspace*{1.5cm}//\\
            \hspace*{1.5cm}// label30\\
            \hspace*{1.5cm}//\\
            \hspace*{1.5cm}this.label30.AutoSize = true;\\
            \hspace*{1.5cm}this.label30.ForeColor = System.Drawing.Color.Red;\\
            \hspace*{1.5cm}this.label30.Location = new System.Drawing.Point(8, 88);\\
            \hspace*{1.5cm}this.label30.Name = "label30";\\
            \hspace*{1.5cm}this.label30.Size = new System.Drawing.Size(39, 16);\\
            \hspace*{1.5cm}this.label30.TabIndex = 35;\\
            \hspace*{1.5cm}this.label30.Text = "Rem:";\\
            \hspace*{1.5cm}//\\
            \hspace*{1.5cm}// pictureBox8\\
            \hspace*{1.5cm}//\\
            \hspace*{1.5cm}this.pictureBox8.BackColor = System.Drawing.Color.White;\\
            \hspace*{1.5cm}this.pictureBox8.Image = ((System.Drawing.Image)(resources.GetObject\\
 ("pictureBox8.Image")));\\
            \hspace*{1.5cm}this.pictureBox8.Location = new System.Drawing.Point(547, 271);\\
            \hspace*{1.5cm}this.pictureBox8.Name = "pictureBox8";\\
            \hspace*{1.5cm}this.pictureBox8.Size = new System.Drawing.Size(83, 31);\\
            \hspace*{1.5cm}this.pictureBox8.TabIndex = 28;\\
            \hspace*{1.5cm}this.pictureBox8.TabStop = false;\\
            \hspace*{1.5cm}//\\
            \hspace*{1.5cm}// pictureBox7\\
            \hspace*{1.5cm}//\\
            \hspace*{1.5cm}this.pictureBox7.BackColor = System.Drawing.Color.White;\\
            \hspace*{1.5cm}this.pictureBox7.Image = ((System.Drawing.Image)(resources.GetObject\\
 ("pictureBox7.Image")));\\
            \hspace*{1.5cm}this.pictureBox7.Location = new System.Drawing.Point(547, 197);\\
            \hspace*{1.5cm}this.pictureBox7.Name = "pictureBox7";\\
            \hspace*{1.5cm}this.pictureBox7.Size = new System.Drawing.Size(83, 33);\\
            \hspace*{1.5cm}this.pictureBox7.TabIndex = 27;\\
            \hspace*{1.5cm}this.pictureBox7.TabStop = false;\\
            \hspace*{1.5cm}//\\
            \hspace*{1.5cm}// textDet$_{-}$minus\\
            \hspace*{1.5cm}//\\
            \hspace*{1.5cm}this.textDet$_{-}$minus.Location = new System.Drawing.Point(636, 275);\\
            \hspace*{1.5cm}this.textDet$_{-}$minus.Name = "textDet$_{-}$minus";\\
            \hspace*{1.5cm}this.textDet$_{-}$minus.Size = new System.Drawing.Size(71, 22);\\
            \hspace*{1.5cm}this.textDet$_{-}$minus.TabIndex = 26;\\
            \hspace*{1.5cm}this.textDet$_{-}$minus.TextAlign = System.Windows.Forms.HorizontalAlignment.Right;\\
            \hspace*{1.5cm}//\\
            \hspace*{1.5cm}// textDet$_{-}$plus\\
            \hspace*{1.5cm}//\\
            \hspace*{1.5cm}this.textDet$_{-}$plus.Location = new System.Drawing.Point(636, 202);\\
            \hspace*{1.5cm}this.textDet$_{-}$plus.Name = "textDet$_{-}$plus";\\
            \hspace*{1.5cm}this.textDet$_{-}$plus.Size = new System.Drawing.Size(71, 22);\\
            \hspace*{1.5cm}this.textDet$_{-}$plus.TabIndex = 25;\\
            \hspace*{1.5cm}this.textDet$_{-}$plus.TextAlign = System.Windows.Forms.HorizontalAlignment.Right;\\
            \hspace*{1.5cm}//\\
            \hspace*{1.5cm}// panel10\\
            \hspace*{1.5cm}//\\
            \hspace*{1.5cm}this.panel10.Controls.Add(this.label23);\\
            \hspace*{1.5cm}this.panel10.Controls.Add(this.dataGridA5);\\
            \hspace*{1.5cm}this.panel10.Location = new System.Drawing.Point(25, 158);\\
            \hspace*{1.5cm}this.panel10.Name = "panel10";\\
            \hspace*{1.5cm}this.panel10.Size = new System.Drawing.Size(269, 297);\\
            \hspace*{1.5cm}this.panel10.TabIndex = 23;\\
            \hspace*{1.5cm}//\\
            \hspace*{1.5cm}// label23\\
            \hspace*{1.5cm}//\\
            \hspace*{1.5cm}this.label23.AutoSize = true;\\
            \hspace*{1.5cm}this.label23.ForeColor = System.Drawing.Color.MidnightBlue;\\
            \hspace*{1.5cm}this.label23.Location = new System.Drawing.Point(27, 16);\\
            \hspace*{1.5cm}this.label23.Name = "label23";\\
            \hspace*{1.5cm}this.label23.Size = new System.Drawing.Size(72, 16);\\
            \hspace*{1.5cm}this.label23.TabIndex = 10;\\
            \hspace*{1.5cm}this.label23.Text = "Matrix A";\\
            \hspace*{1.5cm}//\\
            \hspace*{1.5cm}// dataGridA5\\
            \hspace*{1.5cm}//\\
            \hspace*{1.5cm}this.dataGridA5.BackgroundColor = System.Drawing.Color.WhiteSmoke;\\
            \hspace*{1.5cm}this.dataGridA5.ColumnHeadersHeightSizeMode = System.Windows.Forms.\\
 DataGridViewColumnHeadersHeightSizeMode.AutoSize;\\
            \hspace*{1.5cm}this.dataGridA5.Location = new System.Drawing.Point(12, 49);\\
            \hspace*{1.5cm}this.dataGridA5.Name = "dataGridA5";\\
            \hspace*{1.5cm}this.dataGridA5.Size = new System.Drawing.Size(223, 227);\\
            \hspace*{1.5cm}this.dataGridA5.TabIndex = 0;\\
            \hspace*{1.5cm}//\\
            \hspace*{1.5cm}// textlineA5\\
            \hspace*{1.5cm}//\\
            \hspace*{1.5cm}this.textlineA5.Location = new System.Drawing.Point(371, 27);\\
            \hspace*{1.5cm}this.textlineA5.Name = "textlineA5";\\
            \hspace*{1.5cm}this.textlineA5.Size = new System.Drawing.Size(41, 22);\\
            \hspace*{1.5cm}this.textlineA5.TabIndex = 21;\\
            \hspace*{1.5cm}//\\
            \hspace*{1.5cm}// label22\\
            \hspace*{1.5cm}//\\
            \hspace*{1.5cm}this.label22.Anchor = ((System.Windows.Forms.AnchorStyles)(((System.Windows.\\
 Forms.AnchorStyles.Top $|$ System.Windows.Forms.AnchorStyles.Bottom)$|$ System.Windows.\\
 Forms.AnchorStyles.Left)));\\
            \hspace*{1.5cm}this.label22.AutoSize = true;\\
            \hspace*{1.5cm}this.label22.Location = new System.Drawing.Point(8, 30);\\
            \hspace*{1.5cm}this.label22.Name = "label22";\\
            \hspace*{1.5cm}this.label22.Size = new System.Drawing.Size(357, 16);\\
            \hspace*{1.5cm}this.label22.TabIndex = 20;\\
            \hspace*{1.5cm}this.label22.Text = "Introduce number of lines and columns of matrix A !";\\
            \hspace*{1.5cm}//\\
            \hspace*{1.5cm}// tabPage1\\
            \hspace*{1.5cm}//\\
            \hspace*{1.5cm}this.tabPage1.BackColor = System.Drawing.Color.FromArgb(((int)\\
 (((byte)(192)))), ((int)(((byte)(255)))),((int)(((byte)(192)))));\\
            \hspace*{1.5cm}this.tabPage1.Controls.Add(this.btResetSystem);\\
            \hspace*{1.5cm}this.tabPage1.Controls.Add(this.btComputSystem);\\
            \hspace*{1.5cm}this.tabPage1.Controls.Add(this.btGeneratingMatrices);\\
            \hspace*{1.5cm}this.tabPage1.Controls.Add(this.panel12);\\
            \hspace*{1.5cm}this.tabPage1.Controls.Add(this.panel9);\\
            \hspace*{1.5cm}this.tabPage1.Controls.Add(this.panel11);\\
            \hspace*{1.5cm}this.tabPage1.Controls.Add(this.textk);\\
            \hspace*{1.5cm}this.tabPage1.Controls.Add(this.textlineA6);\\
            \hspace*{1.5cm}this.tabPage1.Controls.Add(this.label39);\\
            \hspace*{1.5cm}this.tabPage1.Controls.Add(this.label40);\\
            \hspace*{1.5cm}this.tabPage1.Location = new System.Drawing.Point(4, 25);\\
            \hspace*{1.5cm}this.tabPage1.Name = "tabPage1";\\
            \hspace*{1.5cm}this.tabPage1.Padding = new System.Windows.Forms.Padding(3);\\
            \hspace*{1.5cm}this.tabPage1.Size = new System.Drawing.Size(992, 667);\\
            \hspace*{1.5cm}this.tabPage1.TabIndex = 5;\\
            \hspace*{1.5cm}this.tabPage1.Text = "Solving of system";\\
            \hspace*{1.5cm}//\\
            \hspace*{1.5cm}// panel12\\
            \hspace*{1.5cm}//\\
            \hspace*{1.5cm}this.panel12.Controls.Add(this.label41);\\
            \hspace*{1.5cm}this.panel12.Controls.Add(this.dataGridX0);\\
            \hspace*{1.5cm}this.panel12.Location = new System.Drawing.Point(293, 103);\\
            \hspace*{1.5cm}this.panel12.Name = "panel12";\\
            \hspace*{1.5cm}this.panel12.Size = new System.Drawing.Size(126, 282);\\
            \hspace*{1.5cm}this.panel12.TabIndex = 64;\\
            \hspace*{1.5cm}//\\
            \hspace*{1.5cm}// label41\\
            \hspace*{1.5cm}//\\
            \hspace*{1.5cm}this.label41.AutoSize = true;\\
            \hspace*{1.5cm}this.label41.ForeColor = System.Drawing.Color.MidnightBlue;\\
            \hspace*{1.5cm}this.label41.Location = new System.Drawing.Point(27, 16);\\
            \hspace*{1.5cm}this.label41.Name = "label41";\\
            \hspace*{1.5cm}this.label41.Size = new System.Drawing.Size(84, 16);\\
            \hspace*{1.5cm}this.label41.TabIndex = 10;\\
            \hspace*{1.5cm}this.label41.Text = "Matrix x(0)";\\
            \hspace*{1.5cm}//\\
            \hspace*{1.5cm}// dataGridX0\\
            \hspace*{1.5cm}//\\
            \hspace*{1.5cm}this.dataGridX0.BackgroundColor = System.Drawing.Color.WhiteSmoke;\\
            \hspace*{1.5cm}this.dataGridX0.ColumnHeadersHeightSizeMode = System.Windows.Forms.\\
 DataGridViewColumnHeadersHeightSizeMode.AutoSize;\\
            \hspace*{1.5cm}this.dataGridX0.Location = new System.Drawing.Point(30, 49);\\
            \hspace*{1.5cm}this.dataGridX0.Name = "dataGridX0";\\
            \hspace*{1.5cm}this.dataGridX0.Size = new System.Drawing.Size(67, 227);\\
            \hspace*{1.5cm}this.dataGridX0.TabIndex = 0;\\
            \hspace*{1.5cm}//\\
            \hspace*{1.5cm}// panel9\\
            \hspace*{1.5cm}//\\
            \hspace*{1.5cm}this.panel9.Controls.Add(this.label$_{-}$k);\\
            \hspace*{1.5cm}this.panel9.Controls.Add(this.label32);\\
            \hspace*{1.5cm}this.panel9.Controls.Add(this.label33);\\
            \hspace*{1.5cm}this.panel9.Controls.Add(this.label37);\\
            \hspace*{1.5cm}this.panel9.Controls.Add(this.dataGridSolutionXk);\\
            \hspace*{1.5cm}this.panel9.Location = new System.Drawing.Point(630, 103);\\
            \hspace*{1.5cm}this.panel9.Name = "panel9";\\
            \hspace*{1.5cm}this.panel9.Size = new System.Drawing.Size(257, 282);\\
            \hspace*{1.5cm}this.panel9.TabIndex = 62;\\
            \hspace*{1.5cm}//\\
            \hspace*{1.5cm}// label$_{-}k$\\
            \hspace*{1.5cm}//\\
            \hspace*{1.5cm}this.label$_{-}k$.ForeColor = System.Drawing.Color.Navy;\\
            \hspace*{1.5cm}this.label$_{-}k$.Location = new System.Drawing.Point(112, 25);\\
            \hspace*{1.5cm}this.label$_{-}k$.Name = "label$_{-}k$";\\
            \hspace*{1.5cm}this.label$_{-}k$.Size = new System.Drawing.Size(19, 18);\\
            \hspace*{1.5cm}this.label$_{-}k$.TabIndex = 20;\\
            \hspace*{1.5cm}this.label$_{-}k$.Text = "k";\\
            \hspace*{1.5cm}//\\
            \hspace*{1.5cm}// label32\\
            \hspace*{1.5cm}//\\
            \hspace*{1.5cm}this.label32.AutoSize = true;\\
            \hspace*{1.5cm}this.label32.Font = new System.Drawing.Font("Microsoft Sans Serif", 9.75F,\\
  System.Drawing.FontStyle.Regular, System. Drawing.GraphicsUnit.Point, ((byte)(0)));\\
            \hspace*{1.5cm}this.label32.ForeColor = System.Drawing.Color.Navy;\\
            \hspace*{1.5cm}this.label32.Location = new System.Drawing.Point(129, 25);\\
            \hspace*{1.5cm}this.label32.Name = "label32";\\
            \hspace*{1.5cm}this.label32.Size = new System.Drawing.Size(12, 16);\\
            \hspace*{1.5cm}this.label32.TabIndex = 19;\\
            \hspace*{1.5cm}this.label32.Text = ")";\\
            \hspace*{1.5cm}//\\
            \hspace*{1.5cm}// label33\\
            \hspace*{1.5cm}//\\
            \hspace*{1.5cm}this.label33.AutoSize = true;\\
            \hspace*{1.5cm}this.label33.ForeColor = System.Drawing.Color.Navy;\\
            \hspace*{1.5cm}this.label33.Location = new System.Drawing.Point(34, 25);\\
            \hspace*{1.5cm}this.label33.Name = "label33";\\
            \hspace*{1.5cm}this.label33.Size = new System.Drawing.Size(60, 16);\\
            \hspace*{1.5cm}this.label33.TabIndex = 18;\\
            \hspace*{1.5cm}this.label33.Text = "Matrix";\\
            \hspace*{1.5cm}//\\
            \hspace*{1.5cm}// label37\\
            \hspace*{1.5cm}//\\
            \hspace*{1.5cm}this.label37.AutoSize = true;\\
            \hspace*{1.5cm}this.label37.ForeColor = System.Drawing.Color.Navy;\\
            \hspace*{1.5cm}this.label37.Location = new System.Drawing.Point(100, 25);\\
            \hspace*{1.5cm}this.label37.Name = "label37";\\
            \hspace*{1.5cm}this.label37.Size = new System.Drawing.Size(18, 16);\\
            \hspace*{1.5cm}this.label37.TabIndex = 14;\\
            \hspace*{1.5cm}this.label37.Text = "x(";\\
            \hspace*{1.5cm}//\\
            \hspace*{1.5cm}// dataGridSolutionXk\\
            \hspace*{1.5cm}//\\
            \hspace*{1.5cm}this.dataGridSolutionXk.BackgroundColor = System.Drawing.Color.WhiteSmoke;\\
            \hspace*{1.5cm}this.dataGridSolutionXk.ColumnHeadersHeightSizeMode = System.Windows.Forms.\\
  DataGridViewColumnHeadersHeightSizeMode.AutoSize;\\
            \hspace*{1.5cm}this.dataGridSolutionXk.Location = new System.Drawing.Point(54, 49);\\
            \hspace*{1.5cm}this.dataGridSolutionXk.Name = "dataGridSolutionXk";\\
            \hspace*{1.5cm}this.dataGridSolutionXk.Size = new System.Drawing.Size(77, 227);\\
            \hspace*{1.5cm}this.dataGridSolutionXk.TabIndex = 8;\\
            \hspace*{1.5cm}//\\
            \hspace*{1.5cm}// panel11\\
            \hspace*{1.5cm}//\\
            \hspace*{1.5cm}this.panel11.Controls.Add(this.label38);\\
            \hspace*{1.5cm}this.panel11.Controls.Add(this.dataGridA6);\\
            \hspace*{1.5cm}this.panel11.Location = new System.Drawing.Point(18, 103);\\
            \hspace*{1.5cm}this.panel11.Name = "panel11";\\
            \hspace*{1.5cm}this.panel11.Size = new System.Drawing.Size(269, 282);\\
            \hspace*{1.5cm}this.panel11.TabIndex = 60;\\
            \hspace*{1.5cm}//\\
            \hspace*{1.5cm}// label38\\
            \hspace*{1.5cm}//\\
            \hspace*{1.5cm}this.label38.AutoSize = true;\\
            \hspace*{1.5cm}this.label38.ForeColor = System.Drawing.Color.MidnightBlue;\\
            \hspace*{1.5cm}this.label38.Location = new System.Drawing.Point(27, 16);\\
            \hspace*{1.5cm}this.label38.Name = "label38";\\
            \hspace*{1.5cm}this.label38.Size = new System.Drawing.Size(72, 16);\\
            \hspace*{1.5cm}this.label38.TabIndex = 10;\\
            \hspace*{1.5cm}this.label38.Text = "Matrix A";\\
            \hspace*{1.5cm}//\\
            \hspace*{1.5cm}// dataGridA6\\
            \hspace*{1.5cm}//\\
            \hspace*{1.5cm}this.dataGridA6.BackgroundColor = System.Drawing.Color.WhiteSmoke;\\
            \hspace*{1.5cm}this.dataGridA6.ColumnHeadersHeightSizeMode = System.Windows.Forms.\\
  DataGridViewColumnHeadersHeightSizeMode.AutoSize;\\
            \hspace*{1.5cm}this.dataGridA6.Location = new System.Drawing.Point(23, 49);\\
            \hspace*{1.5cm}this.dataGridA6.Name = "dataGridA6";\\
            \hspace*{1.5cm}this.dataGridA6.Size = new System.Drawing.Size(223, 227);\\
            \hspace*{1.5cm}this.dataGridA6.TabIndex = 0;\\
            \hspace*{1.5cm}//\\
            \hspace*{1.5cm}// textk\\
            \hspace*{1.5cm}//\\
            \hspace*{1.5cm}this.textk.Location = new System.Drawing.Point(176, 66);\\
            \hspace*{1.5cm}this.textk.Name = "textk";\\
            \hspace*{1.5cm}this.textk.Size = new System.Drawing.Size(41, 22);\\
            \hspace*{1.5cm}this.textk.TabIndex = 58;\\
            \hspace*{1.5cm}//\\
            \hspace*{1.5cm}// textlineA6\\
            \hspace*{1.5cm}//\\
            \hspace*{1.5cm}this.textlineA6.Location = new System.Drawing.Point(378, 26);\\
            \hspace*{1.5cm}this.textlineA6.Name = "textlineA6";\\
            \hspace*{1.5cm}this.textlineA6.Size = new System.Drawing.Size(41, 22);\\
            \hspace*{1.5cm}this.textlineA6.TabIndex = 57;\\
            \hspace*{1.5cm}//\\
            \hspace*{1.5cm}// label39\\
            \hspace*{1.5cm}//\\
            \hspace*{1.5cm}this.label39.Anchor = ((System.Windows.Forms.AnchorStyles)(((System.Windows.\\
  Forms.AnchorStyles.Top $|$ System.Windows.Forms.AnchorStyles.Bottom)$|$ System.Windows.\\
  Forms.AnchorStyles.Left)));\\
            \hspace*{1.5cm}this.label39.AutoSize = true;\\
            \hspace*{1.5cm}this.label39.Location = new System.Drawing.Point(15, 66);\\
            \hspace*{1.5cm}this.label39.Name = "label39";\\
            \hspace*{1.5cm}this.label39.Size = new System.Drawing.Size(155, 16);\\
            \hspace*{1.5cm}this.label39.TabIndex = 56;\\
            \hspace*{1.5cm}this.label39.Text = "Introduce value of k!";\\
            \hspace*{1.5cm}//\\
            \hspace*{1.5cm}// label40\\
            \hspace*{1.5cm}//\\
            \hspace*{1.5cm}this.label40.Anchor = ((System.Windows.Forms.AnchorStyles)(((System.Windows.\\
  Forms.AnchorStyles.Top $|$ System.Windows.Forms.AnchorStyles.Bottom)$|$ System.Windows.\\
  Forms.AnchorStyles.Left)));\\
            \hspace*{1.5cm}this.label40.AutoSize = true;\\
            \hspace*{1.5cm}this.label40.Location = new System.Drawing.Point(15, 29);\\
            \hspace*{1.5cm}this.label40.Name = "label40";\\
            \hspace*{1.5cm}this.label40.Size = new System.Drawing.Size(357, 16);\\
            \hspace*{1.5cm}this.label40.TabIndex = 55;\\
            \hspace*{1.5cm}this.label40.Text = "Introduce number of lines and columns of matrix A !";\\
            \hspace*{1.5cm}//\\
            \hspace*{1.5cm}// btResetSystem\\
            \hspace*{1.5cm}//\\
            \hspace*{1.5cm}this.btResetSystem.BackColor = System.Drawing.Color.Pink;\\
            \hspace*{1.5cm}this.btResetSystem.ForeColor = System.Drawing.Color.Black;\\
            \hspace*{1.5cm}this.btResetSystem.Location = new System.Drawing.Point(461, 314);\\
            \hspace*{1.5cm}this.btResetSystem.Name = "btResetSystem";\\
            \hspace*{1.5cm}this.btResetSystem.Size = new System.Drawing.Size(125,65);\\
            \hspace*{1.5cm}this.btResetSystem.TabIndex = 67;\\
            \hspace*{1.5cm}this.btResetSystem.Text = "Reset values";\\
            \hspace*{1.5cm}this.btResetSystem.UseVisualStyleBackColor = false;\\
            \hspace*{1.5cm}this.btResetSystem.Click + = new System.EventHandler(this.btResetSystem$_{-}$Click);\\
            \hspace*{1.5cm}//\\
            \hspace*{1.5cm}// btComputSystem\\
            \hspace*{1.5cm}//\\
            \hspace*{1.5cm}this.btComputSystem.BackColor = System.Drawing.Color.PowderBlue;\\
            \hspace*{1.5cm}this.btComputSystem.ForeColor = System.Drawing.SystemColors.ControlText;\\
            \hspace*{1.5cm}this.btComputSystem.Location = new System.Drawing.Point(461, 211);\\
            \hspace*{1.5cm}this.btComputSystem.Name = "btComputSystem";\\
            \hspace*{1.5cm}this.btComputSystem.Size = new System.Drawing.Size(125, 65);\\
            \hspace*{1.5cm}this.btComputSystem.TabIndex = 66;\\
            \hspace*{1.5cm}this.btComputSystem.Text = "Computing matrix x(k)";\\
            \hspace*{1.5cm}this.btComputSystem.UseVisualStyleBackColor = false;\\
            \hspace*{1.5cm}this.btComputSystem.Click + = new System.EventHandler(this.btComputSystem$_{-}$Click);\\
            \hspace*{1.5cm}//\\
            \hspace*{1.5cm}// btGeneratingMatrices\\
            \hspace*{1.5cm}//\\
            \hspace*{1.5cm}this.btGeneratingMatrices.BackColor = System.Drawing.Color.PowderBlue;\\
            \hspace*{1.5cm}this.btGeneratingMatrices.ForeColor = System.Drawing.SystemColors.ControlText;\\
            \hspace*{1.5cm}this.btGeneratingMatrices.Location = new System.Drawing.Point(461, 103);\\
            \hspace*{1.5cm}this.btGeneratingMatrices.Name = "btGeneratingMatrices";\\
            \hspace*{1.5cm}this.btGeneratingMatrices.Size = new System.Drawing.Size(125, 65);\\
            \hspace*{1.5cm}this.btGeneratingMatrices.TabIndex = 65;\\
            \hspace*{1.5cm}this.btGeneratingMatrices.Text = "Generating matrix";\\
            \hspace*{1.5cm}this.btGeneratingMatrices.UseVisualStyleBackColor = false;\\
            \hspace*{1.5cm}this.btGeneratingMatrices.Click + = new System.EventHandler(this.\\
  btGeneratingMatrices$_{-}$Click);\\
            \hspace*{1.5cm}//\\
            \hspace*{1.5cm}// Form1\\
            \hspace*{1.5cm}//\\
            \hspace*{1.5cm}this.AutoScaleDimensions = new System.Drawing.SizeF(8F, 16F);\\
            \hspace*{1.5cm}this.AutoScaleMode = System.Windows.Forms.AutoScaleMode.Font;\\
            \hspace*{1.5cm}this.BackColor = System.Drawing.Color.FromArgb(((int)\\
  (((byte)(192)))), ((int)(((byte)(255)))),((int)(((byte)(192)))));\\
            \hspace*{1.5cm}this.ClientSize = new System.Drawing.Size(1024, 742);\\
            \hspace*{1.5cm}this.Controls.Add(this.tabControl1);\\
            \hspace*{1.5cm}this.Font = new System.Drawing.Font("Microsoft Sans Serif", 9.75F, System.Drawing.\\
  FontStyle.Regular, System.Drawing.GraphicsUnit.Point,((byte)(0)));\\
            \hspace*{1.5cm}this.Margin = new System.Windows.Forms.Padding(4);\\
            \hspace*{1.5cm}this.Name = "Form1";\\
            \hspace*{1.5cm}this.Text = "Operations with  matrices in min-plus algebra";\\
            \hspace*{1.5cm}this.tabControl1.ResumeLayout(false);\\
            \hspace*{1.5cm}this.tabAddition.ResumeLayout(false);\\
            \hspace*{1.5cm}this.tabAddition.PerformLayout();\\[0.2cm]
 ((System.ComponentModel.ISupportInitialize)(this.pictureBox1)).EndInit();\\
            \hspace*{1.5cm}this.panel1.ResumeLayout(false);\\
 ((System.ComponentModel.ISupportInitialize)(this.pictureBox2)).EndInit();\\
 ((System.ComponentModel.ISupportInitialize)(this.dataGridResAddition)).EndInit();\\
            \hspace*{1.5cm}this.panel2.ResumeLayout(false);\\
            \hspace*{1.5cm}this.panel2.PerformLayout();\\
 ((System.ComponentModel.ISupportInitialize)(this.dataGridB)).EndInit();\\
 ((System.ComponentModel.ISupportInitialize)(this.dataGridA)).EndInit();\\
            \hspace*{1.5cm}this.tabMultiplication.ResumeLayout(false);\\
            \hspace*{1.5cm}this.tabMultiplication.PerformLayout();\\
 ((System.ComponentModel.ISupportInitialize)(this.pictureBox3)).EndInit();\\
            \hspace*{1.5cm}this.panel3.ResumeLayout(false);\\
 ((System.ComponentModel.ISupportInitialize)(this.pictureBox4)).EndInit();\\
 ((System.ComponentModel.ISupportInitialize)(this.dataGridProduct)).EndInit();\\
            \hspace*{1.5cm}this.panel4.ResumeLayout(false);\\
            \hspace*{1.5cm}this.panel4.PerformLayout();\\
 ((System.ComponentModel.ISupportInitialize)(this.dataGridB2)).EndInit();\\
 ((System.ComponentModel.ISupportInitialize)(this.dataGridA2)).EndInit();\\
            \hspace*{1.5cm}this.tabMultiplication$_{-}$scalar.ResumeLayout(false);\\
            \hspace*{1.5cm}this.tabMultiplication$_{-}$scalar.PerformLayout();\\
 ((System.ComponentModel.ISupportInitialize)(this.pictureBox6)).EndInit();\\
            \hspace*{1.5cm}this.panel6.ResumeLayout(false);\\
 ((System.ComponentModel.ISupportInitialize)(this.pictureBox5)).EndInit();\\
 ((System.ComponentModel.ISupportInitialize)(this.dataGridProductScalar)).EndInit();\\
            \hspace*{1.5cm}this.panel5.ResumeLayout(false);\\
            \hspace*{1.5cm}this.panel5.PerformLayout();\\
 ((System.ComponentModel.ISupportInitialize)(this.dataGridA3)).EndInit();\\
            \hspace*{1.5cm}this.tabLifting$_{-}$at$_{-}$power.ResumeLayout(false);\\
            \hspace*{1.5cm}this.tabLifting$_{-}$at$_{-}$power.PerformLayout();\\
 ((System.ComponentModel.ISupportInitialize)(this.pictureBox9)).EndInit();\\
            \hspace*{1.5cm}this.panel7.ResumeLayout(false);\\
            \hspace*{1.5cm}this.panel7.PerformLayout();\\
 ((System.ComponentModel.ISupportInitialize)(this.dataGridMatrix$_{-}$at$_{-}$power$_{-}n$)).EndInit();\\
            \hspace*{1.5cm}this.panel8.ResumeLayout(false);\\
            \hspace*{1.5cm}this.panel8.PerformLayout();\\
 ((System.ComponentModel.ISupportInitialize)(this.dataGridA4)).EndInit();\\
            \hspace*{1.5cm}this.tabComputDet.ResumeLayout(false);\\
            \hspace*{1.5cm}this.tabComputDet.PerformLayout();\\
 ((System.ComponentModel.ISupportInitialize)(this.pictureBox10)).EndInit();\\
 ((System.ComponentModel.ISupportInitialize)(this.pictureBox8)).EndInit();\\
 ((System.ComponentModel.ISupportInitialize)(this.pictureBox7)).EndInit();\\
            \hspace*{1.5cm}this.panel10.ResumeLayout(false);\\
            \hspace*{1.5cm}this.panel10.PerformLayout();\\
 ((System.ComponentModel.ISupportInitialize)(this.dataGridA5)).EndInit();\\
            \hspace*{1.5cm}this.tabPage1.ResumeLayout(false);\\
            \hspace*{1.5cm}this.tabPage1.PerformLayout();\\
            \hspace*{1.5cm}this.panel12.ResumeLayout(false);\\
            \hspace*{1.5cm}this.panel12.PerformLayout();\\
 ((System.ComponentModel.ISupportInitialize)(this.dataGridX0)).EndInit();\\
            \hspace*{1.5cm}this.panel9.ResumeLayout(false);\\
            \hspace*{1.5cm}this.panel9.PerformLayout();\\
 ((System.ComponentModel.ISupportInitialize)(this.dataGridSolutieXk)).EndInit();\\
            \hspace*{1.5cm}this.panel11.ResumeLayout(false);\\
            \hspace*{1.5cm}this.panel11.PerformLayout();\\
 ((System.ComponentModel.ISupportInitialize)(this.dataGridA6)).EndInit();\\
            \hspace*{1.5cm}this.ResumeLayout(false);\\[0.1cm]
       \hspace*{1cm}$\}$\\[0.1cm]
        \hspace*{1cm}$\#$endregion\\[0.1cm]
        \hspace*{1cm}private System.Windows.Forms.TabControl tabControl1;\\
        \hspace*{1cm}private System.Windows.Forms.TabPage tabAddition;\\
        \hspace*{1cm}private System.Windows.Forms.TabPage tabMultiplication;\\
        \hspace*{1cm}private System.Windows.Forms.TabPage tabMultiplication$_{-}$scalar;\\
        \hspace*{1cm}private System.Windows.Forms.Panel panel2;\\
        \hspace*{1cm}private System.Windows.Forms.DataGridView dataGridB;\\
        \hspace*{1cm}private System.Windows.Forms.DataGridView dataGridA;\\
        \hspace*{1cm}private System.Windows.Forms.TextBox textColumn;\\
        \hspace*{1cm}private System.Windows.Forms.TextBox textLine;\\
        \hspace*{1cm}private System.Windows.Forms.Label label3;\\
        \hspace*{1cm}private System.Windows.Forms.Label label2;\\
        \hspace*{1cm}private System.Windows.Forms.Label label4;\\
        \hspace*{1cm}private System.Windows.Forms.Label label1;\\
        \hspace*{1cm}private System.Windows.Forms.Panel panel1;\\
        \hspace*{1cm}private System.Windows.Forms.DataGridView dataGridResAddition;\\
        \hspace*{1cm}private System.Windows.Forms.PictureBox pictureBox2;\\
        \hspace*{1cm}private System.Windows.Forms.Label label13;\\
        \hspace*{1cm}private System.Windows.Forms.TextBox textcolumnB2;\\
        \hspace*{1cm}private System.Windows.Forms.TextBox textLineB2;\\
        \hspace*{1cm}private System.Windows.Forms.Label label11;\\
        \hspace*{1cm}private System.Windows.Forms.Label label12;\\
        \hspace*{1cm}private System.Windows.Forms.Panel panel3;\\
        \hspace*{1cm}private System.Windows.Forms.DataGridView dataGridProduct;\\
        \hspace*{1cm}private System.Windows.Forms.Panel panel4;\\
        \hspace*{1cm}private System.Windows.Forms.Label label7;\\
        \hspace*{1cm}private System.Windows.Forms.Label label8;\\
        \hspace*{1cm}private System.Windows.Forms.DataGridView dataGridB2;\\
        \hspace*{1cm}private System.Windows.Forms.DataGridView dataGridA2;\\
        \hspace*{1cm}private System.Windows.Forms.TextBox textColumnA2;\\
        \hspace*{1cm}private System.Windows.Forms.TextBox textLineA2;\\
        \hspace*{1cm}private System.Windows.Forms.Label label9;\\
        \hspace*{1cm}private System.Windows.Forms.Label label10;\\
        \hspace*{1cm}private System.Windows.Forms.PictureBox pictureBox4;\\
        \hspace*{1cm}private System.Windows.Forms.Panel panel5;\\
        \hspace*{1cm}private System.Windows.Forms.Label label17;\\
        \hspace*{1cm}private System.Windows.Forms.DataGridView dataGridA3;\\
        \hspace*{1cm}private System.Windows.Forms.TextBox textcolumnA3;\\
        \hspace*{1cm}private System.Windows.Forms.TextBox textLineA3;\\
        \hspace*{1cm}private System.Windows.Forms.Label label14;\\
        \hspace*{1cm}private System.Windows.Forms.Label label15;\\
        \hspace*{1cm}private System.Windows.Forms.TextBox textScalar;\\
        \hspace*{1cm}private System.Windows.Forms.Label label16;\\
        \hspace*{1cm}private System.Windows.Forms.Panel panel6;\\
        \hspace*{1cm}private System.Windows.Forms.PictureBox pictureBox5;\\
        \hspace*{1cm}private System.Windows.Forms.DataGridView dataGridProductScalar;\\
        \hspace*{1cm}private System.Windows.Forms.TabPage tabLifting$_{-}$at$_{-}$power;\\
        \hspace*{1cm}private System.Windows.Forms.Panel panel7;\\
        \hspace*{1cm}private System.Windows.Forms.DataGridView dataGridMatrix$_{-}$at$_{-}$power$_{-}n$;\\
        \hspace*{1cm}private System.Windows.Forms.Panel panel8;\\
        \hspace*{1cm}private System.Windows.Forms.Label label20;\\
        \hspace*{1cm}private System.Windows.Forms.DataGridView dataGridA4;\\
        \hspace*{1cm}private System.Windows.Forms.TextBox textPower;\\
        \hspace*{1cm}private System.Windows.Forms.TextBox textlineA4;\\
        \hspace*{1cm}private System.Windows.Forms.Label label18;\\
        \hspace*{1cm}private System.Windows.Forms.Label label19;\\
        \hspace*{1cm}private System.Windows.Forms.Label label21;\\
        \hspace*{1cm}private System.Windows.Forms.TabPage tabComputDet;\\
        \hspace*{1cm}private System.Windows.Forms.TextBox textlineA5;\\
        \hspace*{1cm}private System.Windows.Forms.Label label22;\\
        \hspace*{1cm}private System.Windows.Forms.Label label25;\\
        \hspace*{1cm}private System.Windows.Forms.Label label26;\\
        \hspace*{1cm}private System.Windows.Forms.Label labelPower;\\
        \hspace*{1cm}private System.Windows.Forms.Label label24;\\
        \hspace*{1cm}private System.Windows.Forms.Panel panel10;\\
        \hspace*{1cm}private System.Windows.Forms.Label label23;\\
        \hspace*{1cm}private System.Windows.Forms.DataGridView dataGridA5;\\
        \hspace*{1cm}private System.Windows.Forms.PictureBox pictureBox8;\\
        \hspace*{1cm}private System.Windows.Forms.PictureBox pictureBox7;\\
        \hspace*{1cm}private System.Windows.Forms.TextBox textDet$_{-}$minus;\\
        \hspace*{1cm}private System.Windows.Forms.TextBox textDet$_{-}$plus;\\
        \hspace*{1cm}private System.Windows.Forms.Label label27;\\
        \hspace*{1cm}private System.Windows.Forms.Label label6;\\
        \hspace*{1cm}private System.Windows.Forms.PictureBox pictureBox1;\\
        \hspace*{1cm}private System.Windows.Forms.Label label5;\\
        \hspace*{1cm}private System.Windows.Forms.Label label28;\\
        \hspace*{1cm}private System.Windows.Forms.Label label29;\\
        \hspace*{1cm}private System.Windows.Forms.Label label30;\\
        \hspace*{1cm}private System.Windows.Forms.Button btGenerating;\\
        \hspace*{1cm}private System.Windows.Forms.Button btReset;\\
        \hspace*{1cm}private System.Windows.Forms.Button btComputSum;\\
        \hspace*{1cm}private System.Windows.Forms.Button btComputProduct;\\
        \hspace*{1cm}private System.Windows.Forms.Button btGenerating2;\\
        \hspace*{1cm}private System.Windows.Forms.Button btResetMultiplication;\\
        \hspace*{1cm}private System.Windows.Forms.Button btReset$_{-}$Product$_{-}$Scalar;\\
        \hspace*{1cm}private System.Windows.Forms.Button btProductScalar;\\
        \hspace*{1cm}private System.Windows.Forms.Button butonsGenerating3;\\
        \hspace*{1cm}private System.Windows.Forms.Button btReset$_{-}$lifting$_{-}$at$_{-}$power;\\
        \hspace*{1cm}private System.Windows.Forms.Button btComputPower;\\
        \hspace*{1cm}private System.Windows.Forms.Button btGeneratingA4;\\
        \hspace*{1cm}private System.Windows.Forms.Button bt$_{-}$Reset$_{-}$values$_{-}$det;\\
        \hspace*{1cm}private System.Windows.Forms.Button btComputingDeterminant;\\
        \hspace*{1cm}private System.Windows.Forms.Button btGeneratingA5;\\
        \hspace*{1cm}private System.Windows.Forms.PictureBox pictureBox3;\\
        \hspace*{1cm}private System.Windows.Forms.PictureBox pictureBox6;\\
        \hspace*{1cm}private System.Windows.Forms.PictureBox pictureBox9;\\
        \hspace*{1cm}private System.Windows.Forms.PictureBox pictureBox10;\\
        \hspace*{1cm}private System.Windows.Forms.TabPage tabPage1;\\
        \hspace*{1cm}private System.Windows.Forms.Panel panel12;\\
        \hspace*{1cm}private System.Windows.Forms.Label label41;\\
        \hspace*{1cm}private System.Windows.Forms.DataGridView dataGridX0;\\
        \hspace*{1cm}private System.Windows.Forms.Panel panel9;\\
        \hspace*{1cm}private System.Windows.Forms.Label label$_{-}k$;\\
        \hspace*{1cm}private System.Windows.Forms.Label label32;\\
        \hspace*{1cm}private System.Windows.Forms.Label label33;\\
        \hspace*{1cm}private System.Windows.Forms.Label label37;\\
        \hspace*{1cm}private System.Windows.Forms.DataGridView dataGridSolutionXk;\\
        \hspace*{1cm}private System.Windows.Forms.Panel panel11;\\
        \hspace*{1cm}private System.Windows.Forms.Label label38;\\
        \hspace*{1cm}private System.Windows.Forms.DataGridView dataGridA6;\\
        \hspace*{1cm}private System.Windows.Forms.TextBox textk;\\
        \hspace*{1cm}private System.Windows.Forms.TextBox textlineA6;\\
        \hspace*{1cm}private System.Windows.Forms.Label label39;\\
        \hspace*{1cm}private System.Windows.Forms.Label label40;\\
        \hspace*{1cm}private System.Windows.Forms.Button btResetSystem;\\
        \hspace*{1cm}private System.Windows.Forms.Button btComputSystem;\\
        \hspace*{1cm}private System.Windows.Forms.Button btGeneratingMatrices;\\
 \hspace*{0.1cm}$\}$\\
$\}$\\

$\bullet~~~~~$ {\it The principal program is constituted from the following lines}.\\[0.1cm]
using System;\\
using System.Collections.Generic;\\
using System.ComponentModel;\\
using System.Data;\\
using System.Drawing;\\
using System.Text;\\
using System.Windows.Forms;\\[0.1cm]
namespace Operations$_{-}$ with$_{-}$matrices\\
$\{$\\
public partial class Form1 : Form\\
$\{$\\
\hspace*{0.5cm}public Form1()\\
\hspace*{0.5cm}$\{$\\
\hspace*{1cm}InitializeComponent();\\
\hspace*{0.5cm}$\}$\\[0.1cm]
\hspace*{0.5cm}public int[,] A5 = new int[50, 50];\\[0.1cm]
\hspace*{0.5cm}private void initMatrixA()\\
$\{$\\[0.1cm]
\hspace*{1cm}// Create an unbound DataGridView by declaring a column
count.\\[0.1cm]
\hspace*{0.5cm}int column = 0;\\
\hspace*{0.5cm}column = Convert.ToInt16(textColumn.Text);\\
\hspace*{0.5cm}dataGridA.ColumnCount = column;\\
\hspace*{0.5cm}dataGridA.AllowUserToOrderColumns = false;\\
\hspace*{0.5cm}dataGridA.AllowUserToAddRows = false;\\
\hspace*{0.5cm}dataGridA.Enabled = true;\\
\hspace*{0.5cm}dataGridA.AutoSizeRowsMode = DataGridViewAutoSizeRowsMode.\\
DisplayedCellsExceptHeaders;\\
\hspace*{0.5cm}dataGridA.ColumnHeadersBorderStyle =DataGridViewHeaderBorderStyle.Raised;\\
\hspace*{0.5cm}dataGridA.CellBorderStyle =DataGridViewCellBorderStyle.Single;\\
\hspace*{0.5cm}dataGridA.GridColor = Color.DodgerBlue;\\
\hspace*{0.5cm}dataGridA.ColumnHeadersVisible = false;\\
\hspace*{0.5cm}dataGridA.RowHeadersVisible = false;\\
\hspace*{0.5cm}dataGridA.BackgroundColor = Color.WhiteSmoke;\\
\hspace*{0.5cm}dataGridA.BorderStyle = BorderStyle.None;\\
\hspace*{0.5cm}dataGridA.AllowUserToResizeColumns = true;\\
\hspace*{0.5cm}DataGridViewCellStyle style = new DataGridViewCellStyle();\\
\hspace*{0.5cm}style.Format = "N0";\\
\hspace*{0.5cm}style.NullValue = null;\\[0.2cm]
\hspace*{1cm}//dataGridA.DefaultCellStyle = style;\\[0.1cm]
\hspace*{1cm}//dataGridA.DefaultCellStyle = Font(12.75F);\\[0.1cm]
\hspace*{1cm}//Set the column header style.\\[0.1cm]
\hspace*{0.5cm}DataGridViewCellStyle columnHeaderStyle = new DataGridViewCellStyle();\\
\hspace*{0.5cm}columnHeaderStyle.BackColor = Color.AntiqueWhite;\\
\hspace*{0.5cm}columnHeaderStyle.Alignment = DataGridViewContentAlignment.MiddleCenter;\\
\hspace*{0.5cm}dataGridA.ColumnHeadersDefaultCellStyle = columnHeaderStyle;\\[0.1cm]
\hspace*{0.5cm}DataGridViewCell CellR = new DataGridViewTextBoxCell();\\
\hspace*{0.5cm}CellR.Style.Alignment = DataGridViewContentAlignment.MiddleRight;\\[0.2cm]
\hspace*{0.5cm}DataGridViewCell CellL = new DataGridViewTextBoxCell();\\
\hspace*{0.5cm}CellL.Style.Alignment = DataGridViewContentAlignment.MiddleLeft;\\[0.1cm]
\hspace*{0.5cm}int line = 0;\\
\hspace*{0.5cm}line = Convert.ToInt16(textLine.Text);\\
\hspace*{0.5cm}dataGridA.RowCount = line;\\[0.1cm]
\hspace*{0.5cm}for (int i = 0; i $<$ column; i++)\\
\hspace*{0.5cm}$\{$\\
\hspace*{1cm}dataGridA.Columns[i].Name = "C" + (i + 1);\\
\hspace*{1cm}dataGridA.Columns[i].CellTemplate = CellR;\\
\hspace*{1cm}dataGridA.Columns[i].Width = 30;\\
\hspace*{1cm}dataGridA.Columns[i].DefaultCellStyle = style;\\
\hspace*{0.5cm}$\}$\\[0.1cm]
\hspace*{0.3cm}$\}$\\
\hspace*{0.3cm}private void initMatrixB()\\
\hspace*{0.3cm}$\{$\\[0.2cm]
\hspace*{0.5cm}// Create an unbound DataGridView by declaring a column count.\\[0.1cm]
\hspace*{0.5cm}int column = 0;\\
\hspace*{0.5cm} column = Convert.ToInt16(textColumn.Text);\\
\hspace*{0.5cm}dataGridB.ColumnCount = column;\\
\hspace*{0.5cm}dataGridB.AllowUserToOrderColumns = false;\\
\hspace*{0.5cm}dataGridB.AllowUserToAddRows = false;\\
\hspace*{0.5cm}dataGridB.Enabled = true;\\
\hspace*{0.5cm}dataGridB.AutoSizeRowsMode = DataGridViewAutoSizeRowsMode.\\
DisplayedCellsExceptHeaders;\\
\hspace*{0.5cm}dataGridB.ColumnHeadersBorderStyle = DataGridViewHeaderBorderStyle.Raised;\\
\hspace*{0.5cm}dataGridB.CellBorderStyle = DataGridViewCellBorderStyle.Single;\\
\hspace*{0.5cm}dataGridB.GridColor = Color.DodgerBlue;\\
\hspace*{0.5cm}dataGridB.ColumnHeadersVisible = false;\\
\hspace*{0.5cm}dataGridB.RowHeadersVisible = false;\\
\hspace*{0.5cm}dataGridB.BackgroundColor = Color.WhiteSmoke;\\
\hspace*{0.5cm}dataGridB.BorderStyle = BorderStyle.None;\\
\hspace*{0.5cm}dataGridB.AllowUserToResizeColumns = true;\\[0.2cm]
\hspace*{1cm}// dataGridB.DefaultCellStyle = Font(12.75F);\\
\hspace*{1cm}// Set the column header style.\\[0.2cm]
\hspace*{0.5cm}DataGridViewCellStyle columnHeaderStyle = new DataGridViewCellStyle();\\
\hspace*{0.5cm}columnHeaderStyle.BackColor = Color.AntiqueWhite;\\
\hspace*{0.5cm}columnHeaderStyle.Alignment = DataGridViewContentAlignment.MiddleCenter;\\
\hspace*{0.5cm}dataGridB.ColumnHeadersDefaultCellStyle = columnHeaderStyle;\\[0.2cm]
\hspace*{0.5cm}DataGridViewCell CellR = new DataGridViewTextBoxCell();\\
\hspace*{0.5cm}CellR.Style.Alignment = DataGridViewContentAlignment.MiddleRight;\\[0.1cm]
\hspace*{0.5cm}DataGridViewCell CellL = new DataGridViewTextBoxCell();\\
\hspace*{0.5cm}CellL.Style.Alignment = DataGridViewContentAlignment.MiddleLeft;\\[0.1cm]
\hspace*{0.5cm}int line = 0;\\
\hspace*{0.5cm}line = Convert.ToInt16(textLine.Text);\\
\hspace*{0.5cm}dataGridB.RowCount = line;\\[0.1cm]
\hspace*{0.5cm}for (int i = 0; i $<$ column; i++)\\
\hspace*{0.5cm}$\{$\\
\hspace*{1cm}dataGridB.Columns[i].Name = "C" + (i + 1);\\
\hspace*{1cm}dataGridB.Columns[i].CellTemplate = CellR;\\
\hspace*{0.5cm}dataGridB.Columns[i].Width = 30;\\
\hspace*{0.5cm}$\}$\\[0.1cm]
\hspace*{0.3cm}$\}$\\
\hspace*{0.3cm}private void initMatrixA2 ()\\
\hspace*{0.3cm}$\{$\\[0.2cm]
\hspace*{0.5cm}// Create an unbound DataGridView by declaring a column count.\\[0.2cm]
\hspace*{0.5cm}int column = 0;\\
\hspace*{0.5cm}column = Convert.ToInt16(textColumnA2.Text);\\
\hspace*{0.5cm}dataGridA2.ColumnCount = column;\\
\hspace*{0.5cm}dataGridA2.AllowUserToOrderColumns = false;\\
\hspace*{0.5cm}dataGridA2.AllowUserToAddRows = false;\\
\hspace*{0.5cm}dataGridA2.Enabled = true;\\
\hspace*{0.5cm}dataGridA2.AutoSizeRowsMode = DataGridViewAutoSizeRowsMode.\\
DisplayedCellsExceptHeaders;\\
\hspace*{0.5cm}dataGridA2.ColumnHeadersBorderStyle = DataGridViewHeaderBorderStyle.Raised;\\
\hspace*{0.5cm}dataGridA2.CellBorderStyle = DataGridViewCellBorderStyle.Single;\\
\hspace*{0.5cm}dataGridA2.GridColor = Color.DodgerBlue;\\
\hspace*{0.5cm}dataGridA2.ColumnHeadersVisible = false;\\
\hspace*{0.5cm}dataGridA2.RowHeadersVisible = false;\\
\hspace*{0.5cm}dataGridA2.BackgroundColor = Color.WhiteSmoke;\\
\hspace*{0.5cm}dataGridA2.BorderStyle = BorderStyle.None;\\
\hspace*{0.5cm}dataGridA2.AllowUserToResizeColumns = true;\\[0.2cm]
\hspace*{1cm}// dataGridA2.DefaultCellStyle = Font(12.75F);\\[0.1cm]
\hspace*{1cm}// Set the column header style.\\[0.1cm]
\hspace*{0.5cm}DataGridViewCellStyle columnHeaderStyle = new DataGridViewCellStyle();\\
\hspace*{0.5cm} columnHeaderStyle.BackColor = Color.AntiqueWhite;\\
\hspace*{0.5cm}columnHeaderStyle.Alignment = DataGridViewContentAlignment.MiddleCenter;\\
\hspace*{0.5cm}dataGridResAddition.ColumnHeadersDefaultCellStyle = columnHeaderStyle;\\[0.1cm]
\hspace*{0.5cm}DataGridViewCell CellR = new DataGridViewTextBoxCell();\\
\hspace*{0.5cm}CellR.Style.Alignment = DataGridViewContentAlignment.MiddleRight;\\[0.1cm]
 \hspace*{0.5cm}DataGridViewCell CellL = new DataGridViewTextBoxCell();\\
 \hspace*{0.5cm}CellL.Style.Alignment = DataGridViewContentAlignment.MiddleLeft;\\[0.1cm]
\hspace*{0.5cm}int line = 0;\\
 \hspace*{0.5cm}line = Convert.ToInt16(textLineA2.Text);\\
 \hspace*{0.5cm}dataGridA2.RowCount = line;\\[0.1cm]
 \hspace*{0.5cm}for (int i = 0; i $<$ column; i++)\\
 \hspace*{0.5cm}$\{$\\
\hspace*{1.5cm}//dataGridA2.Columns[i].Name = "C" + (i + 1);\\
 \hspace*{1cm}dataGridA2.Columns[i].Name = "C" + (i + 1);\\
 \hspace*{1cm}dataGridA2.Columns[i].CellTemplate = CellR;\\
 \hspace*{1cm}dataGridA2.Columns[i].Width = 30;\\
 \hspace*{0.5cm}$\}$\\
\hspace*{0.3cm}$\}$\\
\hspace*{0.3cm}private void initMatrixB2()\\
\hspace*{0.3cm}$\{$\\[0.1cm]
\hspace*{1cm}// Create an unbound DataGridView by declaring a column count.\\[0.1cm]
\hspace*{0.5cm}int column = 0;\\
\hspace*{0.5cm}column = Convert.ToInt16(textcolumnB2.Text);\\
\hspace*{0.5cm}dataGridB2.ColumnCount = column;\\
\hspace*{0.5cm}dataGridB2.AllowUserToOrderColumns = false;\\
\hspace*{0.5cm}dataGridB2.AllowUserToAddRows = false;\\
\hspace*{0.5cm}dataGridB2.Enabled = true;\\
\hspace*{0.5cm}dataGridB2.AutoSizeRowsMode = DataGridViewAutoSizeRowsMode.\\
DisplayedCellsExceptHeaders;\\
\hspace*{0.5cm}dataGridB2.ColumnHeadersBorderStyle = DataGridViewHeaderBorderStyle.Raised;\\
\hspace*{0.5cm}dataGridB2.CellBorderStyle = DataGridViewCellBorderStyle.Single;\\
\hspace*{0.5cm}dataGridB2.GridColor = Color.DodgerBlue;\\
\hspace*{0.5cm}dataGridB2.ColumnHeadersVisible = false;\\
\hspace*{0.5cm}dataGridB2.RowHeadersVisible = false;\\
\hspace*{0.5cm}dataGridB2.BackgroundColor = Color.WhiteSmoke;\\
\hspace*{0.5cm}dataGridB2.BorderStyle = BorderStyle.None;\\
\hspace*{0.5cm}dataGridB2.AllowUserToResizeColumns = true;\\[0.1cm]
\hspace*{1cm}// dataGridB2.DefaultCellStyle = Font(12.75F);\\[0.1cm]
\hspace*{1cm}// Set the column header style.\\[0.1cm]
\hspace*{0.5cm}DataGridViewCellStyle columnHeaderStyle = new DataGridViewCellStyle();\\
\hspace*{0.5cm}columnHeaderStyle.BackColor = Color.AntiqueWhite;\\
\hspace*{0.5cm}columnHeaderStyle.Alignment = DataGridViewContentAlignment.MiddleCenter;\\
\hspace*{0.5cm}dataGridB2.ColumnHeadersDefaultCellStyle = columnHeaderStyle;\\[0.1cm]
\hspace*{0.5cm}DataGridViewCell CellR = new DataGridViewTextBoxCell();\\
\hspace*{0.5cm}CellR.Style.Alignment = DataGridViewContentAlignment.MiddleRight;\\[0.1cm]
\hspace*{0.5cm}DataGridViewCell CellL = new DataGridViewTextBoxCell();\\
\hspace*{0.5cm}CellL.Style.Alignment = DataGridViewContentAlignment.MiddleLeft;\\[0.1cm]
\hspace*{0.5cm}int line = 0;\\
\hspace*{0.5cm}line = Convert.ToInt16(textLineB2.Text);\\
\hspace*{0.5cm}dataGridB2.RowCount = line;\\[0.1cm]
\hspace*{0.5cm}for (int i = 0; i $<$ column; i++)\\
\hspace*{0.5cm}$\{$\\
\hspace*{1cm}dataGridB2.Columns[i].Name = "C" + (i + 1);\\
\hspace*{1cm}dataGridB2.Columns[i].CellTemplate = CellR;\\
\hspace*{1cm}dataGridB2.Columns[i].Width = 30;\\
\hspace*{0.5cm}$\}$\\[0.2cm]
\hspace*{0.3cm}$\}$\\
\hspace*{0.3cm}private void initMatrixA3()\\
\hspace*{0.3cm}$\{$\\[0.1cm]
\hspace*{1cm}// Create an unbound DataGridView by declaring a column count.\\[0.1cm]
\hspace*{0.5cm} int column = 0;\\
\hspace*{0.5cm}column = Convert.ToInt16(textcolumnA3.Text);\\
\hspace*{0.5cm}dataGridA3.ColumnCount = column;\\
\hspace*{0.5cm}dataGridA3.AllowUserToOrderColumns = false;\\
\hspace*{0.5cm}dataGridA3.AllowUserToAddRows = false;\\
\hspace*{0.5cm}dataGridA3.Enabled = true;\\
\hspace*{0.5cm}dataGridA3.AutoSizeRowsMode = DataGridViewAutoSizeRowsMode.\\
DisplayedCellsExceptHeaders;\\
\hspace*{0.5cm}dataGridA3.ColumnHeadersBorderStyle = DataGridViewHeaderBorderStyle.Raised;\\
\hspace*{0.5cm}dataGridA3.CellBorderStyle = DataGridViewCellBorderStyle.Single;\\
\hspace*{0.5cm}dataGridA3.GridColor = Color.DodgerBlue;\\
\hspace*{0.5cm}dataGridA3.ColumnHeadersVisible = false;\\
\hspace*{0.5cm}dataGridA3.RowHeadersVisible = false;\\
\hspace*{0.5cm}dataGridA3.BackgroundColor = Color.WhiteSmoke;\\
\hspace*{0.5cm}dataGridA3.BorderStyle = BorderStyle.None;\\
\hspace*{0.5cm}dataGridA3.AllowUserToResizeColumns = true;\\[0.1cm]
\hspace*{1cm}// dataGridA3.DefaultCellStyle = Font(12.75F);\\[0.1cm]
 \hspace*{1cm}// Set the column header style.\\
\hspace*{0.5cm}DataGridViewCellStyle columnHeaderStyle = new DataGridViewCellStyle();\\
\hspace*{0.5cm}columnHeaderStyle.BackColor = Color.AntiqueWhite;\\
\hspace*{0.5cm}columnHeaderStyle.Alignment = DataGridViewContentAlignment.MiddleCenter;\\
\hspace*{0.5cm} dataGridA3.ColumnHeadersDefaultCellStyle = columnHeaderStyle;\\[0.1cm]
\hspace*{0.5cm}DataGridViewCell CellR = new DataGridViewTextBoxCell();\\
\hspace*{0.5cm}CellR.Style.Alignment = DataGridViewContentAlignment.MiddleRight;\\[0.1cm]
\hspace*{0.5cm}DataGridViewCell CellL = new DataGridViewTextBoxCell();\\
\hspace*{0.5cm}CellL.Style.Alignment = DataGridViewContentAlignment.MiddleLeft;\\[0.1cm]
\hspace*{0.5cm}int line = 0;\\
\hspace*{0.5cm} line = Convert.ToInt16(textLineA3.Text);\\
\hspace*{0.5cm}dataGridA3.RowCount = line;\\[0.1cm]
\hspace*{0.5cm}for (int i = 0; i $<$ column; i++)\\
\hspace*{0.5cm}$\{$\\
\hspace*{1cm}dataGridA3.Columns[i].Name = "C" + (i +1);\\
\hspace*{1cm}dataGridA3.Columns[i].CellTemplate = CellR;\\
\hspace*{1cm}dataGridA3.Columns[i].Width = 30;\\
\hspace*{0.5cm}$\}$\\[0.2cm]
\hspace*{0.3cm}$\}$\\
\hspace*{0.3cm}private void initMatrixA4()\\
\hspace*{0.3cm}$\{$\\[0.1cm]
\hspace*{1cm}// Create an unbound DataGridView by declaring a column count.\\
\hspace*{0.5cm}int column = 0;\\
\hspace*{0.5cm}column = Convert.ToInt16(textlineA4.Text);\\
\hspace*{0.5cm}dataGridA4.ColumnCount = column;\\
\hspace*{0.5cm}dataGridA4.AllowUserToOrderColumns = false;\\
\hspace*{0.5cm}dataGridA4.AllowUserToAddRows = false;\\
\hspace*{0.5cm}dataGridA4.Enabled = true;\\
\hspace*{0.5cm}dataGridA4.AutoSizeRowsMode = DataGridViewAutoSizeRowsMode.\\
DisplayedCellsExceptHeaders;\\
\hspace*{0.5cm}dataGridA4.ColumnHeadersBorderStyle = DataGridViewHeaderBorderStyle.Raised;\\
\hspace*{0.5cm}dataGridA4.CellBorderStyle = DataGridViewCellBorderStyle.Single;\\
\hspace*{0.5cm}dataGridA4.GridColor = Color.DodgerBlue;\\
\hspace*{0.5cm}dataGridA4.ColumnHeadersVisible = false;\\
\hspace*{0.5cm}dataGridA4.RowHeadersVisible = false;\\
\hspace*{0.5cm}dataGridA4.BackgroundColor = Color.WhiteSmoke;\\
\hspace*{0.5cm}dataGridA4.BorderStyle = BorderStyle.None;\\
\hspace*{0.5cm}dataGridA4.AllowUserToResizeColumns = true;\\[0.1cm]
\hspace*{1cm}// dataGridA4.DefaultCellStyle = Font(12.75F);\\[0.1cm]
\hspace*{1cm}// Set the column header style.\\[0.1cm]
\hspace*{0.5cm}DataGridViewCellStyle columnHeaderStyle = new DataGridViewCellStyle();\\
\hspace*{0.5cm}columnHeaderStyle.BackColor = Color.AntiqueWhite;\\
\hspace*{0.5cm}columnHeaderStyle.Alignment = DataGridViewContentAlignment.MiddleCenter;\\
\hspace*{0.5cm}dataGridA4.ColumnHeadersDefaultCellStyle = columnHeaderStyle;\\[0.1cm]
\hspace*{0.5cm}DataGridViewCell CellR = new DataGridViewTextBoxCell();\\
\hspace*{0.5cm}CellR.Style.Alignment = DataGridViewContentAlignment.MiddleRight;\\[0.1cm]
\hspace*{0.5cm}DataGridViewCell CellL = new DataGridViewTextBoxCell();\\
\hspace*{0.5cm}CellL.Style.Alignment = DataGridViewContentAlignment.MiddleLeft;\\[0.1cm]
\hspace*{0.5cm}int line = 0;\\
\hspace*{0.5cm}line = Convert.ToInt16(textlineA4.Text);\\
\hspace*{0.5cm}dataGridA4.RowCount = line;\\[0.1cm]
\hspace*{0.5cm}for (int i = 0; i $<$ column; i++)\\
\hspace*{0.5cm}$\{$\\
\hspace*{1cm}dataGridA4.Columns[i].Name = "C" + (i + 1);\\
\hspace*{1cm}dataGridA4.Columns[i].CellTemplate = CellR;\\
\hspace*{0.5cm}dataGridA4.Columns[i].Width = 30;\\
\hspace*{0.5cm}$\}$\\[0.2cm]
\hspace*{0.3cm}$\}$\\
\hspace*{0.3cm}private void initMatrixA5()\\
\hspace*{0.3cm}$\{$\\[0.1cm]
\hspace*{1cm}// Create an unbound DataGridView by declaring a column count.\\[0.1cm]
\hspace*{0.5cm}int column = 0;\\
\hspace*{0.5cm}column = Convert.ToInt16(textlineA5.Text);\\
\hspace*{0.5cm}dataGridA5.ColumnCount = column;\\
\hspace*{0.5cm}dataGridA5.AllowUserToOrderColumns =false;\\
\hspace*{0.5cm}dataGridA5.AllowUserToAddRows = false;\\
\hspace*{0.5cm}dataGridA5.Enabled = true;\\
\hspace*{0.5cm}dataGridA5.AutoSizeRowsMode = DataGridViewAutoSizeRowsMode.\\
DisplayedCellsExceptHeaders;\\
\hspace*{0.5cm}dataGridA5.ColumnHeadersBorderStyle =
 DataGridViewHeaderBorderStyle.Raised;\\
\hspace*{0.5cm}dataGridA5.CellBorderStyle =
 DataGridViewCellBorderStyle.Single;\\
\hspace*{0.5cm}dataGridA5.GridColor = Color.DodgerBlue;\\
\hspace*{0.5cm}dataGridA5.ColumnHeadersVisible =false;\\
\hspace*{0.5cm}dataGridA5.RowHeadersVisible = false;\\
\hspace*{0.5cm}dataGridA5.BackgroundColor = Color.WhiteSmoke;\\
\hspace*{0.5cm}dataGridA5.BorderStyle = BorderStyle.None;\\
\hspace*{0.5cm}dataGridA5.AllowUserToResizeColumns =true;\\[0.1cm]
\hspace*{1cm}// dataGridA5.DefaultCellStyle = Font(12.75F);\\[0.1cm]
\hspace*{1cm}// Set the column header style.\\[0.1cm]
\hspace*{0.5cm}DataGridViewCellStyle columnHeaderStyle = new DataGridViewCellStyle();\\
\hspace*{0.5cm}columnHeaderStyle.BackColor = Color.AntiqueWhite;\\
\hspace*{0.5cm}columnHeaderStyle.Alignment = DataGridViewContentAlignment.MiddleCenter;\\
\hspace*{0.5cm}dataGridA5.ColumnHeadersDefaultCellStyle = columnHeaderStyle;\\[0.1cm]
\hspace*{0.5cm}DataGridViewCell CellR = new DataGridViewTextBoxCell();\\
\hspace*{0.5cm}CellR.Style.Alignment = DataGridViewContentAlignment.MiddleRight;\\[0.1cm]
\hspace*{0.5cm}DataGridViewCell CellL = new DataGridViewTextBoxCell();\\
\hspace*{0.5cm}CellL.Style.Alignment = DataGridViewContentAlignment.MiddleLeft;\\[0.1cm]
\hspace*{0.5cm}int line = 0;\\
\hspace*{0.5cm}line = Convert.ToInt16(textlineA5.Text);\\
\hspace*{0.5cm}dataGridA5.RowCount = line;\\[0.1cm]
\hspace*{0.5cm}for (int i = 0; i $<$ column; i++)\\
\hspace*{0.5cm}$\{$\\
\hspace*{1cm}dataGridA5.Columns[i].Name = "C" + (i +1);\\
\hspace*{1cm}dataGridA5.Columns[i].CellTemplate =CellR;\\
\hspace*{1cm}dataGridA5.Columns[i].Width = 30;\\
\hspace*{0.5cm}$\}$\\[0.2cm]
\hspace*{0.3cm}$\}$\\
\hspace*{0.3cm}private void initMatrixA6()\\
\hspace*{0.3cm}$\{$\\
\hspace*{0.5cm}int column = 0;\\
\hspace*{0.5cm}column = Convert.ToInt16(textlineA6.Text);\\
\hspace*{0.5cm}dataGridA6.ColumnCount = column;\\
\hspace*{0.5cm}dataGridA6.AllowUserToOrderColumns =false;\\
\hspace*{0.5cm}dataGridA6.AllowUserToAddRows = false;\\
\hspace*{0.5cm}dataGridA6.Enabled = true;\\
\hspace*{0.5cm}dataGridA6.AutoSizeRowsMode = DataGridViewAutoSizeRowsMode.\\
DisplayedCellsExceptHeaders;\\
\hspace*{0.5cm}dataGridA6.ColumnHeadersBorderStyle =
 DataGridViewHeaderBorderStyle.Raised;\\
\hspace*{0.5cm}dataGridA6.CellBorderStyle =
 DataGridViewCellBorderStyle.Single;\\
\hspace*{0.5cm}dataGridA6.GridColor = Color.DodgerBlue;\\
\hspace*{0.5cm}dataGridA6.ColumnHeadersVisible =false;\\
\hspace*{0.5cm}dataGridA6.RowHeadersVisible = false;\\
\hspace*{0.5cm}dataGridA6.BackgroundColor = Color.WhiteSmoke;\\
\hspace*{0.5cm}dataGridA6.BorderStyle = BorderStyle.None;\\
\hspace*{0.5cm}dataGridA6.AllowUserToResizeColumns =true;\\[0.1cm]
\hspace*{0.5cm}DataGridViewCellStyle columnHeaderStyle = new DataGridViewCellStyle();\\
\hspace*{0.5cm}columnHeaderStyle.BackColor = Color.AntiqueWhite;\\
\hspace*{0.5cm}columnHeaderStyle.Alignment =
 DataGridViewContentAlignment.MiddleCenter;\\
\hspace*{0.5cm}dataGridA6.ColumnHeadersDefaultCellStyle = columnHeaderStyle;\\[0.1cm]
\hspace*{0.5cm}DataGridViewCell CellR = new DataGridViewTextBoxCell();\\
\hspace*{0.5cm}CellR.Style.Alignment = DataGridViewContentAlignment.MiddleRight;\\[0.1cm]
\hspace*{0.5cm}DataGridViewCell CellL = new DataGridViewTextBoxCell();\\
\hspace*{0.5cm}CellL.Style.Alignment = DataGridViewContentAlignment.MiddleLeft;\\[0.1cm]
\hspace*{0.5cm}int line = 0;\\
\hspace*{0.5cm}line = Convert.ToInt16(textlineA6.Text);\\
\hspace*{0.5cm}dataGridA6.RowCount = line;\\[0.1cm]
\hspace*{0.5cm}for (int i = 0; i $<$ column; i++)\\
\hspace*{0.5cm}$\{$\\
\hspace*{1cm}dataGridA6.Columns[i].Name = "C" + (i +1);\\
\hspace*{1cm}dataGridA6.Columns[i].CellTemplate =CellR;\\
\hspace*{1cm}dataGridA6.Columns[i].Width = 30;\\
\hspace*{0.5cm}$\}$\\[0.1cm]
\hspace*{0.3cm}$\}$\\
\hspace*{0.3cm}private void initMatrixX0()\\
\hspace*{0.3cm}$\{$\\[0.1cm]
\hspace*{0.5cm}int column = 0;\\
\hspace*{0.5cm}column = 1;\\
\hspace*{0.5cm}dataGridX0.ColumnCount = column;\\
\hspace*{0.5cm}dataGridX0.AllowUserToOrderColumns = false;\\
\hspace*{0.5cm}dataGridX0.AllowUserToAddRows = false;\\
\hspace*{0.5cm}dataGridX0.Enabled = true;\\
\hspace*{0.5cm}dataGridX0.AutoSizeRowsMode = DataGridViewAutoSizeRowsMode.\\
DisplayedCellsExceptHeaders;\\
\hspace*{0.5cm}dataGridX0.ColumnHeadersBorderStyle = DataGridViewHeaderBorderStyle.Raised;\\
\hspace*{0.5cm}dataGridX0.CellBorderStyle = DataGridViewCellBorderStyle.Single;\\
\hspace*{0.5cm}dataGridX0.GridColor = Color.DodgerBlue;\\
\hspace*{0.5cm}dataGridX0.ColumnHeadersVisible = false;\\
\hspace*{0.5cm}dataGridX0.RowHeadersVisible = false;\\
\hspace*{0.5cm}dataGridX0.BackgroundColor = Color.WhiteSmoke;\\
\hspace*{0.5cm}dataGridX0.BorderStyle = BorderStyle.None;\\
\hspace*{0.5cm}dataGridX0.AllowUserToResizeColumns = true;\\[0.1cm]
\hspace*{0.5cm}DataGridViewCellStyle columnHeaderStyle = new DataGridViewCellStyle();\\
\hspace*{0.5cm}columnHeaderStyle.BackColor = Color.AntiqueWhite;\\
\hspace*{0.5cm}columnHeaderStyle.Alignment = DataGridViewContentAlignment.MiddleCenter;\\
 \hspace*{0.5cm}dataGridX0.ColumnHeadersDefaultCellStyle = columnHeaderStyle;\\[0.1cm]
\hspace*{0.5cm}DataGridViewCell CellR = new DataGridViewTextBoxCell();\\
\hspace*{0.5cm}CellR.Style.Alignment = DataGridViewContentAlignment.MiddleRight;\\[0.1cm]
\hspace*{0.5cm}DataGridViewCell CellL = new DataGridViewTextBoxCell();\\
\hspace*{0.5cm}CellL.Style.Alignment = DataGridViewContentAlignment.MiddleLeft;\\[0.1cm]
\hspace*{0.5cm}int line = 0;\\
\hspace*{0.5cm}line = Convert.ToInt16(textlineA6.Text);\\
\hspace*{0.5cm}dataGridX0.RowCount = line;\\[0.1cm]
\hspace*{0.5cm}for (int i = 0; i $<$ column; i++)\\
\hspace*{0.5cm}$\{$\\
\hspace*{1cm}dataGridX0.Columns[i].Name = "C" + (i + 1);\\
\hspace*{1cm}dataGridX0.Columns[i].CellTemplate = CellR;\\
\hspace*{1cm}dataGridX0.Columns[i].Width = 30;\\
\hspace*{0.5cm}$\}$\\[0.1cm]
\hspace*{0.3cm}$\}$\\
\hspace*{0.3cm}private void initMatrixAdditionResult()\\
\hspace*{0.3cm}$\{$\\[0.1cm]
\hspace*{1cm}// Create an unbound DataGridView by declaring a column count.\\[0.1cm]
\hspace*{0.5cm}int column = 0;\\
\hspace*{0.5cm}column = Convert.ToInt16(textcolumn.Text);\\
\hspace*{0.5cm}dataGridResAddition.ColumnCount = column;\\
\hspace*{0.5cm}dataGridResAddition.AllowUserToOrderColumns = false;\\
\hspace*{0.5cm}dataGridResAddition.AllowUserToAddRows = false;\\
\hspace*{0.5cm}dataGridResAddition.Enabled = true;\\
\hspace*{0.5cm}dataGridResAddition.AutoSizeRowsMode = DataGridViewAutoSizeRowsMode.\\
DisplayedCellsExceptHeaders;\\
\hspace*{0.5cm}dataGridResAddition.ColumnHeadersBorderStyle = DataGridViewHeaderBorderStyle.Raised;\\
\hspace*{0.5cm}dataGridResAddition.CellBorderStyle = DataGridViewCellBorderStyle.Single;\\
\hspace*{0.5cm}dataGridResAddition.GridColor = Color.DodgerBlue;\\
\hspace*{0.5cm}dataGridResAddition.ColumnHeadersVisible = false;\\
\hspace*{0.5cm}dataGridResAddition.RowHeadersVisible = false;\\
\hspace*{0.5cm}dataGridResAddition.BackgroundColor = Color.WhiteSmoke;\\
\hspace*{0.5cm}dataGridResAddition.BorderStyle = BorderStyle.None;\\
\hspace*{0.5cm}dataGridResAddition.AllowUserToResizeColumns = true;\\[0.1cm]
\hspace*{1cm}// dataGridResAddition.DefaultCellStyle = Font(12.75F);\\[0.1cm]
\hspace*{1cm}// Set the column header style.\\[0.1cm]
\hspace*{0.5cm}DataGridViewCellStyle columnHeaderStyle = new DataGridViewCellStyle();\\
\hspace*{0.5cm}columnHeaderStyle.BackColor = Color.AntiqueWhite;\\
\hspace*{0.5cm} columnHeaderStyle.Alignment = DataGridViewContentAlignment.MiddleCenter;\\
\hspace*{0.5cm}dataGridResAddition.ColumnHeadersDefaultCellStyle = columnHeaderStyle;\\[0.1cm]
\hspace*{0.5cm}DataGridViewCell CellR = new DataGridViewTextBoxCell();\\
\hspace*{0.5cm}CellR.Style.Alignment = DataGridViewContentAlignment.MiddleRight;\\[0.1cm]
\hspace*{0.5cm}DataGridViewCell CellL = new DataGridViewTextBoxCell();\\
\hspace*{0.5cm}CellL.Style.Alignment = DataGridViewContentAlignment.MiddleLeft;\\[0.1cm]
\hspace*{0.5cm}int line = 0;\\
\hspace*{0.5cm}line = Convert.ToInt16(textLine.Text);\\
\hspace*{0.5cm}dataGridResAddition.RowCount = line;\\[0.1cm]
\hspace*{0.5cm}for (int i = 0; i $<$ column; i++)\\
\hspace*{0.5cm}$\{$\\
\hspace*{1cm}dataGridResAddition.Columns[i].Name = "C" + (i + 1);\\
\hspace*{1cm}dataGridResAddition.Columns[i].CellTemplate = CellR;\\
\hspace*{1cm}dataGridResAddition.Columns[i].Width = 30;\\
\hspace*{0.5cm}$\}$\\[0.1cm]
\hspace*{0.3cm}$\}$\\
\hspace*{0.3cm}private void initMatrixProduct()\\
\hspace*{0.3cm}$\{$\\[0.1cm]
\hspace*{1cm}// Create an unbound DataGridView by declaring a column count.\\
\hspace*{0.5cm}int column = 0;\\
\hspace*{0.5cm}column = Convert.ToInt16(textcolumnB2.Text);\\
\hspace*{0.5cm}dataGridProduct.ColumnCount = column;\\
\hspace*{0.5cm}dataGridProduct.AllowUserToOrderColumns = false;\\
\hspace*{0.5cm}dataGridProduct.AllowUserToAddRows = false;\\
\hspace*{0.5cm}dataGridProduct.Enabled = true;\\
\hspace*{0.5cm}dataGridProduct.AutoSizeRowsMode = DataGridViewAutoSizeRowsMode.\\
DisplayedCellsExceptHeaders;\\
\hspace*{0.5cm}dataGridProduct.ColumnHeadersBorderStyle = DataGridViewHeaderBorderStyle.Raised;\\
\hspace*{0.5cm}dataGridProduct.CellBorderStyle = DataGridViewCellBorderStyle.Single;\\
\hspace*{0.5cm}dataGridProduct.GridColor = Color.DodgerBlue;\\
\hspace*{0.5cm}dataGridProduct.ColumnHeadersVisible = false;\\
\hspace*{0.5cm}dataGridProduct.RowHeadersVisible = false;\\
\hspace*{0.5cm}dataGridProduct.BackgroundColor = Color.WhiteSmoke;\\
\hspace*{0.5cm}dataGridProduct.BorderStyle = BorderStyle.None;\\
\hspace*{0.5cm}dataGridProduct.AllowUserToResizeColumns = true;\\[0.1cm]
\hspace*{1cm}// dataGridProduct.DefaultCellStyle = Font(12.75F);\\[0.1cm]
\hspace*{1cm}// Set the column header style.\\[0.1cm]
\hspace*{0.5cm}DataGridViewCellStyle columnHeaderStyle = new DataGridViewCellStyle();\\
\hspace*{0.5cm}columnHeaderStyle.BackColor = Color.AntiqueWhite;\\
\hspace*{0.5cm} columnHeaderStyle.Alignment = DataGridViewContentAlignment.MiddleCenter;\\
\hspace*{0.5cm}dataGridProduct.ColumnHeadersDefaultCellStyle = columnHeaderStyle;\\[0.1cm]
\hspace*{0.5cm}DataGridViewCell CellR = new DataGridViewTextBoxCell();\\
\hspace*{0.5cm}CellR.Style.Alignment = DataGridViewContentAlignment.MiddleRight;\\[0.1cm]
\hspace*{0.5cm}DataGridViewCell CellL = new DataGridViewTextBoxCell();\\
\hspace*{0.5cm}CellL.Style.Alignment = DataGridViewContentAlignment.MiddleLeft;\\[0.1cm]
\hspace*{0.5cm}int line = 0;\\
\hspace*{0.5cm}line = Convert.ToInt16(textLineA2.Text);\\
\hspace*{0.5cm}dataGridProduct.RowCount = line;\\[0.1cm]
\hspace*{0.5cm}for (int i = 0; i $<$ column; i++)\\
\hspace*{0.5cm}$\{$\\
\hspace*{1cm}dataGridProduct.Columns[i].Name = "C" + (i + 1);\\
\hspace*{1cm}dataGridProduct.Columns[i].CellTemplate = CellR;\\
\hspace*{1cm}dataGridProduct.Columns[i].Width = 30;\\
\hspace*{0.5cm}$\}$\\[0.2cm]
\hspace*{0.3cm}$\}$\\
\hspace*{0.3cm}private void initMatrixScalarProduct()\\
\hspace*{0.3cm}$\{$\\[0.1cm]
\hspace*{1cm}// Create an unbound DataGridView by declaring a column count.\\[0.1cm]
\hspace*{0.5cm}int column = 0;\\
\hspace*{0.5cm}column = Convert.ToInt16(textcolumnA3.Text);\\
\hspace*{0.5cm}dataGridScalarProduct.ColumnCount = column;\\
\hspace*{0.5cm}dataGridScalarProduct.AllowUserToOrderColumns = false;\\
\hspace*{0.5cm}dataGridScalarProduct.AllowUserToAddRows = false;\\
\hspace*{0.5cm}dataGridScalarProduct.Enabled = true;\\
\hspace*{0.5cm}dataGridScalarProduct.AutoSizeRowsMode = DataGridViewAutoSizeRowsMode.\\
DisplayedCellsExceptHeaders;\\
\hspace*{0.5cm}dataGridScalarProduct.ColumnHeadersBorderStyle = DataGridViewHeaderBorderStyle.Raised;\\
\hspace*{0.5cm}dataGridScalarProduct.CellBorderStyle = DataGridViewCellBorderStyle.Single;\\
\hspace*{0.5cm}dataGridScalarProduct.GridColor = Color.DodgerBlue;\\
\hspace*{0.5cm}dataGridScalarProduct.ColumnHeadersVisible = false;\\
\hspace*{0.5cm}dataGridScalarProduct.RowHeadersVisible = false;\\
\hspace*{0.5cm}dataGridScalarProduct.BackgroundColor = Color.WhiteSmoke;\\
\hspace*{0.5cm}dataGridScalarProduct.BorderStyle = BorderStyle.None;\\
\hspace*{0.5cm}dataGridScalarProduct.AllowUserToResizeColumns = true;\\[0.1cm]
\hspace*{1cm}// dataGridProduct.DefaultCellStyle = Font(12.75F);\\[0.1cm]
\hspace*{1cm}// Set the column header style.\\[0.1cm]
\hspace*{0.5cm}DataGridViewCellStyle columnHeaderStyle = new DataGridViewCellStyle();\\
\hspace*{0.5cm}columnHeaderStyle.BackColor = Color.AntiqueWhite;\\
\hspace*{0.5cm}columnHeaderStyle.Alignment = DataGridViewContentAlignment.MiddleCenter;\\
\hspace*{0.5cm}dataGridScalarProduct.ColumnHeadersDefaultCellStyle = columnHeaderStyle;\\[0.1cm]
\hspace*{0.5cm}DataGridViewCell CellR = new DataGridViewTextBoxCell();\\
\hspace*{0.5cm}CellR.Style.Alignment = DataGridViewContentAlignment.MiddleRight;\\[0.1cm]
\hspace*{0.5cm}DataGridViewCell CellL = new DataGridViewTextBoxCell();\\
\hspace*{0.5cm}CellL.Style.Alignment = DataGridViewContentAlignment.MiddleLeft;\\[0.1cm]
\hspace*{0.5cm}int line = 0;\\
\hspace*{0.5cm}line = Convert.ToInt16(textLineA3.Text);\\
\hspace*{0.5cm}dataGridScalarProduct.RowCount = line;\\[0.1cm]
\hspace*{0.5cm}for (int i = 0; i $<$ column; i++)\\
\hspace*{0.5cm}$\{$\\
\hspace*{1cm}dataGridScalarProduct.Columns[i].Name = "C" + (i + 1);\\
\hspace*{1cm}dataGridScalarProduct.Columns[i].CellTemplate = CellR;\\
\hspace*{1cm}dataGridScalarProduct.Columns[i].Width = 30;\\
\hspace*{0.5cm}$\}$\\
\hspace*{0.3cm}$\}$\\
\hspace*{0.3cm}private void initPowerMatrix ()\\
\hspace*{0.3cm}$\{$\\[0.1cm]
\hspace*{1cm}// Create an unbound DataGridView by declaring a column count.\\[0.1cm]
\hspace*{0.5cm}int column = 0;\\
\hspace*{0.5cm}column = Convert.ToInt16(textlineA4.Text);\\
\hspace*{0.5cm}dataGridMatrix$_{-}$at$_{-}$power$_{-}$n.ColumnCount = column;\\
\hspace*{0.5cm}dataGridMatrix$_{-}$at$_{-}$power$_{-}$n.AllowUserToOrderColumns=false;\\
\hspace*{0.5cm}dataGridMatrix$_{-}$at$_{-}$power$_{-}$n.AllowUserToAddRows=false;\\
\hspace*{0.5cm}dataGridMatrix$_{-}$at$_{-}$power$_{-}$n.Enabled =true;\\
\hspace*{0.5cm}dataGridMatrix$_{-}$at$_{-}$power$_{-}$n.AutoSizeRowsMode=DataGridViewAutoSizeRowsMode.\\
DisplayedCellsExceptHeaders;\\
\hspace*{0.5cm}dataGridMatrix$_{-}$at$_{-}$power$_{-}$n.ColumnHeadersBorderStyle = DataGridViewHeader\\
BorderStyle.Raised;\\
\hspace*{0.5cm}dataGridMatrix$_{-}$at$_{-}$power$_{-}$n.CellBorderStyle = DataGridViewCellBorderStyle.Single;\\
\hspace*{0.5cm}dataGridMatrix$_{-}$at$_{-}$power$_{-}$n.GridColor =Color.DodgerBlue;\\
\hspace*{0.5cm}dataGridMatrix$_{-}$at$_{-}$power$_{-}$n.ColumnHeadersVisible=false;\\
\hspace*{0.5cm}dataGridMatrix$_{-}$at$_{-}$power$_{-}$n.RowHeadersVisible=false;\\
\hspace*{0.5cm}dataGridMatrix$_{-}$at$_{-}$power$_{-}$n.BackgroundColor = Color.WhiteSmoke;\\
\hspace*{0.5cm}dataGridMatrix$_{-}$at$_{-}$power$_{-}$n.BorderStyle = BorderStyle.None;\\
\hspace*{0.5cm}dataGridMatrix$_{-}$at$_{-}$power$_{-}$n.AllowUserToResizeColumns=true;\\[0.1cm]
\hspace*{1cm}// dataGridProduct.DefaultCellStyle = Font(12.75F);\\[0.1cm]
\hspace*{1cm}// Set the column header style.\\[0.1cm]
\hspace*{0.5cm}DataGridViewCellStyle columnHeaderStyle = new DataGridViewCellStyle();\\
\hspace*{0.5cm}columnHeaderStyle.BackColor = Color.AntiqueWhite;\\
\hspace*{0.5cm}columnHeaderStyle.Alignment = DataGridViewContentAlignment.MiddleCenter;\\
\hspace*{0.5cm}dataGridMatrix$_{-}$at$_{-}$power$_{-}$n.ColumnHeadersDefaultCellStyle = columnHeaderStyle;\\[0.1cm]
\hspace*{0.5cm}DataGridViewCell CellR = new DataGridViewTextBoxCell();\\
\hspace*{0.5cm}CellR.Style.Alignment = DataGridViewContentAlignment.MiddleRight;\\[0.1cm]
\hspace*{0.5cm}DataGridViewCell CellL = new DataGridViewTextBoxCell();\\
\hspace*{0.5cm}CellL.Style.Alignment = DataGridViewContentAlignment.MiddleLeft;\\[0.1cm]
\hspace*{0.5cm}int line = 0;\\
\hspace*{0.5cm}line = Convert.ToInt16(textlinieA4.Text);\\
\hspace*{0.5cm}dataGridMatrix$_{-}$at$_{-}$power$_{-}$n.RowCount = line;\\[0.1cm]
\hspace*{0.5cm}for (int i = 0; i $<$ column; i++)\\
\hspace*{0.5cm}$\{$\\
\hspace*{1cm}dataGridMatrix$_{-}$at$_{-}$power$_{-}$n.Columns[i].Name = "C" + (i +1);\\
\hspace*{1cm}dataGridMatrix$_{-}$at$_{-}$power$_{-}$n.Columns[i].CellTemplate = CellR;\\
\hspace*{1cm}dataGridMatrix$_{-}$at$_{-}$power$_{-}$n.Columns[i].Width = 30;\\
\hspace*{0.5cm}$\}$\\[0.2cm]
\hspace*{0.3cm}$\}$\\
\hspace*{0.3cm}private void initMatrixSolutionSystem()\\
\hspace*{0.3cm}$\{$\\[0.1cm]
\hspace*{1cm}// Create an unbound DataGridView by declaring a column count.\\[0.1cm]
\hspace*{0.5cm}int column = 0;\\
\hspace*{0.5cm}column = 1;\\
\hspace*{0.5cm}dataGridSolutionXk.ColumnCount = column;\\
\hspace*{0.5cm}dataGridSolutionXk.AllowUserToOrderColumns = false;\\
\hspace*{0.5cm}dataGridSolutionXk.AllowUserToAddRows = false;\\
\hspace*{0.5cm}dataGridSolutionXk.Enabled = true;\\
\hspace*{0.5cm}dataGridSolutionXk.AutoSizeRowsMode = DataGridViewAutoSizeRowsMode.\\
DisplayedCellsExceptHeaders;\\
\hspace*{0.5cm}dataGridSolutionXk.ColumnHeadersBorderStyle = DataGridViewHeader\\
BorderStyle.Raised;\\
\hspace*{0.5cm}dataGridSolutionXk.CellBorderStyle = DataGridViewCellBorderStyle.Single;\\
\hspace*{0.5cm}dataGridSolutionXk.GridColor = Color.DodgerBlue;\\
\hspace*{0.5cm}dataGridSolutionXk.ColumnHeadersVisible = false;\\
\hspace*{0.5cm}dataGridSolutionXk.RowHeadersVisible = false;\\
\hspace*{0.5cm}dataGridSolutionXk.BackgroundColor = Color.WhiteSmoke;\\
\hspace*{0.5cm}dataGridSolutionXk.BorderStyle = BorderStyle.None;\\
\hspace*{0.5cm}dataGridSolutionXk.AllowUserToResizeColumns = true;\\[0.1cm]
\hspace*{0.5cm}DataGridViewCellStyle columnHeaderStyle = new DataGridViewCellStyle();\\
\hspace*{0.5cm}columnHeaderStyle.BackColor = Color.AntiqueWhite;\\
\hspace*{0.5cm}columnHeaderStyle.Alignment = DataGridViewContentAlignment.MiddleCenter;\\
\hspace*{0.5cm}dataGridSolutieXk.ColumnHeadersDefaultCellStyle = columnHeaderStyle;\\[0.1cm]
\hspace*{0.5cm}DataGridViewCell CellR = new DataGridViewTextBoxCell();\\
 \hspace*{0.5cm}CellR.Style.Alignment = DataGridViewContentAlignment.MiddleRight;\\[0.1cm]
\hspace*{0.5cm}DataGridViewCell CellL = new DataGridViewTextBoxCell();\\
 \hspace*{0.5cm}CellL.Style.Alignment = DataGridViewContentAlignment.MiddleLeft;\\[0.1cm]
\hspace*{0.5cm}int line = 0;\\
\hspace*{0.5cm}line = Convert.ToInt16(textlineA6.Text);\\
 \hspace*{0.5cm}dataGridSolutionXk.RowCount = line;\\[0.1cm]
 \hspace*{0.5cm}for (int i = 0; i $<$ column; i++)\\
 \hspace*{0.5cm}$\{$\\
 \hspace*{1cm} dataGridSolutionXk.Columns[i].Name = "C" + (i + 1);\\
 \hspace*{1cm}dataGridSolutionXk.Columns[i].CellTemplate = CellR;\\
 \hspace*{1cm}dataGridSolutionXk.Columns[i].Width = 30;\\[0.1cm]
 \hspace*{0.5cm}$\}$\\[0.2cm]
\hspace*{0.3cm}$\}$\\
\hspace*{0.3cm}private void btGenerating$_{-}$Click(object sender, EventArgs e)\\
\hspace*{0.3cm}$\{$\\[0.1cm]
 \hspace*{0.5cm}if (textLine.Text != "" $\&\&$ textColumn.Text != "")\\
 \hspace*{0.5cm}$\{$\\
 \hspace*{1cm}initMatrixA();\\
 \hspace*{1cm}initMatrixB();\\
 \hspace*{0.5cm}$\}$\\
 \hspace*{0.5cm}else\\
 \hspace*{0.5cm}$\{$\\
 \hspace*{1cm}if (textLine.Text == "" $\&\&$ textColumn.Text == "")\\
 \hspace*{1.5cm}MessageBox.Show("Introduce number of lines and number of columns!");\\
 \hspace*{1cm}else if (textLine.Text == "")\\
 \hspace*{1.5cm}MessageBox.Show("Introduce number of lines!");\\
 \hspace*{1cm}else if (textColumns.Text == "")\\
 \hspace*{1.5cm}MessageBox.Show("Introduce number of columns!");\\
 \hspace*{0.5cm}$\}$\\
\hspace*{0.3cm}$\}$\\
\hspace*{0.3cm}private void btGenerating2$_{-}$Click(object sender, EventArgs e)\\
\hspace*{0.3cm}$\{$\\[0.1cm]
\hspace*{0.5cm}if (textLineA2.Text != "" $\&\&$ textColumnA2.Text != "" $\&\&$ textLineB2.Text != "" $\&\&$ textcolumnB2.Text != "")\\
\hspace*{0.5cm}$\{$\\
\hspace*{1cm}initMatrixA2();\\
\hspace*{1cm}initMatrixB2();\\
\hspace*{0.5cm}$\}$\\
\hspace*{0.5cm}else\\
\hspace*{0.5cm}$\{$\\
\hspace*{1cm} if (textLineA2.Text == "" $\&\&$ textColumnA2.Text == "" $\&\&$ textLineB2.Text == "" $\&\&$ textcolumnB2.Text == "")\\
\hspace*{1.5cm} MessageBox.Show("Introduce number of lines and number of columns!");\\
\hspace*{1cm}else\\
\hspace*{1cm}$\{$\\
\hspace*{1.5cm}if (textLineA2.Text == "" $||$ textLineB2.Text == "")\\
\hspace*{1.8cm}MessageBox.Show("Introduce number of lines for matrix A!");\\
\hspace*{1.5cm}else if (textColumnA2.Text == "")\\
\hspace*{1.8cm}MessageBox.Show("Introduce number of columns for matrix A!");\\
\hspace*{1cm}else if (textcolumnB2.Text == "")\\
\hspace*{1.5cm}MessageBox.Show("Introduce number of columns for matrix B!");\\
\hspace*{1cm}$\}$\\
\hspace*{0.5cm}$\}$\\[0.1cm]
\hspace*{0.3cm}$\}$\\
\hspace*{0.3cm}private void Generating butons3$_{-}$Click(object sender, EventArgs e)\\
\hspace*{0.3cm}$\{$\\[0.1cm]
\hspace*{0.5cm}if (textLineA3.Text != "")\\
\hspace*{0.5cm}$\{$\\
\hspace*{1cm}  initMatrixA3();\\
\hspace*{0.5cm}$\}$\\
\hspace*{0.5cm}else\\
\hspace*{1cm}MessageBox.Show("Introduce number of lines and number
of columns!");\\
\hspace*{0.3cm}$\}$\\
\hspace*{0.3cm}private void GeneratingA4$_{-}$Click(object sender, EventArgs e)\\
\hspace*{0.3cm}$\{$\\[0.1cm]
\hspace*{0.5cm}if (textLineA4.Text != "")\\
\hspace*{0.5cm}$\{$\\
\hspace*{1cm}  initMatrixA4();\\
\hspace*{0.5cm}$\}$\\
\hspace*{0.5cm}else\\
\hspace*{1cm}MessageBox.Show("Introduce number of lines and number
of columns!");\\
\hspace*{0.3cm}$\}$\\
\hspace*{0.3cm}private void GeneratingA5$_{-}$Click(object sender, EventArgs e)\\
\hspace*{0.3cm}$\{$\\[0.1cm]
\hspace*{0.5cm}if (textLineA5.Text != "")\\
\hspace*{0.5cm}$\{$\\
\hspace*{1cm}  initMatrixA5();\\
\hspace*{0.5cm}$\}$\\
\hspace*{0.5cm}else\\
\hspace*{1cm}MessageBox.Show("Introduce number of lines and number of columns!");\\
\hspace*{0.3cm}$\}$\\
\hspace*{0.3cm}private void btGeneratingMatrices$_{-}$Click(object sender, EventArgs e)\\
 \hspace*{0.3cm}$\{$\\
 \hspace*{0.5cm}if (textlineA6.Text != "")\\
 \hspace*{0.5cm}$\{$\\
 \hspace*{1cm}initMatrixA6();\\
 \hspace*{1cm}initMatrixX0();\\
 \hspace*{0.5cm}$\}$\\
 \hspace*{0.5cm}else\\
 \hspace*{1cm}MessageBox.Show("Introduce number of lines and columns!");\\
\hspace*{0.3cm}$\}$\\
 \hspace*{0.3cm}private void btComputeSumClick(object sender, EventArgs e)\\
 \hspace*{0.3cm}$\{$\\[0.1cm]
 \hspace*{0.5cm}initMatrixAdditionResult();\\[0.1cm]
 \hspace*{0.5cm}int line = 0;\\
 \hspace*{0.5cm}int column = 0;\\
 \hspace*{0.5cm}column = Convert.ToInt16(textColumn.Text);\\
 \hspace*{0.5cm}line = Convert.ToInt16(textLine.Text);\\
 \hspace*{0.5cm}int[,] A = new int[line, column];\\
 \hspace*{0.5cm}int[,] B = new int[line, column];\\
 \hspace*{0.5cm}int[,] ResAddition = new int[line, column];\\
 \hspace*{0.5cm}for (int i = 0; i $<$ line; i++)\\
 \hspace*{0.5cm}$\{$\\
 \hspace*{1cm}for (int j = 0; j $<$ column; j++)\\
 \hspace*{1cm}$\{$\\
 \hspace*{1.5cm}MessageBox.Show(dataGridA.Rows[i].Cells[j].Value.ToString() )\\
 \hspace*{1.5cm}if (dataGridA.Rows[i].Cells[j].Value.ToString() == "E")\\
 \hspace*{1.8cm} A[i, j] = Int32.MaxValue;\\
 \hspace*{1.5cm}else\\
 \hspace*{1.8cm} A[i, j] = Convert.ToInt16(dataGridA.Rows[i].Cells[j].Value.ToString());\\
 \hspace*{1.5cm}if (dataGridB.Rows[i].Cells[j].Value.ToString() == "E")\\
 \hspace*{1.8cm}B[i, j] = Int32.MaxValue;\\
 \hspace*{1.5cm}else\\
 \hspace*{1.8cm}B[i, j] = Convert.ToInt16(dataGridB.Rows[i].Cells[j].Value.ToString());\\
\hspace*{1.5cm} if (A[i, j] $<$ B[i, j])\\
 \hspace*{1.8cm}ResAddition[i, j] = A[i, j];\\
 \hspace*{1.5cm}else\\
 \hspace*{1.8cm}ResAddition[i, j] = B[i, j];\\
 \hspace*{1cm}$\}$\\
 \hspace*{0.5cm}$\}$\\
 \hspace*{0.5cm}for (int i = 0; i $<$ line; i++)\\
 \hspace*{0.5cm}$\{$\\
 \hspace*{1cm}for (int j = 0; j $<$ column; j++)\\
 \hspace*{1cm}$\{$\\
 \hspace*{1.5cm}if (ResAddition[i, j] == Int32.MaxValue)\\
 \hspace*{1.8cm}dataGridResAddition.Rows[i].Cells[j].Value = "E";\\
 \hspace*{1.5cm}else\\
 \hspace*{1.8cm}  dataGridResAddition.Rows[i].Cells[j].Value = ResAddition[i, j];\\
 \hspace*{1cm}$\}$\\
 \hspace*{0.5cm}$\}$\\[0.1cm]
\hspace*{0.3cm}$\}$\\
\hspace*{0.3cm}private void btComputationProduct$_{-}$Click(object sender, EventArgs e)\\
\hspace*{0.3cm}$\{$\\[0.1cm]
\hspace*{0.5cm}initMatrixProduct();\\[0.1cm]
\hspace*{0.5cm}int lineA = 0; int lineB = 0;\\
\hspace*{0.5cm}int columnA = 0; int columnB = 0;\\
\hspace*{0.5cm}columnA = Convert.ToInt16(textColumnA2.Text);\\
\hspace*{0.5cm}lineA = Convert.ToInt16(textLineA2.Text);\\
\hspace*{0.5cm}columnB = Convert.ToInt16(textcolumnB2.Text);\\
\hspace*{0.5cm}lineB = Convert.ToInt16(textLineB2.Text);\\
\hspace*{0.5cm}int[,] A2 = new int[lineA, columnA];\\
\hspace*{0.5cm}int[,] B2 = new int[lineB, columnB];\\
\hspace*{0.5cm}int[,] Product = new int[lineA, columnB];\\
\hspace*{0.5cm}int[,] Sum = new int[lineA, columnB];\\
\hspace*{0.5cm}int k;\\
\hspace*{0.5cm}for (int i = 0; i $<$ lineA; i++)\\
\hspace*{0.5cm}$\{$\\
\hspace*{1cm}for (int j = 0; j $<$ columnA; j++)\\
\hspace*{1cm}$\{$\\
\hspace*{1cm}if (dataGridA2.Rows[i].Cells[j].Value.ToString() == "E")\\
\hspace*{1.5cm}A2[i, j] = Int32.MaxValue;\\
\hspace*{1cm}else\\
\hspace*{1.5cm}A2[i, j] = Convert.ToInt16(dataGridA2.Rows[i].Cells[j].Value.ToString());\\
\hspace*{1cm}$\}$\\
\hspace*{0.5cm}$\}$\\[0.1cm]
\hspace*{0.5cm}for (int i = 0; i $<$ lineB; i++)\\
\hspace*{0.5cm}$\{$\\
\hspace*{1cm}for (int j = 0; j $<$ columnB; j++)\\
\hspace*{1cm}$\{$\\
\hspace*{1cm}if (dataGridB2.Rows[i].Cells[j].Value.ToString() == "E")\\
\hspace*{1.5cm}B2[i, j] = Int32.MaxValue;\\
\hspace*{1cm}else\\
\hspace*{1.5cm}B2[i, j] = Convert.ToInt16(dataGridB2.Rows[i].Cells[j].Value.ToString());\\
\hspace*{1cm}$\}$\\
\hspace*{0.5cm}$\}$\\
\hspace*{0.5cm}for (int i = 0; i $<$ lineA; i++)
\hspace*{0.5cm}$\{$\\
\hspace*{1cm}for (int j = 0; j $<$ columnB; j++)\\
\hspace*{1cm}$\{$\\
\hspace*{1.5cm}Product[i, j] = Int32.MaxValue;\\
\hspace*{1.5cm}for (k = 0; k $<$ lineB; k++)\\
\hspace*{1.5cm}$\{$\\
\hspace*{1.8cm} if (A2[i, k] == Int32.MaxValue $||$ B2[k, j] == Int32.MaxValue)\\
\hspace*{1.8cm}Sum[i, j] = Int32.MaxValue;\\
\hspace*{1.5cm} else\\
\hspace*{1.8cm}Sum[i, j] = A2[i, k] + B2[k, j];\\[0.1cm]
\hspace*{1.5cm}if (Sum[i, j] $>$ Product[i, j])\\
\hspace*{1.8cm}Product[i, j] = Product[i, j];\\
\hspace*{1.5cm}else\\
\hspace*{1.8cm}Product[i, j] = Sum[i, j];\\
\hspace*{1.5cm}$\}$\\
\hspace*{1cm}$\}$\\
\hspace*{0.5cm}$\}$\\
\hspace*{0.5cm}for (int i = 0; i $<$ lineA; i++)\\
\hspace*{0.5cm}$\{$\\
\hspace*{1cm}for (int j = 0; j $<$ columnB; j++)\\
\hspace*{1cm}$\{$\\
\hspace*{1.5cm} if (Product[i, j] == Int32.MaxValue)\\
\hspace*{1.8cm} dataGridProduct.Rows[i].Cells[j].Value = "E";\\
\hspace*{1.5cm}else\\
\hspace*{1.8cm}dataGridProduct.Rows[i].Cells[j].Value = Product[i, j];\\
\hspace*{1.5cm}$\}$\\
\hspace*{1cm}$\}$\\[0.1cm]
\hspace*{0.3cm}$\}$\\
\hspace*{0.3cm}private void textColumnA2$_{-}$Leave(object sender, EventArgs e)\\
\hspace*{0.3cm}$\{$\\
\hspace*{0.5cm}textLineB2.Text = textColumnA2.Text;\\
\hspace*{0.3cm}$\}$\\
\hspace*{0.3cm} private void btScalarProduct$_{-}$Click(object sender, EventArgs e)\\
\hspace*{0.3cm}$\{$\\
\hspace*{0.5cm}initMatrixScalarProduct ();\\
\hspace*{0.5cm}if (textScalar.Text != "")\\
\hspace*{0.5cm}$\{$\\
\hspace*{1cm}int lineA = 0;\\
\hspace*{1cm}int columnA = 0;\\
\hspace*{1cm}columnA = Convert.ToInt16(textcolumnA3.Text);\\
\hspace*{1cm}lineA = Convert.ToInt16(textLineA3.Text);\\
\hspace*{1cm}int[,] A3 = new int[lineA, columnA];\\
\hspace*{1cm}int[,] ScalarProduct = new int[lineA, columnA];\\
\hspace*{1cm}int a;\\
\hspace*{1cm}a = Convert.ToInt16(textScalar.Text);\\
\hspace*{1cm}for (int i = 0; i $<$ lineA; i++)\\
\hspace*{1cm}$\{$\\
\hspace*{1.5cm}for (int j = 0; j $<$ columnA; j++)\\
\hspace*{1.5cm}$\{$\\
\hspace*{1.8cm} if (dataGridA3.Rows[i].Cells[j].Value.ToString() == "E")\\
\hspace*{2cm}  A3[i, j] = Int32.MaxValue;\\
\hspace*{1.8cm}else\\
\hspace*{2cm}A3[i, j] =Convert.ToInt16(dataGridA3.Rows[i].Cells[j].Value.ToString());\\
\hspace*{1.5cm}$\}$\\
\hspace*{1cm}$\}$\\
\hspace*{1cm}for (int i = 0; i $<$ lineA; i++)\\
\hspace*{1cm}$\{$\\
\hspace*{1.5cm}for (int j = 0; j $<$ columnA; j++)\\
\hspace*{1.5cm}$\{$\\
\hspace*{1.8cm}if (A3[i, j] == Int32.MaxValue)\\
\hspace*{2cm}Scalar Product [i, j] = Int32.MaxValue;\\
\hspace*{1.8cm}else\\
\hspace*{2cm}Scalar Product [i, j] = (A3[i, j] + a);\\
\hspace*{1.5cm}$\}$\\
\hspace*{1cm}$\}$\\
\hspace*{1cm}for (int i = 0; i $<$ lineA; i++)\\
 \hspace*{1cm}$\{$\\
 \hspace*{1.5cm}for (int j = 0; j $<$ columnA; j++)\\
 \hspace*{1.5cm}$\{$\\
 \hspace*{1.8cm}if (ScalarProduct [i, j] == Int32.MaxValue)\\
 \hspace*{2cm}dataGridScalarProduct.Rows[i].Cells[j].Value = "E";\\
 \hspace*{1.8cm}else\\
 \hspace*{2cm}dataGridScalarProduct.Rows[i].Cells[j].Value = ScalarProduct[i,j];\\
 \hspace*{1.5cm}$\}$\\
 \hspace*{1cm}$\}$\\[0.2cm]
\hspace*{0.5cm}$\}$\\[0.1cm]
\hspace*{0.3cm}$\}$\\
\hspace*{0.3cm}$\}$private void btComputPower$_{-}$Click(object sender, EventArgs e)\\
\hspace*{0.3cm}$\{$\\[0.1cm]
\hspace*{0.5cm}initPowerMatrix ();\\
\hspace*{0.5cm}labelPower.Text = textPower.Text;\\
\hspace*{0.5cm}if (textPower.Text != "")\\
\hspace*{0.5cm}$\{$\\
\hspace*{1cm}int lineA = 0;\\
\hspace*{1cm}int columnA = 0;\\
\hspace*{1cm}columnA =Convert.ToInt16(textlineA4.Text);\\
\hspace*{1cm}lineA = Convert.ToInt16(textlineA4.Text);\\
\hspace*{1cm}int[,] A4 = new int[lineA, columnA];\\
\hspace*{1cm}int[,] B = new int[lineA, columnA];\\
\hspace*{1cm}int[,] Power$_{-}$n = new int[lineA, columnA];\\
\hspace*{1cm}int[,] sum = new int[lineA, columnA];\\
\hspace*{1cm}int a;\\
\hspace*{1cm}a = Convert.ToInt16(textPower.Text);\\
\hspace*{1cm}for (int i = 0; i $<$ lineA; i++)\\
\hspace*{1cm}$\{$\\
\hspace*{1.5cm}for (int j = 0; j $<$ columnA;j++)\\
\hspace*{1.5cm}$\{$\\
\hspace*{1.5cm}if (dataGridA4.Rows[i].Cells[j].Value.ToString() == "E")\\
\hspace*{1.8cm}  A4[i, j] =Int32.MaxValue;\\
\hspace*{1.5cm}else\\
\hspace*{1.8cm}A4[i, j] =
Convert.ToInt16(dataGridA4.Rows[i].Cells[j].Value.ToString());\\
\hspace*{1.5cm}B[i, j] = A4[i, j];\\
\hspace*{1.5cm}$\}$\\
\hspace*{1cm}$\}$\\
\hspace*{1cm}for (int p = 2; p $<=$ a; p++)\\
\hspace*{1cm}$\{$\\
\hspace*{1.5cm}for (int i = 0; i $<$ lineA;i++)\\
\hspace*{1.5cm}$\{$\\
\hspace*{1.5cm} for (int j = 0; j $<$ columnA;j++)\\
\hspace*{1.5cm}$\{$\\
\hspace*{1.8cm}Power$_{-}$n[i, j] = Int32.MaxValue;\\
\hspace*{1.8cm}for (int k = 0; k $<$ lineA; k++)\\
\hspace*{1.8cm}$\{$\\
\hspace*{2.3cm} //if (Power$_{-}$n[i, j] > (Power[i, k] + A4[k,j]))\\
\hspace*{2.3cm}// Power$_{-}$n[i, j] = Power[i, j];\\
\hspace*{2.3cm}//else\\
\hspace*{1.8cm}if (A4[i, k] == Int32.MaxValue $||$ B[k, j] == Int32.MaxValue)\\
\hspace*{2.1cm}sum[i, j] =Int32.MaxValue;\\
\hspace*{1.8cm}else\\
\hspace*{2.3cm}sum[i, j] = A4[i, k] + B[k,j];\\[0.1cm]
\hspace*{1.8cm}if (Power$_{-}$n[i, j] $<$ sum[i,j])\\
\hspace*{2.3cm}Power$_{-}$n[i, j] = Power$_{-}$n[i,j];\\
\hspace*{1.8cm}else\\
\hspace*{2.3cm}Power$_{-}n$[i, j] = sum[i,j];\\[0.1cm]
\hspace*{2.1cm}//sum[i, j] = Power$_{-}$n[i, j];\\
\hspace*{2.1cm}//Power$_{-}$n[i, j] = sum[i, j] + (A4[i,k]$\ast$B[k,j]);\\
\hspace*{1.8cm}$\}$\\
\hspace*{1.9cm}//Power$_{-}$n[i, j] = Power[i,j];\\
\hspace*{1.9cm}//Power = new int [i,j];\\
\hspace*{1.9cm}//Power[i, j] = Power$_{-}$n[i, j];\\
\hspace*{1.9cm}//dataGridScalarProduct.Rows[i].Cells[j].Value = Power$_{-}$n[i, j];\\
\hspace*{1.5cm}$\}$\\
\hspace*{1cm}$\}$\\
\hspace*{1cm}for (int i = 0; i $<$ lineA;i++)\\
\hspace*{1cm}$\{$\\
\hspace*{1.5cm}for (int j = 0; j $<$ columnA;j++)\\
\hspace*{1.5cm}$\{$\\
\hspace*{1.5cm}B[i, j] = Power$_{-}$n[i,j];\\
\hspace*{1.5cm}//dataGridMatrix$_{-}$at$_{-}$power$_{-}$n.Rows[i].Cells[j].Value = Power$_{-}$n[i, j];\\
\hspace*{1cm}$\}$\\
 \hspace*{1cm}$\}$\\
\hspace*{1.2cm}//Power$_{-}$n = B;\\
\hspace*{1cm}$\}$\\
 \hspace*{1cm}for (int i = 0; i $<$ lineA; i++)\\
 \hspace*{1cm}$\{$\\
 \hspace*{1.5cm}for (int j = 0; j $<$ columnA;j++)\\
 \hspace*{1.5cm}$\{$\\
 \hspace*{1.8cm}if (Power$_{-}$n[i, j] == Int32.MaxValue)\\
 \hspace*{2.1cm}dataGridMatrix$_{-}$at$_{-}$power$_{-}$n.Rows[i].Cells[j].Value ="E";\\
 \hspace*{1.8cm}else\\
 \hspace*{2.1cm}dataGridMatrix$_{-}$at$_{-}$Power$_{-}$n.Rows[i].Cells[j].Value = Power$_{-}$n[i,j];\\
 \hspace*{1.5cm}$\}$\\
 \hspace*{1cm}$\}$\\
 \hspace*{0.5cm}$\}$\\[0.2cm]
\hspace*{0.3cm}$\}$\\
\hspace*{0.3cm}private void ComputingDeterminant()\\
\hspace*{0.3cm}$\{$\\[0.1cm]
\hspace*{0.5cm}int lineA = 0; int columnA = 0;\\
\hspace*{0.5cm} lineA = Convert.ToInt16(textlinieA5.Text);\\
\hspace*{0.5cm}columnA = Convert.ToInt16(textlineA5.Text);\\
\hspace*{0.5cm}A5 = new int[lineA, lineA];\\
\hspace*{0.5cm}for (int i = 0; i < lineA; i++)\\
\hspace*{0.5cm}$\{$\\
\hspace*{1cm} for (int j = 0; j < columnA; j++)\\
\hspace*{1cm}$\{$\\
\hspace*{1.5cm}if (dataGridA5.Rows[i].Cells[j].Value.ToString() == "E")\\
\hspace*{1.8cm}A5[i, j] = Int32.MaxValue;\\
\hspace*{1.5cm}else\\
\hspace*{1.8cm} A5[i, j] = Convert.ToInt16(dataGridA5.Rows[i].Cells[j].Value.ToString());\\
\hspace*{1cm}$\}$\\
\hspace*{0.5cm}$\}$\\
\hspace*{0.5cm}int det$_{-}$plus;\\
\hspace*{0.5cm}int det$_{-}$minus;\\
\hspace*{0.5cm}int n = Convert.ToInt16(textlineA5.Text);\\
\hspace*{0.5cm}if (n == 1)\\
\hspace*{0.5cm}$\{$\\
\hspace*{1cm}det$_{-}$plus = A5[0, 0];\\
\hspace*{1cm}textDet$_{-}$plus.Text = det$_{-}$plus.ToString();\\
\hspace*{1cm}textDet$_{-}$minus.Text = "";\\
\hspace*{0.5cm}$\}$\\[0.1cm]
\hspace*{0.5cm}else if (n == 2)\\
\hspace*{0.5cm}{\\
\hspace*{1cm}if (A5[0, 0] == Int32.MaxValue $||$ A5[1, 1] == Int32.MaxValue)\\
\hspace*{1.5cm} det$_{-}$plus = Int32.MaxValue;\\
\hspace*{1cm}else\\
\hspace*{1.5cm}det$_{-}$plus = A5[0, 0] + A5[1,1];\\
\hspace*{1cm}if (A5[0, 1] == Int32.MaxValue $||$ A5[1, 0] == Int32.MaxValue)\\
\hspace*{1.5cm}det$_{-}$minus = Int32.MaxValue;\\
\hspace*{1cm}else\\
\hspace*{1.5cm}det$_{-}$minus = A5[0, 1] + A5[1,0];\\[0.1cm]
\hspace*{1cm}if (det$_{-}$plus == Int32.MaxValue)\\
\hspace*{1.5cm}textDet$_{-}$plus.Text = "E";\\
\hspace*{1cm}else\\
\hspace*{1.5cm}textDet$_{-}$plus.Text = det$_{-}$plus.ToString();\\
\hspace*{1cm}if (det$_{-}$minus == Int32.MaxValue)\\
\hspace*{1.5cm}textDet$_{-}$minus.Text = "E";\\
\hspace*{1cm}else\\
\hspace*{1.5cm}textDet$_{-}$minus.Text = det$_{-}$minus.ToString();\\
\hspace*{0.5cm}$\}$\\[0.1cm]
\hspace*{0.5cm}else if (n == 3)\\
\hspace*{0.5cm}{\\
\hspace*{1cm}int a1;\\
\hspace*{1cm}int a2;\\
\hspace*{1cm}int a3;\\
\hspace*{1cm}if (A5[0, 0] == Int32.MaxValue $||$ A5[1, 1] == Int32.MaxValue $||$ A5[2, 2] == Int32.MaxValue)\\
\hspace*{1.5cm}a1 = Int32.MaxValue;\\
\hspace*{1cm}else\\
\hspace*{1.5cm}a1 = A5[0, 0] + A5[1, 1] + A5[2,2];\\[0.1cm]
\hspace*{1cm}if (A5[1, 0] == Int32.MaxValue $||$ A5[2, 1] == Int32.MaxValue $||$ A5[0, 2] == Int32.MaxValue)\\
\hspace*{1.5cm}a2 = Int32.MaxValue;\\
\hspace*{1cm}else\\
\hspace*{1.5cm}a2 = A5[1, 0] + A5[2, 1] + A5[0,2];\\[0.1cm]
\hspace*{1cm}if (A5[2, 0] == Int32.MaxValue $||$ A5[0, 1] == Int32.MaxValue $||$ A5[1, 2] == Int32.MaxValue)\\
\hspace*{1.5cm}a3 = Int32.MaxValue;\\
\hspace*{1cm}else\\
\hspace*{1.5cm}a3 = A5[2, 0] + A5[0, 1] + A5[1,2];\\[0.1cm]
\hspace*{1cm}det$_{-}$plus = Int32.MaxValue;\\
\hspace*{1cm}if (a1 $<$ det$_{-}$plus)\\
\hspace*{1.5cm}det$_{-}$plus = a1;\\
\hspace*{1cm}if (a2 $<$ det$_{-}$plus)\\
\hspace*{1.5cm}det$_{-}$plus = a2;\\
\hspace*{1cm}if (a3 $<$ det$_{-}$plus)\\
\hspace*{1.5cm}det$_{-}$plus = a3;\\[0.1cm]
\hspace*{1cm}int b1;\\
\hspace*{1cm}int b2;\\
\hspace*{1cm}int b3;\\[0.1cm]
\hspace*{1cm}if (A5[0, 2] == Int32.MaxValue $||$ A5[1, 1] == Int32.MaxValue $||$ A5[2, 0] == Int32.MaxValue)\\
\hspace*{1.5cm}b1 = Int32.MaxValue;\\
\hspace*{1cm}else\\
\hspace*{1.5cm}b1 = A5[0, 2] + A5[1, 1] + A5[2,0];\\
\hspace*{1cm}if (A5[1, 2] == Int32.MaxValue $||$ A5[2, 1] == Int32.MaxValue $||$ A5[0, 0] == Int32.MaxValue)\\
\hspace*{1.5cm}b2 = Int32.MaxValue;\\
\hspace*{1cm}else\\
\hspace*{1.5cm}b2 = A5[1, 2] + A5[2, 1] + A5[0, 0];\\[0.1cm]
\hspace*{1cm}if (A5[0, 1] == Int32.MaxValue $||$ A5[1, 0] == Int32.MaxValue $||$ A5[2, 2] == Int32.MaxValue)\\
\hspace*{1.5cm}b3 = Int32.MaxValue;\\
\hspace*{1cm}else\\
\hspace*{1.5cm}b3 = A5[0, 1] + A5[1, 0] + A5[2, 2];\\[0.1cm]
\hspace*{1cm}det$_{-}$minus = Int32.MaxValue;\\[0.1cm]
\hspace*{1cm}if (b1 $<$ det$_{-}$minus)\\
\hspace*{1.5cm}det$_{-}$minus = b1;\\
\hspace*{1cm}if (b2 $<$ det$_{-}$minus)\\
\hspace*{1.5cm}det$_{-}$minus = b2;\\
\hspace*{1cm}if (b3 $<$ det$_{-}$minus)\\
\hspace*{1.5cm}det$_{-}$minus = b3;\\[0.1cm]
\hspace*{1cm}if (det$_{-}$plus == Int32.MaxValue)\\
\hspace*{1.5cm}textDet$_{-}$plus.Text = "E";\\
\hspace*{1cm}else\\
\hspace*{1.5cm}textDet$_{-}$plus.Text = det$_{-}$plus.ToString();\\
\hspace*{1cm}if (det$_{-}$minus == Int32.MaxValue)\\
\hspace*{1.5cm}textDet$_{-}$minus.Text = "E";\\
\hspace*{1cm}else\\
\hspace*{1.5cm}textDet$_{-}$minus.Text = det$_{-}$minus.ToString();\\
\hspace*{0.5cm}$\}$\\
\hspace*{0.3cm}$\}$\\
\hspace*{0.3cm}private void btComputingDeterminant$_{-}$Click(object sender, EventArgs e)\\
\hspace*{0.3cm}$\{$\\
\hspace*{0.5cm}ComputingDeterminant();\\
\hspace*{0.5cm}$\}$\\[0.2cm]
\hspace*{0.3cm}$\}$\\
\hspace*{0.3cm}private void btComputSystem$_{-}$Click(object sender, EventArgs e)\\
\hspace*{0.3cm}$\{$\\[0.1cm]
\hspace*{0.5cm} initMatrixSolutionSystem();\\
 \hspace*{0.5cm}label$_{-}$k.Text = textk.Text;\\
 \hspace*{0.5cm}if (textk.Text != "")\\
 \hspace*{0.5cm}$\{$\\
 \hspace*{1cm}int lineA = 0;\\
 \hspace*{1cm}int columnA = 0;\\
 \hspace*{1cm}int linex0 = 0;\\
 \hspace*{1cm} columnA = Convert.ToInt16(textlineA6.Text);\\
 \hspace*{1cm}lineA = Convert.ToInt16(textlineA6.Text);\\
 \hspace*{1cm}linex0 = Convert.ToInt16(textlineA6.Text);\\
 \hspace*{1cm}\hspace*{1cm}int[,] A6 = new int[lineA,columnA];\\
 \hspace*{1cm}nt[,] X0 = new int[linex0, 1];\\
 \hspace*{1cm}int[,] Xk = new int[lineA, 1];\\
 \hspace*{1cm}int[,] B = new int[lineA, columnA];\\
 \hspace*{1cm}int[,] Power$_{-}$n = new int[lineA, columnA];\\
 \hspace*{1cm}int[,] sum = new int[lineA, columnA];\\
 \hspace*{1cm}int[,] sum2 = new int[lineA, 1];\\
 \hspace*{1cm}int k;\\
 \hspace*{1cm}k = Convert.ToInt16(textk.Text);\\[0.1cm]
 \hspace*{0.5cm}for (int i = 0; i $<$ lineA; i++)\\
 \hspace*{0.5cm}$\{$\\
 \hspace*{1cm}for (int j = 0; j $<$ columnA; j++)\\
 \hspace*{1cm}$\{$\\
 \hspace*{1.5cm}if (dataGridA6.Rows[i].Cells[j].Value.ToString() == "E")\\
 \hspace*{1.8cm}A6[i, j] =Int32.MaxValue;\\
 \hspace*{1.5cm}else\\
 \hspace*{1.8cm}A6[i, j] = Convert.ToInt16(dataGridA6.Rows[i].Cells[j].Value.ToString());\\
  \hspace*{1.5cm}B[i, j] = A6[i, j];\\
  \hspace*{1cm}$\}$\\
  \hspace*{0.5cm}$\}$
  \hspace*{0.5cm}for (int i = 0; i $<$ linex0; i++)\\
  \hspace*{0.5cm}$\{$\\
\hspace*{1cm}if (dataGridX0.Rows[i].Cells[0].Value.ToString() == "E")\\
\hspace*{1.5cm}X0[i, 0] = Int32.MaxValue;\\
\hspace*{1cm}else\\
\hspace*{1.5cm}X0[i, 0] = Convert.ToInt16(dataGridX0.Rows[i].Cells[0].Value.ToString());\\
 \hspace*{0.5cm}$\}$\\
 \hspace*{0.5cm}if (k == 1)\\
 \hspace*{0.5cm}$\{$\\
 \hspace*{1cm}for (int i = 0; i $<$ lineA; i++)\\
 \hspace*{1cm}$\{$\\
 \hspace*{1.5cm}for (int j = 0; j $<$ columnA; j++)\\
 \hspace*{1.5cm}$\{$\\
 \hspace*{1.8cm}Power$_{-}$n[i, j] = A6[i,j];\\
 \hspace*{1.5cm}$\}$\\
 \hspace*{1cm}$\}$\\
 \hspace*{0.5cm}$\}$\\
 \hspace*{0.5cm}else\\
 \hspace*{0.5cm}$\{$\\
  \hspace*{1cm}for (int p = 2; p $<=$ k; p++)\\
  \hspace*{1cm}$\{$\\
 \hspace*{1.5cm}for (int i = 0; i $<$ lineA; i++)\\
 \hspace*{1.5cm}$\{$\\
 \hspace*{1.5cm}for (int j = 0; j $<$ columnA; j++)\\
 \hspace*{1.5cm}$\{$\\
\hspace*{1.8cm} Power$_{-}$n[i, j] = Int32.MaxValue;\\
\hspace*{1.8cm}for (int h = 0; h $<$ lineA; h++)\\
\hspace*{1.8cm}$\{$\\
\hspace*{1.8cm}if (A6[i, h] == Int32.MaxValue $||$ B[h, j] == Int32.MaxValue)\\
\hspace*{2cm}sum[i, j] = Int32.MaxValue;\\
\hspace*{1.8cm}else\\
\hspace*{2cm}sum[i, j] = A6[i, h] + B[h,j];\\
\hspace*{1.8cm} if (Power$_{-}$n[i, j] $<$ sum[i,j])\\
\hspace*{2cm}  Power$_{-}$n[i, j] = Power$_{-}$n[i,j];\\
\hspace*{1.8cm}else\\
\hspace*{2cm}Power$_{-}$n[i, j] = sum[i,j];\\
\hspace*{1.5cm}$\}$\\
\hspace*{1cm}$\}$\\
\hspace*{0.5cm}$\}$\\
\hspace*{1cm}for (int i = 0; i $<$ lineA; i++)\\
\hspace*{1.5cm}$\{$\\
\hspace*{1.8cm}for (int j = 0; j $<$ columnA;j++)\\
\hspace*{1.5cm}$\{$\\
\hspace*{1.8cm}B[i, j] = Power$_{-}$n[i, j];\\
\hspace*{1.5cm}$\}$\\
\hspace*{1cm}$\}$\\
\hspace*{0.5cm}$\}$\\
$\}$\\
for (int i = 0; i $<$ lineA; i++)\\
$\{$\\
\hspace*{0.5cm}for (int j = 0; j $<$ 1; j++)\\
\hspace*{0.5cm}$\{$\\
\hspace*{0.5cm} Xk[i, j] = Int32.MaxValue;\\
\hspace*{0.5cm}for (int h = 0; h $<$ lineA; h++)\\
\hspace*{0.5cm}$\{$\\
\hspace*{1cm}if (Power$_{-}$n[i, h] == Int32.MaxValue $||$ X0[h,j] == Int32.MaxValue)\\
 \hspace*{1.5cm}sum2[i, j] = Int32.MaxValue;\\
 \hspace*{1cm}else\\
 \hspace*{1.5cm}sum2[i, j] = Power$_{-}$n[i, h] + X0[h, j];\\[0.1cm]
 \hspace*{1.5cm} if (Xk[i, j] $<$ sum2[i,j])\\
 \hspace*{1.5cm}Xk[i, j] = Xk[i,j];\\
 \hspace*{1cm}else\\
 \hspace*{1.5cm}Xk[i, j] = sum2[i,j];\\
 \hspace*{1.5cm}$\}$\\
 \hspace*{1cm}$\}$\\
 \hspace*{0.5cm}$\}$\\
 \hspace*{0.5cm} for (int i = 0; i $<$ lineA; i++)\\
 \hspace*{0.5cm}$\{$\\
 \hspace*{1cm}if (Xk[i, 0] == Int32.MaxValue)\\
 \hspace*{1.5cm}dataGridSolutionXk.Rows[i].Cells[0].Value = "E";\\
 \hspace*{1cm}else\\
 \hspace*{1.5cm}dataGridSolutionXk.Rows[i].Cells[0].Value = Xk[i,0];\\
\hspace*{0.5cm}$\}$\\
\hspace*{0.3cm}$\}$\\
$\}$\\
public void Reset$_{-}$ Values $_{-}$ For $_{-}$ Addition()\\
$\{$\\
 \hspace*{0.5cm}textColumn.ResetText();\\
 \hspace*{0.5cm}textLine.ResetText();\\
 \hspace*{0.5cm}dataGridA.Rows.Clear();\\
 \hspace*{0.5cm}dataGridB.Rows.Clear();\\
 \hspace*{0.5cm}dataGridResAddition.Rows.Clear();\\[0.1cm]
 $\}$\\
public void Reset$_{-}$Values for Multiplication()\\
$\{$\\
\hspace*{0.5cm}textColumnA2.ResetText();\\
\hspace*{0.5cm}textLineA2.ResetText();\\
\hspace*{0.5cm}textcolumnB2.ResetText();\\
\hspace*{0.5cm}textLineB2.ResetText();\\
\hspace*{0.5cm}dataGridA2.Rows.Clear();\\
\hspace*{0.5cm}dataGridB2.Rows.Clear();\\
\hspace*{0.5cm}dataGridProduct.Rows.Clear();\\[0.1cm]
$\}$\\
public void Reset$_{-}$Values$_{-}$for$_{-}$Scalar Multiplication()\\
$\{$\\
\hspace*{0.5cm}textcolumnA3.ResetText();\\
\hspace*{0.5cm}textScalar.ResetText();\\
\hspace*{0.5cm}textLineA3.ResetText();\\
\hspace*{0.5cm}dataGridA3.Rows.Clear();\\
\hspace*{0.5cm}dataGridScalarProduct.Rows.Clear();\\[0.1cm]
$\}$\\
public void Reset$_{-}$Values$_{-}$for$_{-}$Lifting$_{-}$at$_{-}$Power()\\
 $\{$\\
 \hspace*{0.5cm}textlineA4.ResetText();\\
 \hspace*{0.5cm}textPower.ResetText();\\
 \hspace*{0.5cm}dataGridA4.Rows.Clear();\\
 \hspace*{0.5cm}dataGridMatrix$_{-}$at$_{-}$Power$_{-}$n.Rows.Clear();\\[0.1cm]
 $\}$\\
public void Reset$_{-}$Values$_{-}$for$_{-}$Computing$_{-}$Determinant()\\
 $\{$\\
\hspace*{0.5cm}textlineA5.ResetText();\\
\hspace*{0.5cm}Det$_{-}$minus.ResetText();\\
\hspace*{0.5cm}Det$_{-}$plus.ResetText();\\
 \hspace*{0.5cm}dataGridA5.Rows.Clear();\\[0.1cm]
 $\}$\\
public void Reset$_{-}$Values$_{-}$for$_{-}$Computation$_{-}$system()\\
$\{$\\
\hspace*{0.5cm}textlineA6.ResetText();\\
\hspace*{0.5cm}textk.ResetText();\\
\hspace*{0.5cm}label$_{-}$k.Text = "k";\\
\hspace*{0.5cm}dataGridSolutionXk.Rows.Clear();\\
\hspace*{0.5cm}dataGridA6.Rows.Clear();\\
\hspace*{0.5cm}dataGridX0.Rows.Clear();\\[0.2cm]
\hspace*{0.5cm}$\}$\\
\hspace*{0.5cm}private void btReset$_{-}$Click(object sender, EventArgs e)\\
 \hspace*{0.5cm}$\{$\\
 \hspace*{1cm}Reset$_{-}$Values $_{-}$For$_{-}$Addition();\\
 \hspace*{0.5cm}$\}$\\
\hspace*{0.5cm}private void btResetMultiplication$_{-}$Click(object sender, EventArgs e)\\
\hspace*{0.5cm}$\{$\\
\hspace*{1cm}Reset$_{-}$Values$_{-}$for$_{-}$Multiplication ();\\
\hspace*{0.5cm}$\}$\\
\hspace*{0.5cm}private void btReset$_{-}$Scalar$_{-}$Product$_{-}$Click(object sender, EventArgs e)\\
\hspace*{0.5cm}$\{$\\
\hspace*{1cm}Reset$_{-}$Values$_{-}$for$_{-}$Scalar Multiplication();\\
\hspace*{0.5cm}$\}$\\
\hspace*{0.5cm}private void btReset$_{-}$lifting$_{-}$at$_{-}$power$_{-}$Click(object sender, EventArgs e)\\
\hspace*{0.5cm}$\{$\\
\hspace*{1cm}Reset$_{-}$Values$_{-}$for$_{-}$Lifting$_{-}$at$_{-}$Power();\\
\hspace*{0.5cm}$\}$\\
\hspace*{0.5cm}private void btReset$_{-}$values$_{-}$det$_{-}$Click(object sender, EventArgs e)\\
 \hspace*{0.5cm}$\{$\\
 \hspace*{1cm}Reset$_{-}$Values $_{-}$For$_{-}$Computing$_{-}$Determinant();\\
 \hspace*{0.5cm}$\}$\\
\hspace*{0.5cm}private void btResetSystem$_{-}$Click(object sender, EventArgs e)\\
\hspace*{0.5cm}$\{$\\
\hspace*{1cm}Reset$_{-}$Values$_{-}$for$_{-}$Computation$_{-}$System();\\
\hspace*{0.5cm}$\}$\\
\hspace*{0.3cm}$\}$\\
$\}$\\

 We illustrate the utilization of the above programs in the
 solving of the following problems.\\[0.2cm]
{\bf Problem 3.1. Sum of two matrices $A$ and $B$ over ${\bf R}_{min}$}.\\[0.3cm]
{\it Inputs data}\hspace*{8.5cm}{\it Outputs data}\\[0.2cm]
Number of lines $~~~~~~\fbox{3}$\\[0.1cm]
Number of columns $~\fbox{5}$\\[0.1cm]
Matrix $~A$\hspace*{2.5cm} Matrix $~B$\hspace*{4.1cm} Matrix
$~~A\oplus B$  \\[0.2cm]
 $\begin{array}{|c|c|c|c|c|}\hline
 12     & 0 & E & 3 & 5\cr \hline
  -1    & 0 & 2 & E & 7\cr \hline
5    & 3 & 1 & 8 & 0\cr \hline
\end{array} ~~~~~~~~ \begin{array}{|c|c|c|c|c|}\hline
 0     & -3 & 2 & E & E\cr \hline
  -8    & -1 & 3 & 8 & E\cr \hline
E    & 1 & -5 & 2 & 7\cr \hline
\end{array}~~~~~~~~~~~~~~~~\begin{array}{|c|c|c|c|c|}\hline
 0     & -3 & 2 & 3 & 5\cr \hline
  -8    & -1 & 2 & 8 & 7\cr \hline
5    & 1 & -5 & 2 & 0\cr \hline
\end{array}.$\\[0.4cm]
{\bf Problem 3.2. Product of two matrices $A$ and $B$ over ${\bf R}_{min}$}.\\[0.3cm]
{\it Inputs data}\hspace*{10.2cm}{\it Outputs data}\\[0.2cm]
Number of lines of $~A!~~~~~~~\fbox{5}~~~$ Number of lines of $~B!~~~~~~~\fbox{3}$\\[0.1cm]
Number of columns of $~A!~~\fbox{3}~~~$ Number of columns of $~B!~~\fbox{4}$\\[0.1cm]
Matrix $A$\hspace*{2.5cm} Matrix $B$ \hspace*{6.2cm} Matrix $A\otimes B$ \\[0.2cm]
 $\begin{array}{|c|c|c|}\hline
 12     &  -1 & 5\cr \hline
0 & 0 & 3\cr \hline
E & 2 & 1\cr \hline
 3   & E & 8\cr \hline
 5& 7 & 0\cr \hline
\end{array} ~~~~~~~~~~~~~ \begin{array}{|c|c|c|c|}\hline
 1     & 7 & E & 0 \cr \hline
  4    & -2 & 1 & E \cr \hline
3   & 7 & E & 0 \cr \hline
\end{array}~~~~~~~~~~~~~~~~~~~~~~~~~~~~~~~~~~~~~~~~~~~~\begin{array}{|c|c|c|c|}\hline
 3  & -3 & 0 & 5 \cr \hline
  1 & -2 & 1 & 0 \cr \hline
4 & 0 & 3 & 1 \cr \hline 4 & 10 & E & 3 \cr \hline
 3 & 5 & 8 & 0\cr \hline
\end{array}.$\\[0.4cm]
{\bf Problem 3.3. Multiplication with scalar $a\in {\bf R}$ of a matrix $A$ over ${\bf R}_{min}$}.\\[0.3cm]
{\it Inputs data}\hspace*{7.5cm}{\it Outputs data}\\[0.2cm]
Number of lines of $~A!~~~~~~\fbox{3}$\\[0.1cm]
Number of columns of $~A!~~\fbox{4}$\\[0.1cm]
Scalar $a~~~\fbox{-5}$\\[0.2cm]
Matrix $~A$\hspace*{7.5cm} Matrix $~~a\otimes
A$\\[0.2cm]
 $\begin{array}{|c|c|c|c|}\hline
 4     & -7 & 8 &  E\cr \hline
  5    & E & 0 &  8\cr \hline
9    & 2 & 3 & 1\cr \hline
\end{array},~~~~~~~~~~~~~~~~~~~~~~~~~~~~~~~~~~~~~~~~~~~~~~~~~~~\begin{array}{|c|c|c|c|c|}\hline
 -1     & -12 & 3 & E \cr \hline
  0    & E & -5 & 3\cr \hline
4    & -3 & -2 & -4 \cr \hline
\end{array}. $\\[0.4cm]
{\bf Problem 3.4. Compute the matrix $A^{n}$ for $~n\geq 2$, where $A\in M_{7\times 7}({\bf R}_{min})$}\\[0.3cm]
{\it Inputs data}\hspace*{9cm}{\it Outputs data}\\[0.2cm]
Number of lines and columns of $~A!~~~\fbox{7}$\\[0.1cm]
Power of $~A!~~~~~~~~~~\fbox{2}$\\[0.2cm]
Matrix $~A$\hspace*{9cm} Matrix $A^{(2)}$\\[0.2cm]
 $\begin{array}{|c|c|c|c|c|c|c|}\hline
 0 & 4  & E &  E & 7 & E & E\cr \hline
 4 & 0  & 5 & 9 & E & E & E\cr \hline
 E & 5 & 0 & E &  6 & 10 & E\cr \hline
E & 9 & E & 0 & E & 8 & E\cr \hline
 7 & E &  6 & E & 0 & 14 & E\cr \hline
E & E &  10 & 8 & 14 & 0 & 11\cr \hline
 E & E &  E & E & E & 11 & 0\cr \hline
\end{array},~~~~~~~~~~~~~~~~~~~~~~~~~~~~~~~~~~~~~~~~~~~\begin{array}{|c|c|c|c|c|c|c|}\hline
0 & 4  & 9 &  13 & 7 & 21 & E\cr \hline
 4 & 0  & 5 & 9 & 11 & 15 & E\cr \hline
 9 & 5 & 0 & 14 &  6 & 10 & 21\cr \hline
13 & 9 & 14 & 0 & 22 & 8 & 19\cr \hline
 7 & 11 &  6 & 22 & 0 & 14 & 25\cr \hline
21 & 15 &  10 & 8 & 14 & 0 & 11\cr \hline
 E & E &  21 & 19 & 25 & 11 & 0\cr \hline
 \end{array}. $\\

Similarly, we have\\[0.2cm]
$A^{(3)}=\left ( \begin{array}{ccccccc}
 0 & 4  & 9 &  13 & 7 & 19 & 32\cr
 4 & 0  & 5 & 9 & 11 & 15 & 26\cr
 9 & 5 & 0 & 14 &  6 & 10 & 21\cr
13 & 9 & 14 & 0 & 20 & 8 & 19\cr
 7 & 11 &  6 & 20 & 0 & 14 & 25\cr
19 & 15 &  10 & 8 & 14 & 0 & 11\cr
 32 & 26 &  21 & 19 & 25 & 11 & 0\cr
\end{array}\right ),~~~~~~~~A^{(4)}=\left ( \begin{array}{ccccccc}
 0 & 4  & 9 &  13 & 7 & 19 & 30\cr
 4 & 0  & 5 & 9 & 11 & 15 & 26\cr
 9 & 5 & 0 & 14 &  6 & 10 & 21\cr
13 & 9 & 14 & 0 & 20 & 8 & 19\cr
 7 & 11 &  6 & 20 & 0 & 14 & 25\cr
19 & 15 &  10 & 8 & 14 & 0 & 11\cr
 30 & 26 &  21 & 19 & 25 & 11 & 0\cr
 \end{array}\right )$\\[0.2cm]
and $A^{(n)}= A^{(4)},~$ for all $n\geq 5.$\\[0.4cm]
{\bf Problem 3.5. Compute the bideterminant and permanent of $A\in M_{3\times 3}({\bf R}_{min})$}.\\[0.3cm]
{\it Inputs data}\hspace*{7.5cm}{\it Outputs data}\\[0.2cm]
Matrix $~A$\\[0.2cm]
 $\begin{array}{|c|c|c|}\hline
 8 & -6 & -1\cr \hline
 3 & 5 & 1\cr \hline
 9 & -4 & 0 \cr \hline
\end{array}~~~~~~~~~~~~~~~~~~~~~~~~~~~~~~~~~~~~~~~~~~~~~~~~~~~~~~~\begin{array}{c}
\det_{\oplus}^{+}(A)=-2\cr
\det_{\oplus}^{-}(A)=-3\cr
\end{array}. $\\[0.4cm].
Then, the bideterminant of $A$ is $ (\Delta_{1}(A), \Delta_{2}(A))
= ( -2, -3)$ and its permanent is $perm(A)=(-2)\oplus (-3)
=\min\{-2, -3\} = -3.$ \\[0.5cm]
{\bf Problem 3.6. Solve the system $~X(k+1)=AX(k)~(k\geq 0)$ with $A\in M_{6\times 6}({\bf R}_{min})$}.\\[0.3cm]
{\it Inputs data}\hspace*{9.5cm}{\it Outputs data}\\[0.2cm]
Number of lines and columns of $~A!~~~~~\fbox{6}$\\[0.1cm]
Value for $~k!~~~\fbox{1}$\\[0.2cm]
Matrix $~A!~~~~~~~~~~~~~~~~~~~~~~~~~~~~~~~~~$ Matrix $X(0)~~~~~~~~~~~~~~~~~~~~~~~~~~$ Matrix $X(1)$\\[0.2cm]
$\begin{array}{|c|c|c|c|c|c|}\hline

0 & 10 & 12 & 2  & 14 &E\cr\hline

E & 0 & E & E  & 5 & E \cr \hline

 E & 4 & 0 & E  & E & E \cr \hline

E & E & E & 0  & E & -3 \cr \hline

E & E & E & E & 0 & 4 \cr \hline

5 & E & E & 6  & E & 0 \cr \hline
\end{array}~~~~~~~~~~~~~~~~~~ \begin{array}{|c|}\hline
 10 \cr \hline
  -2\cr \hline
3\cr \hline
5\cr \hline
9\cr \hline
4\cr \hline
\end{array}~~~~~~~~~~~~~~~~~~~~~~~~~~~~~~~~~~~~~\begin{array}{|c|}\hline
7 \cr \hline
  -2\cr \hline
2\cr \hline
1\cr \hline
8\cr \hline
4\cr \hline
\end{array}.$\\

Similarly,  $~X(2)=(3~~ -2~~~ 2~~~ 1~~~ 8~~~ 4)^{T}~$ and
$~X(k)=X(2)~$ for all $k\geq 3.$\\

Author's adresses\\

West University of Timi\c soara,\\
Seminar of Algebra. Department of Mathematics,\\
 Bd. V. P{\^a}rvan, no. 4, 300223, Timi\c soara, Romania\\
\hspace*{0.7cm} E-mail: ivangm31@yahoo.com; ivan@math.uvt.ro\\
\end{document}